\documentclass[11pt,reqno]{amsart}

\usepackage[top=1in, bottom=1in, left=1in, right=1in]{geometry}

\usepackage{latexsym}
\usepackage{amsmath}
\usepackage{amssymb}
\usepackage{mathrsfs}
\usepackage[T1]{fontenc}
\usepackage[latin1]{inputenc}
\usepackage{mathtools}
\usepackage{dsfont}
\usepackage[nocompress]{cite}
\usepackage{microtype}
\usepackage{comment}

\usepackage[shortlabels]{enumitem}
\usepackage{aliascnt}
\usepackage[bookmarks=true,pdfstartview=FitH, pdfborder={0 0 0}, colorlinks=true,citecolor=red, linkcolor=blue]{hyperref}

\theoremstyle{plain}
\newtheorem{thm}{Theorem}[section]
\newaliascnt{cor}{thm}
\newaliascnt{prop}{thm}
\newaliascnt{lem}{thm}
\newaliascnt{claim}{thm}
\newtheorem{cor}[cor]{Corollary}
\newtheorem{prop}[prop]{Proposition}
\newtheorem{lem}[lem]{Lemma}
\newtheorem{claim}[claim]{Claim}
\aliascntresetthe{cor}
\aliascntresetthe{prop}
\aliascntresetthe{lem}
\aliascntresetthe{claim}

\theoremstyle{definition}
\newaliascnt{defn}{thm}
\newaliascnt{asu}{thm}
\newaliascnt{con}{thm}
\newtheorem{defn}[defn]{Definition}

\aliascntresetthe{defn}
\aliascntresetthe{asu}
\aliascntresetthe{con}

\newcounter{stp}
\newcounter{stpi}
\newcounter{stpci}
\newcounter{stpiii}

\theoremstyle{remark}
\newaliascnt{rem}{thm}
\newaliascnt{exa}{thm}
\newaliascnt{masu}{thm}
\newaliascnt{nota}{thm}
\newaliascnt{sett}{thm}
\newtheorem{rem}[rem]{Remark}

\aliascntresetthe{rem}
\aliascntresetthe{exa}
\aliascntresetthe{masu}
\aliascntresetthe{nota}
\aliascntresetthe{sett}

\numberwithin{equation}{section}

\setlist[enumerate]{font = \normalfont}

\newcommand{\R}{\mathbb{R}}

\newcommand{\N}{\mathbb{N}}

\newcommand{\I}{\mathbb{I}}
\newcommand{\bS}{\mathbb{S}}

\newcommand{\dsone}{\mathds{1}}
\newcommand{\rW}{\mathrm{W}}
\newcommand{\rL}{\mathrm{L}}

\newcommand{\rH}{\mathrm{H}}

\newcommand{\rC}{\mathrm{C}}

\newcommand{\rD}{\mathrm{D}}

\newcommand{\SO}{\mathrm{SO}}
\newcommand{\rd}{\,\mathrm{d}}

\newcommand{\dist}{\mathrm{dist}}

\newcommand{\cS}{\mathcal{S}}
\newcommand{\cF}{\mathcal{F}}

\newcommand{\cQ}{\mathcal{Q}}
\newcommand{\cK}{\mathcal{K}}

\newcommand{\cT}{\mathcal{T}}

\newcommand{\cD}{\mathcal{D}}
\newcommand{\cH}{\mathcal{H}}
\newcommand{\cV}{\mathcal{V}}
\newcommand{\cP}{\mathcal{P}}

\newcommand{\cX}{\mathcal{X}}

\newcommand{\eps}{\varepsilon}
\newcommand{\del}{\partial}

\DeclareMathOperator{\mdiv}{div}
\DeclareMathOperator{\Char}{Char}
\DeclareMathOperator{\Span}{Span}
\DeclareMathOperator{\supp}{supp}

\newcommand{\tin}{\enspace \text{in} \enspace}
\newcommand{\ton}{\enspace \text{on} \enspace}
\newcommand{\tfor}{\enspace \text{for} \enspace}
\newcommand{\tforall}{\enspace \text{for all} \enspace}
\newcommand{\tand}{\enspace \text{and} \enspace}

\newcommand{\tas}{\enspace \text{as} \enspace}

\newcommand{\bu}{\bfu}

\newcommand{\bfx}{\boldsymbol{x}}
\newcommand{\bfy}{\boldsymbol{y}}
\newcommand{\bfh}{\boldsymbol{h}}
\newcommand{\bfw}{\boldsymbol{w}}
\newcommand{\bfq}{\boldsymbol{q}}
\newcommand{\bfr}{\boldsymbol{r}}
\newcommand{\bfz}{\boldsymbol{z}}
\newcommand{\bfY}{\boldsymbol{Y}}
\newcommand{\bfX}{\boldsymbol{X}}
\newcommand{\bfu}{\boldsymbol{u}}
\newcommand{\bfU}{\boldsymbol{U}}
\newcommand{\bfM}{\boldsymbol{M}}
\newcommand{\bfR}{\boldsymbol{R}}
\newcommand{\bfB}{\boldsymbol{B}}
\newcommand{\bfZ}{\boldsymbol{Z}}
\newcommand{\bfO}{\boldsymbol{O}}
\newcommand{\bfd}{\boldsymbol{d}}
\newcommand{\bfa}{\boldsymbol{a}}
\newcommand{\bfm}{\boldsymbol{m}}
\newcommand{\bb}{\boldsymbol{b}}
\newcommand{\bfn}{\boldsymbol{n}}
\newcommand{\bfomega}{\boldsymbol{\omega}}
\newcommand{\bff}{\boldsymbol{f}}
\newcommand{\bfell}{\boldsymbol{\ell}}
\newcommand{\bfv}{\boldsymbol{v}}
\newcommand{\bfpsi}{\boldsymbol{\psi}}

\newcommand{\bfp}{\boldsymbol{p}}
\newcommand{\bfzeta}{\boldsymbol{\zeta}}
\newcommand{\bftau}{\boldsymbol{\tau}}
\newcommand{\bfphi}{\boldsymbol{\phi}}
\newcommand{\bfe}{\boldsymbol{e}}

\newcommand{\bfxi}{\boldsymbol{\xi}}
\newcommand{\tbfd}{\Tilde{\bfd}}

\newcommand{\srd}{\rd}

\newcommand{\bJ}{\mathbb{J}}
\newcommand{\bD}{\mathbb{D}}
\newcommand{\bT}{\mathbb{T}}

\begin{document}

\title[Existence of global weak solutions for Ericksen--Leslie fluid--rigid body interaction]{Existence of global weak solutions for a simplified Ericksen--Leslie fluid interacting with a rigid body}

\author{Felix Brandt}
\address{Department of Mathematics, University of Southern California, Los Angeles, CA 90089, USA.}
\email{fbrandt@usc.edu}
\author{Vinh Nguyen}
\address{Department of Mathematics, University of California, Berkeley, CA 94720, USA.}
\email{vnguyen26@berkeley.edu}


\begin{abstract}
In this paper, we consider a three-dimensional nematic liquid crystal governed by the simplified Ericksen--Leslie system interacting with a moving rigid body. 
We prove the existence of global weak solutions, where we especially treat the case of contact with the container wall. 
The director field $\bfd$ satisfies $|\bfd| = 1$ as well as $d_3 \ge 0$, and its non-homogeneous Dirichlet boundary condition on the body interface, referred to as strong anchoring, is transported with the rigid motion. 
The moving-domain problem with the above director condition is obtained as a limit of a Ginzburg--Landau regularization.
An extension of the director to the body domain yields global compactness, while a compactness argument as in the seminal work by Lin and Wang \cite{LW:16} yields strong convergence away from the interface.
If the body is in contact with the wall, we impose a perfectly inelastic sticking law. 
Since the resulting domain after contact is not necessarily Lipschitz, we prove weak existence for the simplified Ericksen--Leslie model on arbitrary bounded open sets, which is of independent interest.
\end{abstract}
\maketitle

\section{Introduction}

The Ericksen--Leslie theory couples incompressible flow to a unit-vector director describing the local orientation of a nematic liquid crystal. 
We refer to Ericksen \cite{Ericksen:62} and Leslie \cite{Leslie:68}, the monographs \cite{BE:94,dGP:93}, and the surveys \cite{Ball:17,WZZ:21} for a comprehensive overview of liquid crystal models. 
The pointwise constraint $|\bfd|=1$ present in Ericksen--Leslie models is analytically delicate in three dimensions, because the natural energy bounds give only weak $\rH^1$-compactness, which does not identify the quadratic Ericksen stress and allows concentration of the director energy.

The interaction problem of an Ericksen--Leslie model with a rigid body body combines the Ericksen--Leslie description of a nematic flow with Newton's laws, giving rise to a moving boundary problem in which the fluid velocity, the molecular orientation, expressed by the director filed, and the rigid motion are coupled.
Mechanically, a significant difficulty arises due to the so-called anchoring at the surface of the rigid body.
It produces viscous forces and additionally distorts the director field, generating elastic forces.
Mathematically, this is expressed by a non-homogeneous Dirichlet boundary condition of the director field.
Another important aspect in this context of the possibility of contact of the rigid body with the outer boundary of the fluid container.

In this paper, we study this coupling for the simplified Ericksen--Leslie system introduced by Lin~\cite{Lin:89}. 
The fluid is subject to no-slip boundary conditions, while the director $\bfd = (d_1,d_2,d_3)$ satisfies the condition $|\bfd| = 1$ as well as $d_3 \ge 0$, and its boundary values on the body are fixed in material coordinate.
Our main result, \autoref{exist:main}, constructs a finite-energy weak solution up to the first wall contact and continues it afterwards under permanent sticking, i.e., the rigid body is assumed to remain attached to the outer fluid boundary after contact.
Moreover, in \autoref{cor:quantitative separation main}, we provide an explicit lower bound for the contact time. 

The solution to the interaction problem of the Ericksen-Leslie model including the aforementioned director condition is obtained by passing to the limit in a Ginzburg--Landau regularization.
On fixed domains, the corresponding regularized model was analyzed by Lin and Liu~\cite{LL:95}.
In their seminal work \cite{LW:16}, Lin and Wang proved global existence of a weak solution to the simplified Ericksen--Leslie model with director condition.
The essential input from~\cite{LW:16} for the present work is a strong compactness result for approximate harmonic maps in the interior of the fluid domain, see \autoref{lem:LW compactness input}.
An additional difficulty in this article is that we need to extend the director through the body.
This extension satisfies a transport equation on the whole fluid--rigid body container, which paves the way for global strong~$\rL^2$-compactness.

Since the domain resulting from a contact of the rigid body with the outer boundary may have cusps or several components, we prove in \autoref{thm:arbitrary open constrained} the existence of a weak solution to the simplified Ericksen--Leslie model with director condition on arbitrary open sets.
The proof uses the variational Dirichlet Laplacian, compactness of the zero-extension, and a compactness argument \`a la Lin--Wang in the interior, without imposing regularity on the boundary.

For the analysis of Ericksen--Leslie models on fixed domains, we refer to \cite{FHM:20} for smooth relaxation limits, \cite{KW:25} for weak solutions in the axisymmetric setting, \cite{WXL:13,CRRS:17,HNPS:16,HLW:25} for strong solutions, or \cite{HLLW:16,Zuo:24} for singularity formation and regularity.
A different description of liquid crystals is due to Beris and Edwards \cite{BE:94}, where the molecular orientation is described by tensors that are symmetric and have trace zero.
The mathematical analysis of such $Q$-tensor models on fixed domains has been carried out in \cite{PZ:11,PZ:12,ADL:14,Wil:15,FRSZ:14,MS:22,HHW:24}.

In the context of interaction problems of viscous incompressible fluids with a rigid body, we refer, e.g., to \cite{DE:99,DE:00,CSMT:00,GLS:00,MST:02,F3,FHN:08} for existence results of weak solutions, while strong well-posedness and the long-time behavior of such fluid--rigid interaction problems has, for example, been discussed in \cite{GGH:13,EMT:23}.
Recent work concerns regularity of weak solutions \cite{MNR:24}, small-body and vanishing-body limits \cite{BN:24,HeSu:25,FRZ:23}, and multi-particle or long-time regimes \cite{bravin2025long}.
The behavior near contact is especially sensitive to the boundary condition.
Under no slip, non-collision and blow-up mechanisms were studied in \cite{HT:09, HT:10}, while the case with slip was analyzed in \cite{GVH:14, GHVW:15}. 

For liquid crystals coupled to a rigid body, Binz et al.~\cite{BBHR:24} proved local strong well-posedness and global strong existence near equilibria for the simplified Ericksen--Leslie system with Neumann director conditions. 
Brandt et al.~\cite{BHR:26} proved similar results in the context of the general $Q$-tensor model.
Geng et al.~\cite{GRZ:24} constructed global weak solutions for a rigid body in a $Q$-tensor fluid, including the discussion of wall contact. 
Mensah~\cite{Men:26} considered a two-dimensional simplified Ericksen--Leslie fluid coupled to a flexible shell. 

It is worthwhile comparing our result to the previous works \cite{BBHR:24, GRZ:24}.
In fact, we allow for general data in the natural finite-energy class, so in particular, compared to \cite{BBHR:24}, we consider lower regularity for the initial data, and we do not require a smallness assumption to infer the global-in-time existence of a solution.
In addition, we discuss the case of contact of the rigid body with the fluid boundary, which has been omitted in the work \cite{BBHR:24} on strong solutions.
Let us emphasize that we address the simplified Ericksen--Leslie model including the director condition $|\bfd| = 1$.
This condition as well as the aforementioned anchoring of the rigid body, which is expressed via non-homogeneous Dirichlet boundary conditions, and which requires a delicate extension procedure, pose significant additional challenges compared to the work \cite{GRZ:24} in the context of the $Q$-tensor model.
To the best of our knowledge, this work is the first result on the global-in-time existence of solutions to the fluid--rigid body interaction problem of an Ericksen--Leslie liquid crystal model in three dimensions.

The rest of the paper is organized as follows. 
\autoref{sec:problem} introduces the model, the weak formulations, and the main results, including energy estimates for the rigid displacement. 
In \autoref{sec:GL-approx}, we construct solutions to the moving domain-problem of the Ginzburg--Landau regularization by removing a Galerkin approximation, a penalization of rigidity, and an artificial viscosity. 
This section also develops the fixed-domain Ginzburg--Landau theory on arbitrary bounded open sets.
\autoref{sec:eps to 0} then carries out the singular limit~$\eps\to0$ on the moving domain.
In \autoref{sec:continuation of sol}, the solution obtained from the singular limit passage is continued up to the first contact time and then also continued after impact, completing the proof of \autoref{exist:main}. 
Auxiliary results are proved in \autoref{sec:proofs aux lemmas}.

\section{Statement of the problem and main results}\label{sec:problem}

First of all, let us introduce some notations. We use $C$ for a general positive constant, $C(a,b)$ for a positive constant that may depend on $a$ and $b$ and may change from line to line. For a measurable set $Q$, let
\[
    \Char(Q):=\{\chi\in\rL^\infty(Q):\chi\in\{0,1\}\text{ a.e. in }Q\}
\]
denote the class of characteristic functions on $Q$.
For a fixed spatial domain $G$, we abbreviate~${\rL^p_tX_{\bfx}:=\rL^p(0,T;X(G))}$ and $\rL^p_{t,\bfx}:=\rL^p((0,T)\times G)$, and we use the same abbreviations for norms. Here $G=\Omega$ unless another domain is specified, and moving spatial domains are always displayed explicitly. We set $\rH^{-1}(G):=(\rH_0^1(G))'$. Time slices of Bochner functions are evaluated at their precise representatives on a common set of full measure, unless a continuous representative is specified.
\subsection{The simplified Ericksen--Leslie model coupled with a rigid body}

We consider a rigid body interacting with a simplified Ericksen--Leslie fluid. The fluid density, viscosity, and solid density are normalized to one, and vector-valued quantities are written in boldface.

\subsubsection{Geometric setting}\label{sssec:geometry}

Let $\Omega \subset \R^3$ be an open, bounded, and connected domain with boundary~$\partial\Omega$ of class $\rC^2$, representing the container in which the system evolves. A rigid body occupies, at initial time, an open, bounded, and connected domain $\cS_0 \subset \Omega$ with $\rC^2$-boundary, whose center of mass is located at the origin, i.e.,
\begin{equation*}
    \frac{1}{|\cS_0|}\int_{\cS_0} \bfx \srd\bfx = 0.
\end{equation*}
Let $T>0$ be a given final time. At each time $t \in (0,T)$, the position of the rigid body is characterized by its center of mass $\bfh(t) \in \R^3$ and a rotation matrix $\bfO(t) \in \SO(3)$, so that the solid domain is given by
\begin{equation*}
    \cS(t) = \bigl\{\bfh(t) + \bfO(t)\bfy : \bfy \in \cS_0\bigr\}.
\end{equation*}
Before contact, $\overline{\cS(t)}\subset\Omega$, and the fluid domain at time $t$ is
\begin{equation*}
    \cF(t) \coloneqq \Omega \setminus \overline{\cS(t)},
\end{equation*}
with initial fluid domain $\cF_0 \coloneqq \Omega \setminus \overline{\cS_0}$.

The rigid velocity is
\begin{equation*}
    \bfu^{\cS}(t,\bfx) = \bfell(t) + \bfomega(t) \times (\bfx - \bfh(t)), \quad (t,\bfx) \in (0,T) \times \cS(t),
\end{equation*}
where $\bfell(t) = \bfh'(t)$ denotes the translational velocity and $\bfomega(t)$ the angular velocity.

To prescribe the director anchoring on the moving interface, let $\bfZ(t,\bfy)$ be the rigid flow map. Denote by $\bfm(t)$ the skew-symmetric matrix satisfying $\bfm(t)\bfx=\bfomega(t)\times\bfx$. Then
\begin{equation}\label{eq:CP flow map}
    \left\{
    \begin{aligned}
        \del_t \bfZ(t,\bfy)
        &= \bfm(t)\bigl(\bfZ(t,\bfy) - \bfh(t)\bigr) + \bfell(t), &&\tin (0,T) \times \R^3,\\
        \bfZ(0,\bfy)
        &= \bfy, &&\tfor \bfy \in \R^3.
    \end{aligned}
    \right.
\end{equation}
Let $\bfO(t)\in\SO(3)$ solve $\bfO'(t)=\bfm(t)\bfO(t)$ with $\bfO(0)=\I_3$. Since $\bfh(0)=0$, the solution of \eqref{eq:CP flow map} is
\[
    \bfZ(t,\bfy)  = \bfh(t) + \bfO(t) \bfy.
\]
Write $\bfZ_t(\bfy):=\bfZ(t,\bfy)$. Then $\cS(t)=\bfZ_t(\cS_0)$.
The inverse $(\bfZ_t)^{-1}$ solves
\begin{equation}\label{eq:CP inverse flow map}
    \left\{
    \begin{aligned}
        \del_t (\bfZ_t)^{-1}(\bfx)
        &= -\bfO^\top(t) \bfm(t) \bfO(t) (\bfZ_t)^{-1}(\bfx) - \bfO^\top(t) \bfell(t), &&\tin (0,T) \times \R^3,\\
        (\bfZ_0)^{-1}(\bfx)
        &= \bfx, &&\tfor \bfx \in \R^3.
    \end{aligned}
    \right.
\end{equation}

\subsubsection{Simplified Ericksen--Leslie equations}

The fluid phase is governed by the simplified Ericksen--Leslie system introduced by Lin \cite{Lin:89} and Lin--Liu \cite{LL:95}. The unknowns in the fluid region are the velocity $\bfu^{\cF}$, the pressure $\pi$, and the director field $\bfd$. They are defined on $(0,T)\times\cF(t)$. Their values lie in $\R^3$, $\R$, and $\bS^2$, respectively, where
\begin{equation*}
    \bS^2 = \bigl\{\bff \in \R^3 : |\bff| = 1\bigr\}
\end{equation*}
denotes the unit sphere in $\R^3$. The condition $|\bfd|=1$ is the director constraint.

The equations in the fluid domain are
\begin{equation}\label{eq:simplified EL}
    \left\{
    \begin{aligned}
        \partial_t \bfu^{\cF} + (\bfu^{\cF} \cdot \nabla)\bfu^{\cF} - \Delta\bfu^{\cF} + \nabla\pi
        &= -\mdiv\bigl((\nabla\bfd)^\top \nabla\bfd\bigr), &&\text{in } (0,T) \times \cF(t),\\
        \mdiv\bfu^{\cF} &= 0, &&\text{in } (0,T) \times \cF(t),\\
        \partial_t \bfd + (\bfu^{\cF} \cdot \nabla)\bfd - \Delta\bfd
        &= |\nabla\bfd|^2\,\bfd, &&\text{in } (0,T) \times \cF(t),\\
        |\bfd| &= 1, &&\text{in } (0,T) \times \cF(t).
    \end{aligned}
    \right.
\end{equation}
For $|\bfd|=1$, the scalar $|\nabla\bfd|^2$ is the Lagrange multiplier associated with the unit-length constraint in~\eqref{eq:simplified EL}$_3$, and $-\mdiv((\nabla\bfd)^\top\nabla\bfd)$ in \eqref{eq:simplified EL}$_1$ is the Ericksen elastic force. The isotropic part of the elastic stress has been absorbed into the pressure.
\subsubsection{Stress tensor}

To formulate the fluid--structure coupling through Newton's laws, we introduce the relevant stress tensors. Let $\bD(\bfu^{\cF}) = \frac{1}{2}(\nabla\bfu^{\cF} + (\nabla\bfu^{\cF})^\top)$ denote the symmetric part of the gradient.
The Newtonian stress tensor is given by
\begin{equation*}
    \bT_N(\bfu^{\cF},\pi) = 2\,\bD(\bfu^{\cF}) - \pi\,\I_3,
\end{equation*}
where $\I_3$ denotes the $3 \times 3$ identity matrix. The Ericksen elastic stress is
\begin{equation*}
    \bT_E(\bfd) = -(\nabla\bfd)^\top\nabla\bfd.
\end{equation*}
The total Cauchy stress is
\begin{equation}\label{eq:stress tensor}
    \bT(\bfu^{\cF},\pi,\bfd) \coloneqq \bT_N(\bfu^{\cF},\pi) + \bT_E(\bfd) = 2\,\bD(\bfu^{\cF}) - \pi\,\I_3 - (\nabla\bfd)^\top\nabla\bfd.
\end{equation}
\subsubsection{Rigid body dynamics}

With unit solid density, the mass and center of mass are
\begin{equation*}
    m = |\cS(t)| \tand \bfh(t) = \frac{1}{m} \int_{\cS(t)} \bfx \srd \bfx,
\end{equation*}
and the moment of inertia is defined by
\begin{equation*}
    [\bJ(t) \bfa] \cdot \bb = \int_{\cS(t)} (\bfa \times (\bfx-\bfh(t))) \cdot (\bb \times (\bfx-\bfh(t))) \srd \bfx,
    \qquad \bfa,\bb\in\R^3.
\end{equation*}
By Newton's laws of motion, the translational and rotational dynamics of the rigid body satisfy
\begin{equation}\label{eq:rigid body}
    \begin{aligned}
        m\,\bfell'(t) &= -\int_{\partial\cS(t)} \bT(\bfu^{\cF},\pi,\bfd)\cdot\bfn\srd\sigma,\\[6pt]
        \bJ(t)\,\bfomega'(t) &= \bJ(t)\,\bfomega(t) \times \bfomega(t) - \int_{\partial\cS(t)} (\bfx - \bfh(t)) \times \bT(\bfu^{\cF},\pi,\bfd)\cdot\bfn\srd\sigma,
    \end{aligned}
\end{equation}
where $\bfn$ denotes the unit outward normal to the fluid domain $\cF(t)$ on the interface $\partial\cS(t)$ (equivalently, the unit inward normal to $\cS(t)$), and $\mathrm{d}\sigma$ is the surface measure.
\subsubsection{Coupling, boundary, and initial conditions}

At the fluid--solid interface and at the container boundary, we impose the no-slip conditions
\begin{equation*}
    \bfu^{\cF} = \bfu^{\cS} \ton (0,T) \times \partial\cS(t) \tand
    \bfu^{\cF} = 0 \ton (0,T) \times \partial\Omega.
\end{equation*}
We denote the upper hemisphere and its closed convex hull respectively by
\begin{equation*}
    \bS_+^2 \coloneqq \bigl\{\bff = (f_1,f_2,f_3) \in \R^3 : |\bff| = 1,\; f_3 \ge 0\bigr\},
\end{equation*}
\begin{equation}\label{e:Kplus}
 K_+:=\bigl\{\bff\in\R^3:|\bff|\le1,\ f_3\ge0\bigr\}.   
\end{equation}
The anchoring is described by a single finite-energy map $\bfd_0\in \rH^1(\Omega;\bS_+^2)$
which, as in \cite{LW:16}, serves both as initial director and as Dirichlet datum. Its traces on $\partial\Omega$ and $\partial\cS_0$ are understood in $\rH^{1/2}$. Only $\bfd_0|_{\cF_0}$ is the physical initial director. The values of $\bfd_0$ on $\cS_0$ prescribe the inner boundary datum in material coordinates and define the extension of the director through the solid. The Dirichlet energy of this extension in the solid is invariant under the rigid motion and drops out of the energy balance. Anchoring maps in $\rH^2(\Omega;\bS^2_+)$ are used only within the approximation scheme. This additional regularity is removed by the hemisphere-preserving approximation of \autoref{lem:hemisphere preparation}.
With the flow map $\bfZ_t$ and its inverse $(\bfZ_t)^{-1}$ as introduced in \eqref{eq:CP inverse flow map}, the interface and boundary conditions for the director field are respectively given by
\begin{equation*}
    \begin{aligned}
        \bfd(t,\bfx)
        &= \bfd_0((\bfZ_t)^{-1}(\bfx)), &&\tfor (t,\bfx) \in (0,T)\times \partial \cS(t), \tand\\
        \bfd(t,\bfx)
        &= \bfd_0(\bfx), &&\tfor (t,\bfx) \in (0,T)\times\partial \Omega.
    \end{aligned}
\end{equation*}
The anchoring values on $\partial\cS(t)$ are transported with the body but not rotated. In particular, they remain in $\bS^2_+$, and the transported lifting used in \autoref{sec:GL-approx} is constant along particle paths.
Finally, we prescribe the initial director in the fluid and the initial position of the body by
\begin{equation*}
    \bfd(0,\cdot)=\bfd_0\quad\text{in }\cF_0,
    \qquad \bfh(0)=0,
\end{equation*}
and we use a unified initial velocity on the whole container,
\begin{equation}\label{eq:compatibility}
 \bfu_0\in\cH(\Omega),\qquad
 \bfu_0(\bfx)=\boldsymbol\ell_0+\boldsymbol\omega_0\times\bfx
 \quad\text{for a.e.\ }\bfx\in\cS_0,
\end{equation}
where the space $\cH(\Omega)$ will be made precise in \eqref{eq:init data space fluid vel} below.
Thus the physical initial fluid velocity is~$\bfu_0|_{\cF_0}$, while the restriction to the solid determines the initial translational and angular velocities~$\boldsymbol\ell_0,\boldsymbol\omega_0$ uniquely. In particular, no boundary trace of the finite-energy datum $\bfu_0$ is assumed. If one additionally assumes $\rH^1$-regularity of the initial velocity together with the classical boundary traces, then this unified formulation is consistent with the usual no-slip conditions on $\partial\cS_0$ and $\partial\Omega$.

\subsection{Weak solution and statements of main results}\label{ssec:main result}

We first introduce the function spaces and the notion of weak solution for the constrained system on an arbitrary bounded open set (\autoref{def:postcontact constrained weak solution}), which governs the evolution after contact and whose spaces are reused later, and then the weak formulation of the coupled fluid--rigid-body problem before contact (\autoref{def:weak sol}). Our main results are \autoref{exist:main} on global existence through wall contact, \autoref{prop:impact law} on the impact law, and \autoref{thm:arbitrary open constrained} on existence on arbitrary bounded open sets.

We use a whole-container formulation for the velocity. Let
\begin{equation}\label{eq:init data space fluid vel}
    \cH(\Omega)
    :=\Bigl\{\bfv\in\rL^2(\Omega;\R^3):\mdiv\bfv=0\text{ in }\cD'(\Omega),\
    \bfv\cdot\bfn=0\text{ in }\rH^{-1/2}(\partial\Omega)\Bigr\},
\end{equation}
and set
\begin{equation*}
    \cV(\Omega):=\bigl\{\bfv\in\rH_0^1(\Omega;\R^3):\mdiv\bfv=0\bigr\}.
\end{equation*}
For any connected open set $\cS\subset\Omega$ representing a rigid configuration, set
\begin{equation}\label{eq:rigid kinetic space general}
    \cH_{\cS}(\Omega)
    :=\bigl\{\bfv\in\cH(\Omega):\bD(\bfv)=0\text{ in }\cD'(\cS)\bigr\}.
\end{equation}
The last condition is meaningful for $\bfv\in\rL^2(\Omega;\R^3)$ in the distributional sense and, because $\cS$ is connected, implies that the restriction of $\bfv$ to $\cS$ is a rigid velocity. In particular, for $\bfu_0\in\cH_{\cS_0}(\Omega)$, the initial translational and angular velocities $\bfell_0$ and $\bfomega_0$ appearing in \eqref{eq:compatibility} are the unique coefficients of this rigid restriction.

For any bounded open set $G\subset\R^3$, denote
\begin{equation}\label{eq:arbitrary domain solenoidal spaces}
 \cD_G:=\{\bfv\in\rC_c^\infty(G;\R^3):\mdiv\bfv=0\},\qquad
 \cH_G:=\overline{\cD_G}^{\,\rL^2(G)},\qquad
 \cV_G:=\overline{\cD_G}^{\,\rH_0^1(G)}.
\end{equation}
These definitions require no regularity of $\partial G$ and will be used both for the residual contact domain and for the arbitrary-open-set existence theorem below. For given $0<T< \infty$, define
\[
 \cX_G:=\{\bfpsi\in \rH^1(0,T;\rL^2(G;\R^3))\cap \rL^2(0,T;\rH_0^1(G;\R^3)):\bfpsi(T)=0\}.
\]
We consider the system restricted to the fixed domain $G$ given by
\begin{equation}\label{eq:simplified EL-fixed dom}
    \left\{
    \begin{aligned}
        \partial_t \bfu + (\bfu \cdot \nabla)\bfu - \Delta\bfu + \nabla\pi
        &= -\mdiv\bigl((\nabla\bfd)^\top \nabla\bfd\bigr), &&\text{in } (0,T) \times G,\\
        \mdiv\bfu &= 0, &&\text{in } (0,T) \times G,\\
        \partial_t \bfd + (\bfu \cdot \nabla)\bfd - \Delta\bfd
        &= |\nabla\bfd|^2\,\bfd, &&\text{in } (0,T) \times G,\\
        |\bfd| &= 1, &&\text{in } (0,T) \times G.
    \end{aligned}
    \right.
\end{equation}

\begin{defn}\label{def:postcontact constrained weak solution}
Let $G\subset\R^3$ be bounded and open, $\bfu_0\in\cH_G$, and
$\bfd_0\in\rH^1(G;\bS_+^2)$. A weak solution to \eqref{eq:simplified EL-fixed dom} on $[0,T]$ with initial data $(\bfu_0,\bfd_0)$ is a pair $(\bfu,\bfd)$ with
\begin{align*}
 &\bfu\in \rL^\infty(0,T;\cH_G)\cap \rL^2(0,T;\cV_G),\\
 &\bfd\in \rL^\infty(0,T;\rH^1(G;\bS_+^2)),\qquad
   \bfd-\bfd_0\in \rL^\infty(0,T;\rH_0^1(G;\R^3)),
\end{align*}
whose tension $\bftau:=\Delta \bfd+|\nabla \bfd|^2\bfd$, defined in the sense of distributions, belongs to $\rL^2((0,T)\times G;\R^3)$, such that the following hold.
\begin{enumerate}[label = (\roman*)]
    \item For every $\bfzeta\in \rC_c^1([0,T);\cD_G)$ and every $\bfpsi\in\cX_G$,
\begin{align*}
 &\int_0^T\int_G\bigl[-\bfu\cdot\partial_t\bfzeta-(\bfu\otimes \bfu):\nabla\bfzeta
       +2\bD(\bfu):\bD(\bfzeta)\bigr] \srd\bfx\srd t\\
 &\qquad\qquad=\int_G\bfu_0\cdot\bfzeta(0)\srd\bfx
 +\int_0^T\int_G(\nabla \bfd)^\top\nabla \bfd:\nabla\bfzeta\srd\bfx\srd t,\\
 &\int_0^T\int_G\bigl[-\bfd\cdot\partial_t\bfpsi+(\bfu\cdot\nabla)\bfd\cdot\bfpsi
       -\bftau\cdot\bfpsi\bigr] \srd\bfx\srd t
 =\int_G\bfd_0\cdot\bfpsi(0)\srd\bfx.
\end{align*}
\item For a.e.\ $t\in(0,T)$,
\begin{equation}\label{eq:postcontact constrained energy}
 \begin{aligned}
&\frac12\|\bfu(t)\|_{\rL^2(G)}^2+\frac12\|\nabla\bfd(t)\|_{\rL^2(G)}^2 +2\int_0^t\int_G|\bD(\bfu)|^2\srd\bfx\srd s\\ &\qquad+\int_0^t\int_G|\bftau|^2\srd\bfx\srd s
\le\frac12\|\bfu_0\|_{\rL^2(G)}^2+\frac12\|\nabla\bfd_0\|_{\rL^2(G)}^2.
\end{aligned}
\end{equation}
\end{enumerate}
A \emph{global weak solution} is a pair defined on $[0,\infty)$ which is a weak solution on $[0,T]$ for every $T>0$.
\end{defn}
For a set $A\subset\Omega$ and $r>0$, write
\[
    A_r:=\{\bfx\in\Omega:\dist(\bfx,A)<r\}.
\]
Whenever a characteristic function $\chi$ is used to encode a solid, an open set $S$ is part of the data and we use the pointwise representative $\chi=\dsone_S$. We then set $\cS(\chi):=S$. Thus $\cS(\chi)_r$ denotes the open~$r$-collar of the represented solid.
For an open set $A\subset\Omega$, define
\[
    \cK(A):=\{\bfv\in\cV(\Omega):\bD(\bfv)=0\text{ a.e.\ in }A\}.
\]
For a represented characteristic function $\chi=\dsone_{\cS(\chi)}$, set
\[
    \cK(\chi):=\cK(\cS(\chi)),\qquad
    \cK_r(\chi):=\cK(\cS(\chi)_r).
\]
For a spacetime characteristic function $\varphi$, we write
\[
\rL^2(0,T;\cK(\varphi))
:=\{\bfv\in \rL^2(0,T;\cV(\Omega)):\varphi\bD(\bfv)=0\text{ a.e.\ in }(0,T)\times\Omega\}.
\]
The unified velocity is
\begin{equation*}
    \bfu(t,\bfx)=\begin{cases}
      \bfu^{\cF}(t,\bfx),&\bfx\in\cF(t),\\
      \bfu^{\cS}(t,\bfx),&\bfx\in\cS(t).
    \end{cases}
\end{equation*}
For a rigid trajectory $\cS(t)=\bfh(t)+\bfO(t)\cS_0$, recall that $\cF(t)=\Omega\setminus\overline{\cS(t)}$ is the corresponding fluid region, and introduce the fluid spacetime
\[
    \cQ_{\cF}:=\{(t,\bfx)\in(0,T)\times\Omega:\bfx\in\cF(t)\}.
\]
The standard momentum test class is
\begin{equation}\label{eq:standard momentum tests}
\cT_{\bfu}:=\bigl\{\bfzeta\in \rC_c^1([0,T);\rC_c^\infty(\Omega;\R^3)):
\mdiv\bfzeta=0,\ \bD(\bfzeta(t))=0\text{ on }\cS(t)\text{ for every }t\in[0,T]\bigr\}.
\end{equation}

For the director equation, we use the molecular field and the moving energy test space
\begin{equation}\label{eq:director energy tests}
\begin{aligned}
\cX_{\bfd}(\varphi):=\bigl\{&\bfpsi\in \rH^1(0,T;\rL^2(\Omega;\R^3))\cap \rL^2(0,T;\rH_0^1(\Omega;\R^3)):\ \bfpsi(T)=0,\\
&\bfpsi(t)=0\ \text{a.e.\ in }\cS(t),\quad
 \bfpsi(t)|_{\cF(t)}\in \rH_0^1(\cF(t);\R^3)\ \text{for a.e.\ }t\bigr\}.
\end{aligned}
\end{equation}
The smooth localized class
\begin{equation}\label{eq:smooth loc test class dir}
    \cT_{\bfd}:=
    \bigl\{\bfpsi\in\rC_c^1([0,T)\times\Omega;\R^3):
    \supp\bfpsi\cap((0,T)\times\partial\Omega)=\varnothing,\
    \supp\bfpsi\Subset \cQ_{\cF}\cup(\{0\}\times\cF_0)\bigr\}
\end{equation}
is dense in \eqref{eq:director energy tests} (see \autoref{lem:density director tests}).  The tension formulation below is well defined on this larger energy test space, whereas the formulation containing $|\nabla\bfd|^2\bfd$ is used only for localized smooth tests.

The pre-contact weak solution is defined as follows.

\begin{defn}\label{def:weak sol}
Let $\bfu_0\in\cH_{\cS_0}(\Omega)$, $\bfd_0\in \rH^1(\Omega;\bS_+^2)$, and $\varphi_0=\dsone_{\cS_0}$. A pre-contact weak solution of \eqref{eq:simplified EL} coupled to \eqref{eq:rigid body}, with the stated boundary and initial conditions, on $[0,T]$ is a triple $(\bfu,\bfd,\varphi)$ together with a rigid trajectory satisfying
\[
    \bfh\in\rW^{1,\infty}(0,T;\R^3),\qquad
    \bfO\in\rW^{1,\infty}(0,T;\R^{3\times3}),\qquad
    \bfO(t)\in\SO(3)\ \text{for every }t,
\]
such that
\begin{align*}
    &\bfh(0)=0,\qquad \bfO(0)=\I_3,\qquad
      \varphi(t)=\dsone_{\cS(t)},\qquad \cS(t)=\bfh(t)+\bfO(t)\cS_0,
      \qquad \overline{\cS(t)}\subset\Omega,\\
    &\bfu\in\rL^\infty(0,T;\cH(\Omega))\cap\rL^2(0,T;\cK(\varphi)),\\
    &\bfd\in\rL^\infty(0,T;\rH^1(\cF(t);\bS_+^2)),\\
    &\bftau:=\Delta\bfd+|\nabla\bfd|^2\bfd\in\rL^2(\cQ_{\cF};\R^3)
      \quad\text{in the sense of distributions},
\end{align*}
and the following hold.
\begin{enumerate}[label = (\roman*)]
    \item The restriction of $\bfu$ to $\cS(t)$ is the rigid velocity
\[
    \bfu(t,\bfx)=\bfh'(t)+\bfomega(t)\times(\bfx-\bfh(t))
    \quad\text{for a.e.\ }\bfx\in\cS(t),
\]
where $\bfO'=\bfm\bfO$ and $\bfm\bfx=\bfomega\times\bfx$.
\item The director satisfies the strong anchoring condition in the trace sense,
\begin{equation}\label{eq:weak trace anchoring}
    \operatorname{Tr}_{\partial\Omega}\bfd(t)=\bfd_0,\qquad
    \operatorname{Tr}_{\partial\cS(t)}\bfd(t,\bfx)=\bfd_0\bigl(\bfO(t)^\top(\bfx-\bfh(t))\bigr)\quad\text{for a.e.\ }t.
\end{equation}
\item For every $\bfzeta\in\cT_{\bfu}$, $\bfpsi\in\cX_{\bfd}(\varphi)$, and $\eta\in\rC_c^1([0,T)\times\Omega)$,
\begin{align}
 &\int_0^T\int_\Omega\Bigl[-\bfu\cdot\partial_t\bfzeta-(\bfu\otimes\bfu):\nabla\bfzeta
     +2\bD(\bfu):\bD(\bfzeta)\Bigr]\srd\bfx\srd t \notag\\
 &\qquad \qquad=\int_\Omega\bfu_0\cdot\bfzeta(0)\srd\bfx
 +\int_0^T\int_{\cF(t)}(\nabla\bfd)^\top\nabla\bfd:\nabla\bfzeta\srd\bfx\srd t,
 \label{eq:weak sol u}\\[3pt]
 &\int_0^T\int_{\cF(t)}\Bigl[-\bfd\cdot\partial_t\bfpsi
 +(\bfu\cdot\nabla)\bfd\cdot\bfpsi-\bftau\cdot\bfpsi\Bigr]\srd\bfx\srd t
 =\int_{\cF_0}\bfd_0\cdot\bfpsi(0)\srd\bfx,
 \label{eq:weak sol d}\\[3pt]
 &\int_0^T\int_\Omega\varphi\,[\partial_t\eta+\bfu\cdot\nabla\eta]\srd\bfx\srd t
 =-\int_\Omega\varphi_0\eta(0)\srd\bfx.
 \label{eq:weak sol solid}
\end{align}
\item Finally, for a.e.\ $t\in(0,T)$, we have the energy inequality
\begin{equation}\label{eq:energy-main}
    E(t)+2\int_0^t\int_\Omega|\bD(\bfu)|^2\srd\bfx\srd s
    +\int_0^t\int_{\cF(s)}|\bftau|^2\srd\bfx\srd s
    \le E(0),
\end{equation}
where
\begin{align*}
    E(t)&:=\frac12\int_\Omega|\bfu(t)|^2\srd\bfx
    +\frac12\int_{\cF(t)}|\nabla\bfd(t)|^2\srd\bfx\quad\text{for a.e.\ }t\in(0,T),\\
    E(0)&:=\frac12\int_\Omega|\bfu_0|^2\srd\bfx
    +\frac12\int_{\cF_0}|\nabla\bfd_0|^2\srd\bfx.
\end{align*}
\end{enumerate}
\end{defn}

Note that $\bfh(0)=0$ in the above definition means precisely that the center of mass at time zero is the origin. For $\bfpsi\in\cT_{\bfd}$, \eqref{eq:weak sol d} is equivalently
\[
\int_0^T\int_{\cF(t)}\Bigl[-\bfd\cdot\partial_t\bfpsi
 +(\bfu\cdot\nabla)\bfd\cdot\bfpsi+\nabla\bfd:\nabla\bfpsi
 -|\nabla\bfd|^2\bfd\cdot\bfpsi\Bigr]\srd\bfx\srd t
 =\int_{\cF_0}\bfd_0\cdot\bfpsi(0)\srd\bfx.
\]

For data satisfying the hypotheses below, the initial distance to the wall and the initial energy~$E(0)$ in \eqref{eq:energy-main} are denoted by
\begin{equation}\label{eq:delta_0 & E_0}
 \delta_0:=\dist(\cS_0,\partial\Omega),\tand
 E_0:=\frac12\|\bfu_0\|_{\rL^2(\Omega)}^2
      +\frac12\|\nabla\bfd_0\|_{\rL^2(\cF_0)}^2.
\end{equation}
For a bounded open set $G$ and $(\bfv,\bfq)\in\rL^2(G;\R^3)\times\rH^1(G;\R^3)$, write
\begin{equation*}
 \mathcal E_G(\bfv,\bfq)
 :=\frac12\|\bfv\|_{\rL^2(G)}^2
 +\frac12\|\nabla\bfq\|_{\rL^2(G)}^2.
\end{equation*}
For a collision-free rigid trajectory and a director satisfying the material-coordinate anchoring, we denote by $\overline\bfd$ the material extension
\begin{equation}\label{eq:material extension notation main}
 \overline\bfd(t,\bfx):=
 \begin{cases}
  \bfd(t,\bfx),&\bfx\in\cF(t),\\
  \bfd_0\bigl(\bfO(t)^\top(\bfx-\bfh(t))\bigr),&\bfx\in\cS(t).
 \end{cases}
\end{equation}

The main result provides global weak solutions under the perfectly inelastic sticking law.

\begin{thm}\label{exist:main}
Assume $\partial\Omega,\partial\cS_0\in\rC^2$, $\cS_0$ is connected,
$\delta_0>0$, and
\[
 \bfu_0\in\cH_{\cS_0}(\Omega),\qquad
 \bfd_0\in\rH^1(\Omega;\bS_+^2).
\]
Then there exist $T_{\rm col}\in(0,\infty]$, a rigid trajectory
$\cS(t)=\bfh(t)+\bfO(t)\cS_0$ and a triple $(\bfu,\bfd,\varphi)$ on~$[0,T_{\rm col})$, and, if $T_{\rm col}<\infty$, a pair $(\bfu,\bfd)$ on $[T_{\rm col},\infty)$, denoted by the same symbols, with the following properties.
\begin{enumerate}[(i)]
\item For every $T<T_{\rm col}$, the triple $(\bfu,\bfd,\varphi)$ is a weak solution on $[0,T]$ in the sense of \autoref{def:weak sol}, and $\dist(\cS(t),\partial\Omega)>0$ for $0\le t<T_{\rm col}$.
\item If $T_{\rm col}<\infty$, the rigid motion extends continuously to $T_{\rm col}$ and, with $\cS_c:=\cS(T_{\rm col})$ and $\Omega_c:=\Omega\setminus\overline{\cS_c}$,
\begin{equation}\label{eq:collision alternative main}
    \dist(\cS_c,\partial\Omega)=0.
\end{equation}
As $t\to T_{\rm col}$, the material extension \eqref{eq:material extension notation main} converges to some $\bfd_c\in\rH^1(\Omega;\bS_+^2)$, strongly in~$\rL^2(\Omega)$ and weakly in $\rH^1(\Omega)$, and there is a unique $\bfu_c^+\in\cH_{\Omega_c}$ such that $(\bfu(t),\bfzeta)_{\rL^2(\Omega)}\to(\bfu_c^+,\bfzeta)_{\rL^2(\Omega_c)}$ for every $\bfzeta\in\cD_{\Omega_c}$.
\item If $T_{\rm col}<\infty$, the body stays at $\cS_c$ for $t\ge T_{\rm col}$, the velocity vanishes in $\cS_c$, and $(\bfu,\bfd)(T_{\rm col}+\cdot\,)|_{\Omega_c}$ is a global weak solution on $\Omega_c$ in the sense of \autoref{def:postcontact constrained weak solution} with initial data $(\bfu_c^+,\bfd_c|_{\Omega_c})$. Moreover, with $\cF(s):=\Omega_c$ for $s\ge T_{\rm col}$ and the dissipation
\[
 \mathcal D(0,t):=2\int_0^t\int_\Omega|\bD(\bfu)|^2\srd\bfx\srd s+\int_0^t\int_{\cF(s)}|\bftau|^2\srd\bfx\srd s,
\]
we have
\begin{equation}\label{eq:global constrained energy main}
 \mathcal E_{\Omega_c}(\bfu(t),\bfd(t))+\mathcal D(0,t)\le E_0\qquad\text{for a.e.\ }t>T_{\rm col}.
\end{equation}
\end{enumerate}
\end{thm}

The first contact time is bounded from below explicitly. By \autoref{cor:quantitative separation main} in \autoref{ssec:uniform precontact}, $T_{\rm col}\ge T_{\rm ed}$, where $T_{\rm ed}>0$, defined in \eqref{eq:T ed main theorem}, depends only on $\delta_0$, $E_0$, $\Omega$, and $\cS_0$.

If contact occurs, the impact law takes the following form. For a connected rigid configuration~${\cS\subset\Omega}$, let $G_{\cS}:=\Omega\setminus\overline{\cS}$, let $P_{\cH_{G_{\cS}}}$ be the $\rL^2(G_{\cS})$-orthogonal projection onto $\cH_{G_{\cS}}$, and define the impact dissipation of $\bfv\in\cH_{\cS}(\Omega)$ by
\begin{equation}\label{eq:impact dissipation main}
 \mathcal J_{\rm imp}(\cS,\bfv)
 :=\frac12\|\bfv\|_{\rL^2(\cS)}^2
 +\frac12\bigl\|\bfv|_{G_{\cS}}
 -P_{\cH_{G_{\cS}}}(\bfv|_{G_{\cS}})\bigr\|_{\rL^2(G_{\cS})}^2.
\end{equation}
It is the kinetic energy of the body plus the kinetic energy of the part of the fluid velocity that is not admissible in the residual fluid domain~$G_{\cS}$.

\begin{prop}\label{prop:impact law}
Let $T_{\rm col}<\infty$ in \autoref{exist:main}, and let $\bfu(t_j)\rightharpoonup\bfu_c^-$ in $\rL^2(\Omega)$ along a sequence~${t_j\to T_{\rm col}}$ in the full-measure set of times at which $\bfu(t_j)$ is rigid on $\cS(t_j)$ and the energy inequality \eqref{eq:energy-main} holds. Then
\begin{equation}\label{eq:impact projection main}
    \bfu_c^-\in\cH_{\cS_c}(\Omega)\tand \bfu_c^+=P_{\cH_{\Omega_c}}(\bfu_c^-|_{\Omega_c}),
\end{equation}
and consequently
\begin{equation}\label{eq:impact energy decomposition main}
    \frac12\|\bfu_c^-\|_{\rL^2(\Omega)}^2=\frac12\|\bfu_c^+\|_{\rL^2(\Omega_c)}^2+\mathcal J_{\rm imp}(\cS_c,\bfu_c^-).
\end{equation}
Moreover, the sequence $(t_j)$ can be chosen such that, in addition,
\begin{equation}\label{eq:impact energy drop main}
 \mathcal E_{\Omega_c}(\bfu_c^+,\bfd_c)+\mathcal J_{\rm imp}(\cS_c,\bfu_c^-)\le\operatorname*{ess\,liminf}_{t\to T_{\rm col}}E(t),
\end{equation}
and the energy inequality \eqref{eq:global constrained energy main} holds with $\mathcal J_{\rm imp}(\cS_c,\bfu_c^-)$ added to its left-hand side.
\end{prop}

\begin{rem}\label{rem:main results}
In \autoref{prop:impact law}, the rigid velocity need not converge as $t\to T_{\rm col}$, so the weak pre-impact trace $\bfu_c^-$ and the impact dissipation $\mathcal J_{\rm imp}(\cS_c,\bfu_c^-)$ may depend on the chosen sequence, whereas the post-impact state $(\bfu_c^+,\bfd_c)$ does not. 
We observe that the inequality \eqref{eq:impact energy drop main} may be strict, and the difference indicates a possible loss of compactness of the kinetic or elastic energy at contact.
\end{rem}

The post-contact construction uses the following fixed-domain existence theorem, which is independent of the rigid-body geometry and requires no boundary regularity.

\begin{thm}\label{thm:arbitrary open constrained}
Consider a bounded and open set $G\subset\R^3$ with no regularity assumption on $\partial G$, and recall the spaces $\cD_G,\cH_G,\cV_G$ from \eqref{eq:arbitrary domain solenoidal spaces}. 
If $\bfu_0\in\cH_G$ and $\bfd_0\in \rH^1(G;\bS_+^2)$, then there exists a global weak solution in the sense of \autoref{def:postcontact constrained weak solution}.
\end{thm}

\subsection{Energy control of the rigid displacement}\label{ssec:uniform precontact}

In this subsection, we discuss a lower bound for the collision time $T_{\rm col}$ by the time $T_{\rm ed}$ associated with the energy dissipation, see \eqref{eq:T ed main theorem} below.

We observe that the energy inequality controls the motion of the body in two ways. 
The kinetic energy bounds the rigid velocity, and the viscous dissipation bounds its time integral upon using Poincar\'e's and Korn's inequalities, see \eqref{eq:Korn Poincare constant} below.
The estimates of this subsection are used at every level of the construction.
They yield collision-free intervals that are independent of the approximation parameters.
Throughout the paper, we use the following two constants.
First, by $C_K=C_K(\Omega)>0$ we denote a Korn--Poincar\'e constant, i.e.,
\begin{equation}\label{eq:Korn Poincare constant}
    \|\bfv\|_{\rL^2(\Omega)}\le C_K\|\bD(\bfv)\|_{\rL^2(\Omega)} \tfor \bfv\in\cV(\Omega).
\end{equation}
Next, we define
\begin{equation}\label{eq:rigid velocity constant}
 C_{\cS_0}:=m^{-1/2}+\operatorname{diam}(\cS_0)\,\lambda_{\min}\bigl(\bJ(0)\bigr)^{-1/2},
\end{equation}
where $\lambda_{\min}(\bJ(0))>0$ denotes the smallest eigenvalue of the inertia tensor of $\cS_0$. For $E_0>0$, let~$T_{\rm ed}>0$ be the unique solution of
\begin{equation}\label{eq:T ed main theorem}
 T_{\rm ed}\bigl(1-e^{-4T_{\rm ed}/C_K^2}\bigr)
 =\frac{2\delta_0^2}{C_{\cS_0}^2C_K^2E_0}.
\end{equation}
The existence of $T_{\rm ed}>0$ is due to the fact that the left-hand side of \eqref{eq:T ed main theorem} increases strictly from zero to infinity. 
If $E_0=0$, we set~$T_{\rm ed}:=+\infty$.

\begin{lem}\label{lem:rigid velocity energy control}
Let $\cS(t)=\bfh(t)+\bfO(t)\cS_0$ denote the position of the rigid body, where we assume that $\bfh$ and $\bfO$ are Lipschitz continuous, and where $\bfO'=\bfm\bfO$ and $\bfm\bfx=\bfomega\times\bfx$, capturing the motion of the rigid body.
Moreover, suppose that the fluid and rigid body velocities coincide, i.e.,
\[
   \bfu(t,\bfx)=\bfh'(t)+\bfomega(t)\times(\bfx-\bfh(t))\qquad\text{for a.e.\ }\bfx\in\cS(t).
\]
Then, with $C_{\cS_0}$ from \eqref{eq:rigid velocity constant}, for a.e.\ $t$, we obtain the bound
\begin{equation*}
   \sup_{\bfy\in\cS_0}\left|\frac{\rd}{\rd t}\bigl(\bfh(t)+\bfO(t)\bfy\bigr)\right|
   \le C_{\cS_0}\|\bfu(t)\|_{\rL^2(\cS(t))}.
\end{equation*}
\end{lem}

\begin{proof}
First, since by assumption the center of mass of $\cS_0$ is the origin and the density equals one, we get
\[
   \int_{\cS(t)}|\bfu|^2\srd\bfx
   =m|\bfh'|^2+\bfomega\cdot\bJ(t)\bfomega .
\]
Moreover, $\bJ(t)=\bfO(t)\bJ(0)\bfO(t)^\top$ has the same eigenvalues as $\bJ(0)$, so that $|\bfh'|\le m^{-1/2}\|\bfu\|_{\rL^2(\cS(t))}$ and $|\bfomega|\le\lambda_{\min}(\bJ(0))^{-1/2}\|\bfu\|_{\rL^2(\cS(t))}$. As the origin lies in the convex hull of $\cS_0$, we have $|\bfy|\le\operatorname{diam}(\cS_0)$ for $\bfy\in\cS_0$, and therefore
\[
   \left|\frac{\rd}{\rd t}\bigl(\bfh(t)+\bfO(t)\bfy\bigr)\right|
   =|\bfh'+\bfomega\times\bfO(t)\bfy|
   \le |\bfh'|+\operatorname{diam}(\cS_0)|\bfomega|
   \le C_{\cS_0}\|\bfu(t)\|_{\rL^2(\cS(t))}.\qedhere
\]
\end{proof}

\begin{lem}\label{lem:dissipative rigid displacement}
Let $T>0$, and let $\cS(t)$ and $\bfu$ be as in \autoref{lem:rigid velocity energy control} on $[0,T]$, with $\bfh(0)=0$ and~${\bfO(0)=\I_3}$, such that $\overline{\cS(t)}\subset\Omega$ and $\bfu\in\rL^\infty(0,T;\cH(\Omega))\cap\rL^2(0,T;\cV(\Omega))$. Assume that, for a.e.~$t\in(0,T)$,
\begin{equation}\label{eq:kinetic energy dissipation hypothesis}
 \frac12\|\bfu(t)\|_{\rL^2(\Omega)}^2
 +2\int_0^t\|\bD(\bfu(s))\|_{\rL^2(\Omega)}^2\srd s
 \le E_0.
\end{equation}
Then, for every $0\le t\le T$,
\begin{equation}\label{eq:coupled rigid displacement}
 d_{\rm H}(\cS(t),\cS_0)
 \le C_{\cS_0}C_K
 \sqrt{\frac{E_0t}{2}\bigl(1-e^{-4t/C_K^2}\bigr)}
 \le C_{\cS_0}\min\left\{\sqrt{2E_0}\,t,\ C_K\sqrt{\frac{E_0t}{2}}\right\},
\end{equation}
and consequently
\begin{equation}\label{eq:coupled separation}
 \dist(\cS(t),\partial\Omega)
 \ge \delta_0-C_{\cS_0}C_K
 \sqrt{\frac{E_0t}{2}\bigl(1-e^{-4t/C_K^2}\bigr)}.
\end{equation}
\end{lem}

\begin{proof}
Set $Y(t):=\int_0^t\|\bfu(s)\|_{\rL^2(\Omega)}^2\srd s$.
It readily follows that $Y'(t)=\|\bfu(t)\|_{\rL^2(\Omega)}^2$ for a.e.\ $t$, and the assumption \eqref{eq:kinetic energy dissipation hypothesis} together with the Korn--Poincar\'e inequality \eqref{eq:Korn Poincare constant} yields, for a.e.\ $t\in(0,T)$, the estimate $\frac12Y'(t)+\frac{2}{C_K^2}Y(t)\le E_0$.
Multiplying by $e^{4t/C_K^2}$ integrating, and using $Y(0)=0$, we get
\[
 Y(t)\le \frac{C_K^2E_0}{2}\bigl(1-e^{-4t/C_K^2}\bigr).
\]
By \autoref{lem:rigid velocity energy control} and the Cauchy--Schwarz inequality in time, every $\bfy\in\cS_0$ satisfies the estimate
\[
    |\bfh(t)+\bfO(t)\bfy-\bfy| \le C_{\cS_0}\int_0^t\|\bfu(s)\|_{\rL^2(\Omega)}\srd s \le C_{\cS_0}\sqrt{tY(t)},
\]
proving the first inequality in \eqref{eq:coupled rigid displacement}.
The second inequality in \eqref{eq:coupled rigid displacement} then follows from the elementary estimate $1-e^{-x}\le\min\{x,1\}$ for $x\ge0$. 
Finally, $\dist(\cS(t),\partial\Omega)\ge\delta_0-d_{\rm H}(\cS(t),\cS_0)$.
\end{proof}

Every energy inequality in this paper contains the kinetic energy and the viscous dissipation on its left-hand side, so \autoref{lem:dissipative rigid displacement} applies at all approximation levels. The construction mostly uses the interval of length
\begin{equation}\label{eq:T star main theorem}
 T_*:=
 \begin{cases}
 \dfrac{\delta_0}{2C_{\cS_0}\sqrt{2E_0}},&E_0>0,\\[6pt]
 +\infty,&E_0=0.
 \end{cases}
\end{equation}

\begin{cor}\label{cor:uniform separation}
Under the assumptions of \autoref{lem:dissipative rigid displacement},
\[
   \dist(\cS(t),\partial\Omega)\ge \frac{\delta_0}{2}
   \quad(0\le t\le \min\{T,T_*\}),
   \qquad
   \dist(\cS(t),\partial\Omega)>0
   \quad(0\le t<\min\{T,T_{\rm ed}\}).
\]
\end{cor}

\begin{proof}
The first assertion follows from the bound $C_{\cS_0}\sqrt{2E_0}\,t$ in \eqref{eq:coupled rigid displacement}.
The second follows from~\eqref{eq:coupled separation}, the definition \eqref{eq:T ed main theorem}, and the strict monotonicity of $t\mapsto t(1-e^{-4t/C_K^2})$.
\end{proof}

Applied to the solution constructed in this paper, \autoref{lem:dissipative rigid displacement} gives the following explicit separation estimate. Its proof is completed in \autoref{ssec:maximal constrained}.

\begin{cor}\label{cor:quantitative separation main}
The solution of \autoref{exist:main} satisfies, for every $0\le t<T_{\rm col}$,
\begin{equation}\label{eq:energy dissipation rigid displacement}
 \dist(\cS(t),\partial\Omega)
 \ge \delta_0-C_{\cS_0}C_K
 \sqrt{\frac{E_0t}{2}\bigl(1-e^{-4t/C_K^2}\bigr)}.
\end{equation}
In particular,
\begin{equation}\label{eq:sharp collision lower bound main}
    T_{\rm col}\ge T_{\rm ed},
\end{equation}
and if $E_0=0$, then the body is at rest and $T_{\rm col}=\infty$.
\end{cor}

\begin{rem}\label{rem:T ed comparison}
Since $1-e^{-x}<\min\{x,1\}$ for $x>0$, the definition \eqref{eq:T ed main theorem} implies, for $E_0>0$,
\[
 T_{\rm ed}>\max\left\{\frac{\delta_0}{C_{\cS_0}\sqrt{2E_0}},\ \frac{2\delta_0^2}{C_{\cS_0}^2C_K^2E_0}\right\}.
\]
\end{rem}

\section{The Ginzburg--Landau approximation}\label{sec:GL-approx}

In this section, $\eps>0$ is fixed, and we construct weak solutions of the Ginzburg--Landau relaxation on the moving domain, see \autoref{thm:weak sol GL}. The construction uses three further approximation levels, namely a Galerkin approximation of the velocity (index $k$), a penalization of the rigidity constraint and of the anchoring condition in the solid (index $n$), and an $\rH^3$ artificial viscosity regularization (parameter $\delta$). These are removed in the order $k\to\infty$, $n\to\infty$, and $\delta\to0$. For fixed $\delta$, the approximate velocities are Lipschitz and vanish on $\partial\Omega$, so the body cannot reach the wall in finite time, see \eqref{eq:lipschitz flow separation}. Separation uniformly in $\delta$ then follows from \autoref{ssec:uniform precontact}.

Define
\begin{equation*}
    F_{\eps}(\bfd):=\frac{1}{4\eps^2}(1-|\bfd|^2)^2 \enspace \text{ and } \enspace \bff_\eps (\bfd) \coloneqq \nabla_{\bfd}F_{\eps}(\bfd)=-\frac{1}{\eps^2}(1-|\bfd|^2)\bfd.
\end{equation*}
The Ginzburg--Landau approximation then reads
\begin{equation}\label{eq:simplified EL-GL}
    \left\{
    \begin{aligned}
        \del_t \bfu^{\cF} + (\bfu^{\cF} \cdot \nabla)\bfu^{\cF} + \nabla\pi
        &= \Delta\bfu^{\cF} - \mdiv\bigl((\nabla\bfd)^\top\nabla\bfd\bigr), &&\tin (0,T)\times\cF(t),\\
        \mdiv\bfu^{\cF}
        &= 0, &&\tin (0,T)\times\cF(t),\\
        \del_t\bfd + (\bfu^{\cF}\cdot\nabla)\bfd
        &= \Delta\bfd + \dfrac{1}{\eps^2}\bigl(1-|\bfd|^2\bigr)\bfd, &&\tin (0,T)\times\cF(t).
    \end{aligned}
    \right.
\end{equation}
The stress tensor $\bT(\bfu^\cF,\pi,\bfd)$ remains as in~\eqref{eq:stress tensor}, and the rigid body dynamics are governed by
\begin{equation*}
    \begin{aligned}
        m\,\bfell'(t)
        &= -\int_{\del\cS(t)} \bT(\bfu^{\cF},\pi,\bfd)\cdot\bfn\srd\sigma,\\
        \bJ(t)\,\bfomega'(t)
        &= \bJ(t)\,\bfomega(t)\times\bfomega(t) - \int_{\del\cS(t)} (\bfx-\bfh(t))\times\bT(\bfu^{\cF},\pi,\bfd)\cdot\bfn\srd\sigma.
    \end{aligned}
\end{equation*}
The coupling and boundary conditions are
\begin{equation*}
    \begin{aligned}
        \bfu^\cF
        &= \bfu^\cS \text{ on } (0,T)\times\del\cS(t),
        \qquad &&\bfu^\cF = 0 \text{ on } (0,T)\times\del\Omega,\\[2pt]
        \bfd(t,\bfx)
        &= \bfd_0((\bfZ_t)^{-1}(\bfx)) \text{ for } t \in (0,T), \enspace \bfx \in \del \cS(t), \qquad &&\bfd = \bfd_0 \text{ on } (0,T)\times\del\Omega,
    \end{aligned}
\end{equation*}
where we invoked the flow map $\bfZ_t$ as well as its inverse $(\bfZ_t)^{-1}$ from \eqref{eq:CP flow map} and \eqref{eq:CP inverse flow map}, respectively.
The initial conditions are
\begin{equation}\label{eq:init conds GL}
    \bfu^\cF(0,\bfx) = \bfu_0(\bfx) \;\tin\cF_0,\quad
    \bfd(0,\bfx) = \bfd_0(\bfx) \;\tin\cF_0,\quad
    \bfh(0) = 0,\quad
    \bfh'(0) = \bfell_0,\quad
    \bfomega(0) = \bfomega_0.
\end{equation}

For the function spaces defined in \autoref{ssec:main result}, the weak solution of the approximate problem~\mbox{\eqref{eq:simplified EL-GL}--\eqref{eq:init conds GL}} is defined as follows.

\begin{defn}\label{def:weak sol GL}
Let $\bfu_0\in\cH_{\cS_0}(\Omega)$, $\bfd_0\in\rH^2(\Omega;\bS_+^2)$, and $\varphi_0=\dsone_{\cS_0}$. 
A weak Ginzburg--Landau solution consists of a rigid trajectory and a triple $(\bfu,\bfd,\varphi)$ satisfying all geometric, kinematic, velocity, anchoring, and solid-transport requirements of \autoref{def:weak sol}. Moreover, $\bfd\in\rL^\infty(0,T;\rH^1(\cF(t);\R^3))$ satisfies $|\bfd|\le1$ and $d_3\ge0$ as well as
\begin{equation*}
\int_0^T\int_{\cF(t)}\Bigl[-\bfd\cdot\partial_t\bfpsi+(\bfu\cdot\nabla)\bfd\cdot\bfpsi
+\nabla\bfd:\nabla\bfpsi+\bff_\eps(\bfd)\cdot\bfpsi\Bigr]\srd\bfx\srd t
=\int_{\cF_0}\bfd_0\cdot\bfpsi(0)\srd\bfx
\end{equation*}
for every $\bfpsi\in\cX_{\bfd}(\varphi)$. 
Furthermore, for every $\bfzeta\in\cT_{\bfu}$ and $\eta\in\rC_c^1([0,T)\times\Omega)$, the momentum and transport identities read
\begin{align*}
 &\int_0^T\int_\Omega\Bigl[-\bfu\cdot\partial_t\bfzeta-(\bfu\otimes\bfu):\nabla\bfzeta
 +2\bD(\bfu):\bD(\bfzeta)\Bigr]\srd\bfx\srd t\\
 &\qquad=\int_\Omega\bfu_0\cdot\bfzeta(0)\srd\bfx
 +\int_0^T\int_{\cF(t)}(\nabla\bfd)^\top\nabla\bfd:\nabla\bfzeta\srd\bfx\srd t,
 \\
 &\int_0^T\int_\Omega\varphi[\partial_t\eta+\bfu\cdot\nabla\eta]\srd\bfx\srd t
 =-\int_\Omega\varphi_0\eta(0)\srd\bfx.
\end{align*}
The strong anchoring condition \eqref{eq:weak trace anchoring} is imposed in the trace sense.
\end{defn}

We are in a position to state the result on the existence of a weak solution to the Ginzburg--Landau system for fixed $\eps > 0$.

\begin{thm}\label{thm:weak sol GL}
Assume $\bfu_0\in\cH_{\cS_0}(\Omega)$ and $\bfd_0\in \rH^2(\Omega;\bS_+^2)$, let $E_0$ be as in \eqref{eq:delta_0 & E_0}, and let $T_*$ be given by \eqref{eq:T star main theorem}. If $E_0=0$, we replace $T_*$ by an arbitrary finite time. Then, for every $\eps>0$, the Ginzburg--Landau system \eqref{eq:simplified EL-GL}--\eqref{eq:init conds GL} admits a weak solution on $[0,T_*]$ in the sense of \autoref{def:weak sol GL} such that $\dist(\cS(t),\partial\Omega)\ge\delta_0/2$ for $0\le t\le T_*$ and which satisfies, for a.e.\ $t\in(0,T_*)$, the energy inequality
\begin{equation}\label{eq:energy-main GL}
 E^{\rm GL}_\eps(t)+2\int_0^t\int_\Omega|\bD(\bfu)|^2\srd\bfx\srd s
 +\int_0^t\int_{\cF(s)}\left|\Delta\bfd-\bff_\eps(\bfd)\right|^2\srd\bfx\srd s
 \le E^{\rm GL}_\eps(0),
\end{equation}
where
\begin{equation*}
 E^{\rm GL}_\eps(t):=\frac12\int_\Omega|\bfu(t)|^2\srd\bfx
 +\frac12\int_{\cF(t)}|\nabla\bfd(t)|^2\srd\bfx
 +\int_{\cF(t)}F_\eps(\bfd(t))\srd\bfx.
\end{equation*}
\end{thm}
We note here that the initial value $E^{\rm GL}_\eps(0)$ in the above theorem is independent of $\eps$ since $|\bfd_0|=1$.

\medskip
The material anchoring condition is imposed by the following transported lifting, which also preserves the fixed outer boundary datum. Let $\bfX_{s,t}$ be the volume-preserving flow generated by
$\bfu\in \rL^1(0,T;\rW^{1,\infty}(\Omega;\R^3))\cap \rC([0,T];\cV(\Omega))$, and put $\bfY(t):=\bfX_{t,0}$. Since the trace of $\bfu$ vanishes on $\partial\Omega$, the flow fixes the outer boundary. Define
\begin{equation}\label{eq:phi & b}
    \varphi(t):=\varphi_0\circ\bfY(t),\qquad
 \bb(t):=\bfd_0\circ\bfY(t).
\end{equation}
Then we have
\[
\partial_t\varphi+\bfu\cdot\nabla\varphi = \partial_t\bb+\bfu\cdot\nabla\bb = 0,\qquad \varphi(t)=\dsone_{\bfX_{0,t}(\cS_0)},
\]
and
$\bb(t,\bfX_{0,t}(\bfy))=\bfd_0(\bfy)$ for $\bfy\in\partial\cS_0$, and $\bb=\bfd_0$ on $\partial\Omega$. 

Let us define $X_u:=\cV(\Omega)\cap\rH^3(\Omega;\R^3),$ then $X_u$ is a closed subspace of $\rH^3(\Omega;\R^3)$.
We can fix a Hilbert inner product $a_u(\cdot,\cdot)$ on $X_u$, which is equivalent to the
$\rH^3$-inner product. We write $\|\bfv\|_{X_u}^2:=a_u(\bfv,\bfv).$
The regularizing term is $\delta a_u(\bfu,\bfv)$. Its energy contribution controls the full $X_u$-norm and hence, by the embedding $X_u\hookrightarrow\rW^{1,\infty}(\Omega;\R^3)$, whose norm we denote by $C_X$, the Lipschitz norm used to construct the approximate flow maps.

\subsection{Penalization and artificial viscosity}
The rigidity of the body and the anchoring of the director inside the solid are enforced by penalization. The dissipation term $n\varphi|\bD(\bfu)|^2$ penalizes deformations of the solid, and the energy term $\frac n2\varphi|\bfd-\bb|^2$ penalizes deviations of the director from the transported lifting. The artificial viscosity $\delta a_u$ provides the Lipschitz regularity of the velocity that is needed to transport $\varphi$ and $\bb$.

\begin{defn}\label{def:weak sol penalized}
Let $\bfu_0 \in \cH_{\cS_0}(\Omega)$, $\bfd_0 \in \rH^2(\Omega;\bS_+^2)$, and let $\varphi_0$ denote the characteristic function of~$\cS_0$. A quadruple $(\bfu^n,\bfd^n,\varphi^n,\bb^n)$ with
\begin{align*}
    \bfu^n &\in \rL^{\infty}(0,T;\cH(\Omega)) \cap \rL^2(0,T;\cV(\Omega)\cap\rH^3(\Omega;\R^3)),\\[3pt]
    \bfd^n &\in \rL^\infty(0,T;\rH^1(\Omega;\R^3)) \cap \rL^2(0,T;\rH^2(\Omega;\R^3)),\\ 
    \bfd^n - \bfd_0 &\in \rL^\infty(0,T;\rH_0^1(\Omega;\R^3)), \quad \bfd^n(t,\bfx)\in K_+ \quad  \text{for a.e.\ } (t,\bfx),\\[3pt]
    \varphi^n &\in \Char([0,T]\times\Omega) \cap \rC([0,T];\rL^p(\Omega)) \quad\text{for every }p\in[1,\infty),\\[3pt]
    \bb^n &\in \rC([0,T];\rL^p(\Omega;\R^3)) \quad\text{for every }p\in[1,\infty)
\end{align*}
is called a \emph{weak solution to the penalized problem} if,
for all test functions
    \begin{align*}
        \bfzeta &\in \rH^1(0,T;\rL^2(\Omega;\R^3))\cap\rL^2(0,T;\cV(\Omega)\cap\rH^3(\Omega;\R^3)),\quad \bfzeta(T)=0,\\
        \bfpsi &\in \rH^1(0,T;\rL^2(\Omega;\R^3))\cap\rL^2(0,T;\rH_0^1(\Omega;\R^3)),\quad \bfpsi(T)=0,\\
        \eta &\in \rC^1([0,T]\times\overline{\Omega}),\quad \eta(T)=0,\\
        \bfxi &\in \rC^1([0,T]\times\overline{\Omega};\R^3),\quad \bfxi(T)=0,
    \end{align*}
    one has
    \begin{align}
        &\int_0^T\int_{\Omega}\bigl[-\bfu^n\cdot\del_t\bfzeta - (\bfu^n\otimes\bfu^n):\nabla\bfzeta + (2+n\varphi^n)\bD(\bfu^n):\bD(\bfzeta)\bigr]\srd\bfx\srd t + \delta\int_0^T a_u(\bfu^n,\bfzeta)\srd t \notag\\
        &\hspace{2cm} = \int_{\Omega}\bfu_0\cdot\bfzeta(0)\srd\bfx - \int_0^T\int_{\Omega}(\bfzeta\cdot\nabla)\bfd^n\cdot\Delta \bfd^n\srd\bfx\srd t,
        \label{eq:penalized fluid}\\[8pt]
        &\int_0^T\int_{\Omega}\Bigl[-\bfd^n\cdot\del_t\bfpsi + (\bfu^n\cdot\nabla)\bfd^n\cdot\bfpsi - \Delta\bfd^n\cdot\bfpsi + \bff_\eps (\bfd^n) \cdot\bfpsi\Bigr]\srd\bfx\srd t \notag\\
        &\hspace{2cm} = \int_\Omega \bfd_0 \cdot \bfpsi(0) \srd \bfx - \int_0^T \int_\Omega n \varphi^n (\bfd^n - \bb^n) \cdot \bfpsi \srd \bfx \srd t,
        \label{eq:penalized d}\\[8pt]
        &\int_0^T\int_{\Omega}\varphi^n\bigl[\partial_t \eta + (\bfu^n\cdot\nabla)\eta\bigr]\srd\bfx\srd t = -\int_{\Omega}\varphi_0\,\eta(0)\srd\bfx, \tand
        \label{eq:penalized solid}\\
        &\int_0^T\int_{\Omega}\bb^n \cdot \bigl[\partial_t\bfxi + (\bfu^n\cdot\nabla)\bfxi\bigr]\srd\bfx\srd t = -\int_{\Omega}\bfd_0 \cdot\bfxi(0)\srd\bfx.
        \label{eq:penalized lifting}
    \end{align}
\end{defn}

We define $\bfM^n$ by
\begin{equation}\label{eq:n-M-def}
    \bfM^n:=
    \Delta\bfd^n
    +\frac{1}{\eps^2}(1-|\bfd^n|^2)\bfd^n
    -n\varphi^n(\bfd^n-\bb^n).
\end{equation}

\begin{thm}\label{thm:ex of sol to penalized probl}
Let $\bfu_0 \in \cH_{\cS_0}(\Omega)$, $\bfd_0 \in \rH^2(\Omega;\bS_+^2)$, and let $\varphi_0$ be the characteristic function of~$\cS_0$. Then, for every $T>0$, the penalized problem~\eqref{eq:penalized fluid}--\eqref{eq:penalized lifting} admits a solution $(\bfu^n,\bfd^n,\varphi^n,\bb^n)$ on $[0,T]$ in the sense of \autoref{def:weak sol penalized} which satisfies, for a.e.\ $t\in(0,T)$, the energy inequality
\begin{equation}\label{eq:n-energy-inequality}
    \begin{aligned}
        E^{n,\delta,\eps}(t)
        &+\int_0^t\int_{\Omega}(2+n\varphi^n)|\bD(\bfu^n)|^2\srd\bfx\srd s
        +\delta\int_0^t\|\bfu^n(s)\|_{X_u}^2\srd s\\
        &+\int_0^t\int_{\Omega}|\bfM^n|^2\srd\bfx\srd s
        \le E^{n,\delta,\eps}(0)=E_{\rm pen,0},
    \end{aligned}
\end{equation}
where the approximate energy $E^{n,\delta,\eps}$ is defined by
\begin{equation*}
    \begin{aligned}
    E^{n,\delta,\eps}(t):=
        &\frac12\int_{\Omega}|\bfu^n|^2\srd\bfx
        +\int_{\Omega}\left(
            \frac12|\nabla\bfd^n|^2
            +F_\eps(\bfd^n)
        \right)\srd\bfx +\frac n2\int_{\Omega}\varphi^n|\bfd^n-\bb^n|^2\srd\bfx.
    \end{aligned}
\end{equation*}
Here
\[
    E_{\rm pen,0}:=\frac12\int_\Omega|\bfu_0|^2\srd\bfx
    +\frac12\int_\Omega|\nabla\bfd_0|^2\srd\bfx,
\]
which is independent of $n$, $\delta$, and $\eps$ because $|\bfd_0|=1$ and the initial director penalty vanishes. Moreover, $\varphi^n(t)=\dsone_{\cS^n(t)}$ with $\cS^n(t):=\bfX^n_{0,t}(\cS_0)$, where $\bfX^n$ denotes the flow of $\bfu^n$, and
\begin{equation}\label{eq:lipschitz flow separation}
 \dist(\cS^n(t),\partial\Omega)\ge\rho_\delta(t):=\delta_0\exp\Bigl(-C_X\sqrt{tE_{\rm pen,0}/\delta}\Bigr)
 \qquad(0\le t\le T).
\end{equation}
\end{thm}

\subsubsection{Semi-Galerkin approximation}\label{ssec:semi-Galerkin approx}

At the first approximation level, the parameters $n\in\N$, $\delta>0$, and $\eps>0$ are fixed, and only the Galerkin dimension $k\in\N$ varies. As in the semi-Galerkin scheme of \cite[Sec.~2]{LL:95}, only the velocity is discretized, while the director solves a parabolic problem. The moving body requires two modifications. The solid indicator and the inhomogeneous director datum are both transported by the flow of the Galerkin velocity, the latter through the lifting $\bb^k$ in~\eqref{eq:k-lift}.

\smallskip
Let us recall that $\cF_0=\Omega\setminus\overline{\cS_0}$ and $\varphi_0:=\dsone_{\cS_0}$. 
At this approximation level we work with the whole-container anchoring map
\[
    \bfd_0\in \rH^2(\Omega;\R^3),
    \qquad |\bfd_0| = 1 \quad \text{a.e.\ in }\Omega.
\]

\smallskip

We also recall the Hilbert space $X_u =\cV(\Omega)\cap\rH^3(\Omega;\R^3)$  and the inner product $a_u$ fixed above. 
The usual solenoidal density theorem on the $\rC^2$-container implies that the embedding $X_u\hookrightarrow \cH(\Omega)$
is dense, continuous, and compact.  The spectral theorem for this compact Hilbert embedding follows that there exist $0<\lambda_1\le\lambda_2\le\cdots$,
$\lambda_j\to\infty$, and an $\rL^2(\Omega)$-orthonormal basis
$\{\bfv_j\}_{j\in\N}\subset X_u$ of $\cH(\Omega)$ satisfying
\begin{equation*}
    a_u(\bfv_j,\bfw)=\lambda_j(\bfv_j,\bfw)_{\rL^2(\Omega)}
    \qquad\text{for every }\bfw\in X_u.
\end{equation*}
Consequently, if $\cP_k$ denotes the $\rL^2$-orthogonal projection onto
$\cV_k:=\Span\{\bfv_1,\ldots,\bfv_k\}$, then
\begin{equation}\label{eq:velocity spectral projection}
    \|\cP_k\bfw\|_{X_u}\le \|\bfw\|_{X_u},
    \qquad
    \cP_k\bfw\rightarrow\bfw\quad\text{strongly in }X_u
    \quad(\bfw\in X_u).
\end{equation}
Indeed, writing $\bfw=\sum_j c_j\bfv_j$ in $\rL^2(\Omega)$ we have $a_u(\bfw,\bfw)=\sum_j\lambda_j|c_j|^2$.  
In particular, we imply that $\bfv_j\in X_u\hookrightarrow\rW^{1,\infty}(\Omega;\R^3)$, while the finite spans $\cV_k$ are dense in $\cH(\Omega)$ with respect to the~$\rL^2$-norm and in $X_u$ with respect to the $\rH^3$-norm.
The approximate velocity has the form
\begin{equation}\label{eq:approx vel}
    \bfu^k(t,\bfx)=\sum_{j=1}^k a_j^k(t)\bfv_j(\bfx).
\end{equation}
In particular, we have $\mdiv\bfu^k=0$ in $\Omega$ and $\bfu^k=0$ on $\partial\Omega$.

\smallskip

For a given velocity $\bfu^k\in \rC([0,T];\cV_k)$, let
$\bfX^k_{s,t}\colon\Omega\to\Omega$ be the corresponding flow map, i.e.,
\begin{equation}\label{eq:k-flow-map}
    \frac{\rd}{\rd t}\bfX^k_{s,t}(\bfy)
    =\bfu^k(t,\bfX^k_{s,t}(\bfy)),
    \qquad
    \bfX^k_{s,s}(\bfy)=\bfy.
\end{equation}
We denote the inverse map by
\begin{equation}\label{eq:k-inverse flow-map}
    \bfY^k(t,\bfx):=\bfX^k_{t,0}(\bfx)=(\bfX^k_{0,t})^{-1}(\bfx).
\end{equation}
Note that the flow map $\bfX^k_{s,t}$ is volume preserving and fixes $\partial\Omega$. The transported solid indicator, the approximate solid, and the transported lifting are respectively defined by
\begin{align}
    \varphi^k(t,\bfx)&:=\varphi_0(\bfY^k(t,\bfx)),\qquad
    \cS^k(t):=\bfX^k_{0,t}(\cS_0),\qquad
    \varphi^k(t)=\dsone_{\cS^k(t)},\label{eq:k-phi}\\
    \bb^k(t,\bfx)&:=\bfd_0(\bfY^k(t,\bfx)).\label{eq:k-lift}
\end{align}
Denote $\rD_t^k:=\del_t+\bfu^k\cdot\nabla$, we have, in the sense of distributions,
\begin{equation}\label{eq:k transport identities}
    \rD_t^k\varphi^k=0,\qquad
    \rD_t^k\bb^k=0
    \quad\tin(0,T)\times\Omega,\qquad
    \varphi^k(0)=\varphi_0,\quad
    \bb^k(0)=\bfd_0.
\end{equation}
Since the flow fixes $\partial\Omega$, we have $\bb^k=\bfd_0$ on $(0,T)\times\partial\Omega$, and $\bb^k(t,\bfX^k_{0,t}(\bfy))=\bfd_0(\bfy)$ for $\bfy\in\partial\cS_0$. 

\smallskip

For fixed $\bfu^k$, $\varphi^k$, and $\bb^k$, the director field
$\bfd^k$ is defined as the solution of the semilinear parabolic problem
\begin{equation}\label{eq:k-director}
    \left\{
    \begin{aligned}
        \del_t\bfd^k+(\bfu^k\cdot\nabla)\bfd^k
        &= \Delta\bfd^k
        +\frac{1}{\eps^2}(1-|\bfd^k|^2)\bfd^k
        -n\varphi^k(\bfd^k-\bb^k)
        \qquad\tin (0,T)\times\Omega,\\
        \bfd^k&=\bfd_0\hspace{2.5cm}
    \ton ((0,T)\times\partial\Omega) \cup (\{0\}\times \Omega).
    \end{aligned}
    \right.
\end{equation}
We define the approximate molecular field $\bfM^{n,k}$ by
\begin{equation}\label{eq:k-M-def}
    \bfM^{n,k}:=
    \Delta\bfd^k
    +\frac{1}{\eps^2}(1-|\bfd^k|^2)\bfd^k
    -n\varphi^k(\bfd^k-\bb^k).
\end{equation}
Then \eqref{eq:k-director} is equivalently
\begin{equation}\label{eq:k-material-d-H}
    \rD_t^k\bfd^k=\bfM^{n,k}.
\end{equation}
Since $\rD_t^k\bb^k=0$, the lifted field $\bfw^k:=\bfd^k-\bb^k$, which vanishes on $\partial\Omega$, satisfies $\rD_t^k\bfw^k=\bfM^{n,k}$.

\smallskip

The velocity coefficients are determined by the following finite-dimensional weak formulation.
For all $\bfv\in\cV_k$, it holds that
\begin{equation}\label{eq:k-momentum}
    \begin{aligned}
        &(\del_t\bfu^k,\bfv)_{\rL^2(\Omega)}+\int_{\Omega} ((\bfu^k\cdot\nabla)\bfu^k)\cdot\bfv\srd\bfx+\int_{\Omega}(2+n\varphi^k)\bD(\bfu^k):\bD(\bfv)\srd\bfx +\delta\,a_u(\bfu^k,\bfv) \\
        &\quad=-\int_{\Omega}(\bfv\cdot\nabla)\bfd^k\cdot\bigl(\bfM^{n,k}+n\varphi^k(\bfd^k-\bb^k)\bigr)\srd\bfx.
    \end{aligned}
\end{equation}
The initial condition is $\bfu^k(0)=\cP_k\bfu_0,$
where $\cP_k\colon\cH(\Omega)\to\cV_k$ denotes the $\rL^2$ orthogonal projection.

\begin{defn}
Let $n\in\N$, $\delta>0$, $\eps>0$, and $k\in\N$ be fixed. A quadruple
$(\bfu^k,\bfd^k,\varphi^k,\bb^k)$ is called a solution to the \emph{semi-Galerkin approximation}
if
\[
    \bfu^k\in \rW^{1,2}(0,T;\cV_k)\cap \rC([0,T];\cV_k),
\]
$\varphi^k$ and $\bb^k$ are given by \eqref{eq:k-phi} and
\eqref{eq:k-lift}, $\bfd^k$ solves \eqref{eq:k-director}, and
\eqref{eq:k-momentum} holds for all $\bfv\in\cV_k$.
\end{defn}

\smallskip

Define the approximate energy $E_k^{n,\delta,\eps}$ by
\begin{equation*}
    \begin{aligned}
    E_k^{n,\delta,\eps}(t):=
        &\frac12\int_{\Omega}|\bfu^k|^2\srd\bfx
        +\int_{\Omega}\left(
            \frac12|\nabla\bfd^k|^2
            +\frac{1}{4\eps^2}(1-|\bfd^k|^2)^2
        \right)\srd\bfx+\frac n2\int_{\Omega}\varphi^k|\bfd^k-\bb^k|^2\srd\bfx.
    \end{aligned}
\end{equation*}

The semi-Galerkin system satisfies the following exact energy identity.

\begin{prop}\label{prop:k-energy-identity}
Every sufficiently smooth semi-Galerkin approximation satisfies, for all
$t\in[0,T]$,
\begin{equation}\label{eq:k-energy-identity}
    \begin{aligned}
        E_k^{n,\delta,\eps}(t)
        &+\int_0^t\int_{\Omega}(2+n\varphi^k)|\bD(\bfu^k)|^2\srd\bfx\srd s
        +\delta\int_0^t\|\bfu^k(s)\|_{X_u}^2\srd s\\
        &+\int_0^t\int_{\Omega}|\bfM^{n,k}|^2\srd\bfx\srd s
        =E_k^{n,\delta,\eps}(0).
    \end{aligned}
\end{equation}
In particular, the estimate is uniform with respect to $n$, provided
$E_k^{n,\delta,\eps}(0)$ is uniformly bounded with respect to $n$.
\end{prop}

\begin{proof}
Testing \eqref{eq:k-momentum} with $\bfv=\bfu^k$ we have
\begin{equation}\label{eq:k-energy-proof-fluid}
    \begin{aligned}
        \frac12\frac{\rd}{\rd t}\int_{\Omega}|\bfu^k|^2\srd\bfx
        &+\int_{\Omega}(2+n\varphi^k)|\bD(\bfu^k)|^2\srd\bfx
        +\delta\|\bfu^k\|_{X_u}^2 \\
        &=-\int_{\Omega}
        (\bfu^k\cdot\nabla)\bfd^k\cdot
        \bigl(\bfM^{n,k}+n\varphi^k(\bfd^k-\bb^k)\bigr)\srd\bfx.
    \end{aligned}
\end{equation}
Using the boundary condition $\bfd^k=\bfd_0$ on $\partial\Omega$, and hence
$\del_t\bfd^k=0$ on $\partial\Omega$, we obtain
\begin{equation*}
    \begin{aligned}
    \frac{\rd}{\rd t}\int_{\Omega}\left(
        \frac12|\nabla\bfd^k|^2+F_{\eps}(\bfd^k)
    \right)\srd\bfx
    =-\int_{\Omega}
    \bigl(\bfM^{n,k}+n\varphi^k(\bfd^k-\bb^k)\bigr)
    \cdot\del_t\bfd^k\srd\bfx.
    \end{aligned}
\end{equation*}
Since $\del_t\bfd^k=\bfM^{n,k}-(\bfu^k\cdot\nabla)\bfd^k$, this yields
\begin{equation}\label{eq:k-director-energy-expanded}
    \begin{aligned}
    \frac{\rd}{\rd t}\int_{\Omega}\left(
        \frac12|\nabla\bfd^k|^2+F_{\eps}(\bfd^k)
    \right)\srd\bfx
    =&-\int_{\Omega}|\bfM^{n,k}|^2\srd\bfx
      -n\int_{\Omega}\varphi^k(\bfd^k-\bb^k)\cdot\bfM^{n,k}\srd\bfx \\
     &+\int_{\Omega}(\bfu^k\cdot\nabla)\bfd^k\cdot
        \bigl(\bfM^{n,k}+n\varphi^k(\bfd^k-\bb^k)\bigr)\srd\bfx.
    \end{aligned}
\end{equation}
Since $\rD_t^k\varphi^k=0$, $\mdiv\bfu^k=0$, and $\bfu^k=0$ on $\partial\Omega$,
we have
\begin{equation*}
    \frac{\rd}{\rd t}\frac n2\int_{\Omega}\varphi^k|\bfd^k-\bb^k|^2\srd\bfx
    =n\int_{\Omega}\varphi^k(\bfd^k-\bb^k)
    \cdot \rD_t^k(\bfd^k-\bb^k)\srd\bfx.
\end{equation*}
Then, combining with \eqref{eq:k transport identities} and \eqref{eq:k-material-d-H}, we get
\begin{equation}\label{eq:k-penalty-energy-derivative}
    \frac{\rd}{\rd t}\frac n2\int_{\Omega}\varphi^k|\bfd^k-\bb^k|^2\srd\bfx
    =n\int_{\Omega}\varphi^k(\bfd^k-\bb^k)\cdot\bfM^{n,k}\srd\bfx.
\end{equation}
Adding \eqref{eq:k-energy-proof-fluid}--\eqref{eq:k-penalty-energy-derivative}, the elastic convection terms cancel,
and the two $n$-weighted terms cancel, so
\[
    \frac{\rd}{\rd t}E_k^{n,\delta,\eps}(t)
    +\int_{\Omega}(2+n\varphi^k)|\bD(\bfu^k)|^2\srd\bfx
    +\delta\|\bfu^k\|_{X_u}^2
    +\int_{\Omega}|\bfM^{n,k}|^2\srd\bfx=0.
\]
Integrating this equation over $(0,t)$ we obtain \eqref{eq:k-energy-identity}, which follows the last statement.
\end{proof}

\smallskip

Next we state the result on the solution existence for the semi-Galerkin system, its proof is given in \autoref{ssec:aux ests Galerkin}.

\begin{prop}\label{prop:k-weak sols}
Let $\bfu_0 \in \cH_{\cS_0}(\Omega)$, $\bfd_0 \in \rH^2(\Omega;\bS_+^2)$, denote by $\varphi_0$ the characteristic function of~$\cS_0$, and assume that the boundaries $\del \Omega$ and $\del \cS_0$ are of class $\rC^2$. For every $k,n\in\N$, $\delta,\eps>0$, and $T>0$, the semi-Galerkin system comprising \eqref{eq:k transport identities}, \eqref{eq:k-director}, and~\eqref{eq:k-momentum} has a solution on $[0,T]$ which satisfies \eqref{eq:k-energy-identity} and $\dist(\cS^k(t),\partial\Omega)\ge\rho_\delta(t)$ for $0\le t\le T$, with $\rho_\delta$ as in \eqref{eq:lipschitz flow separation}.
\end{prop}

We are now ready to prove the existence of a weak solution to the penalized problem~\mbox{\eqref{eq:penalized fluid}--\eqref{eq:penalized lifting}} as asserted in \autoref{thm:ex of sol to penalized probl}.

\begin{proof}[Proof of \autoref{thm:ex of sol to penalized probl}]
The parameters $n$, $\delta$, and $\eps$ are fixed throughout the proof. We will pass to the limit $k\to\infty$ in the semi-Galerkin system. Fix $T>0$, and let $(\bfu^k,\bfd^k,\varphi^k,\bb^k)$ be the solutions on $[0,T]$ provided by \autoref{prop:k-weak sols}, with corresponding flows $\bfX^k_{s,t}$, inverse flows $\bfY^k$, and molecular fields $\bfM^{n,k}$ as in \eqref{eq:k-M-def}. We recall the transport identities \eqref{eq:k transport identities} and the director equation \eqref{eq:k-material-d-H}.

By the Galerkin energy identity \eqref{eq:k-energy-identity} and the director equation \eqref{eq:k-director}, there exists a constant $C>0$ which is independent of $k$ such that
\begin{equation}\label{recall:energy}
    \begin{aligned}
        &\|\bfu^k\|_{\rL^\infty(0,T;\cH(\Omega))}
        +\|\bfu^k\|_{\rL^2(0,T;\cV(\Omega)\cap \rH^3(\Omega))}
        +\|\bfd^k\|_{\rL^\infty(0,T;\rH^1(\Omega))}
        \\
        &\quad
        +\|\bfd^k\|_{\rL^2(0,T;\rH^2(\Omega))}
        +\|\bfM^{n,k}\|_{\rL^2((0,T)\times\Omega)}
        +\sqrt n\,
        \|(\varphi^k)^{1/2}(\bfd^k-\bb^k)\|_{\rL^\infty(0,T;\rL^2(\Omega))}
        \leq C.
    \end{aligned}
\end{equation}
The required time-derivative bounds are stated next and proved in \autoref{ssec:aux ests Galerkin}.

\begin{claim}\label{claim:t-regs u and d}
There exists a constant $C>0$, independent of $k$, such that
\begin{equation}\label{eq:t-regs u & d}
    \|\del_t\bfu^k\|_{\rL^2(0,T;X_u')} \leq C \tand \|\del_t\bfd^k\|_{\rL^2(0,T;\rL^2(\Omega;\R^3))} \leq C.
\end{equation}
In particular, $\del_t\bfd^k$ is bounded in $\rL^2(0,T;(\rH^1(\Omega;\R^3))')$.
\end{claim}

The uniform bounds~\eqref{recall:energy}, the fact that $\varphi^k$ is uniformly bounded by one as a characteristic function, and the Banach--Alaoglu theorem provide subsequences (not relabeled) such that
\begin{align}
    \bfu^k
    &\rightarrow \bfu \text{ weakly in }\rL^2(0,T; \rH^3(\Omega;\R^3)) \text{ and weak-* in }\rL^{\infty}(0,T; \cH(\Omega)),\label{eq:conv uk}\\
    \bfd^k
    &\rightarrow \bfd \text{ weakly in }\rL^2(0,T; \rH^2(\Omega;\R^3))  \text{ and weak-* in }\rL^{\infty}(0,T; \rH^1(\Omega;\R^3)),\label{eq:conv dk}\\
    \bfM^{n,k}
    &\rightarrow \bfM \text{ weakly in }\rL^2(0,T; \rL^2(\Omega;\R^3)), \enspace
    \varphi^k
    \rightarrow \varphi \text{ weak-* in }\rL^{\infty}(0,T; \rL^{\infty}(\Omega)),\label{eq:conv phik}
\end{align}
for some limit functions
\[
\begin{split}
& \bfu\in \rL^2(0,T;\rH^3(\Omega;\R^3))\cap \rL^\infty(0,T;\cH(\Omega)),\\
&\bfd\in \rL^2(0,T; \rH^2(\Omega;\R^3)) \cap \rL^\infty(0,T;\rH^1(\Omega;\R^3)), \\
&\bfM\in \rL^2((0,T)\times\Omega;\R^3), \quad \tand\quad \varphi\in \rL^\infty((0,T)\times\Omega).
\end{split}
\]

By \eqref{recall:energy} and \eqref{eq:t-regs u & d}, the
Banach--Alaoglu theorem and the Aubin--Lions--Simon compactness lemma imply,
up to a non-relabeled subsequence,
\begin{align}
    \bfu^k
    &\to \bfu
    \quad\text{strongly in }\rL^2(0,T;\rH^s(\Omega;\R^3))
    \quad\text{for every }s<3,
    \label{eq:strong conv uk}\\
    \bfd^k
    &\to \bfd
    \quad\text{strongly in }\rL^2(0,T;\rH^{2-\eta}(\Omega;\R^3))
    \quad\text{for every }\eta>0.
    \label{eq:strong conv dk}
\end{align}

We now identify the limits of the transported solid indicator~$\varphi^k$ and of the transported lifting~$\bb^k$. By \eqref{eq:strong conv uk}, choosing $s\in(\frac52,3)$, and
using $\rH^s(\Omega;\R^3)\hookrightarrow \rW^{1,\infty}(\Omega;\R^3)$, we obtain
\begin{equation}\label{eq:uk strong W1inf}
\begin{aligned}
    \bfu^k \to\bfu \text{ strongly in } \rL^1\bigl(0,T;\rW^{1,\infty}(\Omega;\R^3)\bigr).
\end{aligned}
\end{equation}
Indeed, the strong convergence first holds in $\rL^2(0,T;\rW^{1,\infty}(\Omega;\R^3))$, and thus in $\rL^1(0,T;\rW^{1,\infty}(\Omega;\R^3))$.

Let $\bfX_{s,t}$ denote the flow generated by the limiting velocity $\bfu$,
that is,
\[
    \frac{\rd}{\rd t}\bfX_{s,t}(\bfy)
    =
    \bfu(t,\bfX_{s,t}(\bfy)),
    \qquad
    \bfX_{s,s}(\bfy)=\bfy.
\]
We recall that $\bfY(t,\bfx)=\bfX_{t,0}(\bfx)$.
Since $\bfu\in \rL^1(0,T;\rW^{1,\infty}(\Omega;\R^3))$, this flow is
well-defined and unique. Moreover, the stability estimate for flows associated
with Sobolev--Lipschitz vector fields, together with
\eqref{eq:uk strong W1inf}, gives
\begin{equation}\label{eq:flow convergence k}
    \bfX^k_{s,t}\to \bfX_{s,t}
    \quad\text{and}\quad
    \bfY^k(t,\cdot)\to \bfY(t,\cdot)
\end{equation}
uniformly on $[0,T]^2\times\overline\Omega$.
More precisely, for $0\leq s,t\leq T$,
\[
    \sup_{\bfy\in\Omega}
    |\bfX^k_{s,t}(\bfy)-\bfX_{s,t}(\bfy)|
    \leq
    C
    \int_0^T
    \|\bfu^k(\tau)-\bfu(\tau)\|_{\rL^\infty(\Omega)}
    \srd\tau,
\]
where $C$ depends on $\sup_k
    \int_0^T
    \|\nabla\bfu^k(\tau)\|_{\rL^\infty(\Omega)}
    \srd\tau$, but is independent of $k$.

For the above limiting velocity field $\bfu$, we define $\varphi(t,\bfx):=\varphi_0(\bfY(t,\bfx))$ and $\bb(t,\bfx):=\bfd_0(\bfY(t,\bfx))$ as in \eqref{eq:phi & b}.
This identifies the weak-* limit.

First, since $\varphi^k(t,\bfx)=\varphi_0(\bfY^k(t,\bfx))
    =
    \dsone_{\cS_0}(\bfY^k(t,\bfx))$, we claim that
\begin{equation}\label{eq:phi strong conv}
    \varphi^k\to\varphi
    \quad\text{strongly in }
    \rC([0,T];\rL^p(\Omega))
    \quad\text{for every }1\leq p<\infty.
\end{equation}
Indeed, let
\[
    \rho_k
    :=
    \sup_{t\in[0,T]}
    \|\bfY^k(t,\cdot)-\bfY(t,\cdot)\|_{\rL^\infty(\Omega)}.
\]
Then $\rho_k\to0$. For each $t\in[0,T]$, the set on which
$\varphi^k(t)\neq\varphi(t)$ is contained in
\[
    \left\{
        \bfx\in\Omega:
        \dist(\bfY(t,\bfx),\del\cS_0)\leq \rho_k
    \right\}.
\]
Since $\bfu$ is divergence-free and satisfies the no-slip condition on
$\del\Omega$, the flow $\bfY(t,\cdot)$ is measure-preserving. Hence
\[
    \bigl|
    \{\bfx\in\Omega:\varphi^k(t,\bfx)\neq\varphi(t,\bfx)\}
    \bigr|
    \leq
    \bigl|
    \{\bfy\in\Omega:\dist(\bfy,\del\cS_0)\leq \rho_k\}
    \bigr|.
\]
The right-hand side tends to zero because $\del\cS_0$ has zero Lebesgue
measure. Therefore
\[
    \sup_{t\in[0,T]}
    \|\varphi^k(t)-\varphi(t)\|_{\rL^p(\Omega)}^p
    \leq
    \bigl|
    \{\bfy\in\Omega:\dist(\bfy,\del\cS_0)\leq \rho_k\}
    \bigr|
    \to0,
\]
which proves \eqref{eq:phi strong conv}. In particular,
\[
    \varphi(t,\cdot)=\dsone_{\cS(t)}
    \quad\text{for all }t\in[0,T],
    \qquad
    \cS(t):=\bfX_{0,t}(\cS_0).
\]
Thus, the weak-* limit of $\varphi^k$ is again a characteristic function.

Next, we treat the transported lifting. Since $\bfd_0\in \rH^2(\Omega;\R^3)$,
and since in the dimensions under consideration $\rH^2(\Omega;\R^3)\hookrightarrow \rC^{0,\alpha}(\overline\Omega;\R^3)$ for some $\alpha\in(0,\frac{1}{2})$, we may choose the H\"older continuous representative
of $\bfd_0$. Using \eqref{eq:flow convergence k}, we obtain
\[
    \|\bb^k(t)-\bb(t)\|_{\rL^\infty(\Omega)}=
    \|\bfd_0(\bfY^k(t,\cdot))-\bfd_0(\bfY(t,\cdot))\|_{\rL^\infty(\Omega)}
    \leq
    [\bfd_0]_{\rC^{0,\alpha}(\overline\Omega)}
    \|\bfY^k(t,\cdot)-\bfY(t,\cdot)\|_{\rL^\infty(\Omega)}^\alpha.
\]
Consequently,
\begin{equation}\label{eq:b strong conv}
    \bb^k\to\bb
    \quad\text{strongly in }
    \rC([0,T];\rL^\infty(\Omega;\R^3)),
\end{equation}
and hence also strongly in $\rC([0,T];\rL^p(\Omega;\R^3))$ for every $1\leq p<\infty$.
Moreover, since $|\bfd_0|=1$ a.e.\ in~$\Omega$, the transported lift satisfies
\begin{equation*}
    |\bb^k|=1
    \quad\text{and}\quad
    |\bb|=1
    \quad\text{a.e.\ in }(0,T)\times\Omega.
\end{equation*}

Finally, passing to the limit in the transport equations
\eqref{eq:k transport identities}, using the strong convergence of~$\bfu^k$,~$\varphi^k$, and~$\bb^k$, yields
\begin{equation*}
    \del_t\varphi+\mdiv(\varphi\bfu)=0,
    \qquad
    \del_t\bb+\bfu\cdot\nabla\bb=0
\end{equation*}
in $\cD'((0,T)\times\Omega)$, with initial data $ \varphi(0)=\varphi_0,
    \qquad
    \bb(0)=\bfd_0.$
Since $\mdiv\bfu=0$, the first equation is equivalently $\del_t\varphi+\bfu\cdot\nabla\varphi=0$
in the distributional sense.

Now we pass to the limit in the momentum equation.
Fix $m\in\N$ and let $\bfzeta\in\rC^1([0,T];\cV(\Omega)\cap\rH^3(\Omega;\R^3))$ be of the form $\bfzeta(t,\bfx) = \sum_{i=1}^m a^i(t)\,\bfv_i(\bfx)$ with $a^i(T)=0$.  For every $k\ge m$ we may substitute $\bfzeta$ in place of $\bfv$ in~\eqref{eq:k-momentum}, integrate over $(0,T)$, and integrate by parts in time to obtain
\begin{equation*}
    \begin{aligned}
        &\int_0^T\int_{\Omega}\Bigl[-\bfu^k\cdot\del_t\bfzeta - (\bfu^k\otimes\bfu^k):\nabla\bfzeta + (2+n\varphi^k)\bD(\bfu^k):\bD(\bfzeta)\\
        &\quad\quad \qquad+ (\bfzeta\cdot\nabla)\bfd^k\cdot(\bfM^{n,k} + n\varphi^k(\bfd^k - \bb^k))\Bigr]\srd\bfx\srd t+\delta\int_0^T a_u(\bfu^k,\bfzeta)\srd t
        = \int_{\Omega}\bfu^k(0)\cdot\bfzeta(0)\srd\bfx.
    \end{aligned}
\end{equation*}
Since $\bfu^k(0)=\cP_k\bfu_0$ and $\cP_k\bfu_0\to\bfu_0$ strongly in $\rL^2(\Omega)$, the initial term converges to $\int_\Omega\bfu_0\cdot\bfzeta(0)\srd\bfx$. Passing to the remaining terms using the weak and weak-* convergences~\eqref{eq:conv uk}--\eqref{eq:conv phik} of~$\bfu^k$,~$\bfd^k$,~$\bfM^{n,k}$, and~$\varphi^k$, and the strong convergences \eqref{eq:strong conv uk}--\eqref{eq:strong conv dk} of $\bfu^k$ and $\bfd^k$ as well as \eqref{eq:phi strong conv}--\eqref{eq:b strong conv} of $\varphi^k$ and~$\bb^k$, we arrive at the identity
\begin{equation*}
    \begin{aligned}
        &\int_0^T\int_{\Omega}\Bigl[-\bfu\cdot\del_t\bfzeta - (\bfu\otimes\bfu):\nabla\bfzeta + (2+n\varphi)\bD(\bfu):\bD(\bfzeta)\\
        &\hspace{3cm} + (\bfzeta\cdot\nabla)\bfd\cdot(\bfM + n\varphi (\bfd - \bb))\Bigr]\srd\bfx\srd t+\delta\int_0^T a_u(\bfu,\bfzeta)\srd t
        = \int_{\Omega}\bfu_0\cdot\bfzeta(0)\srd\bfx.
    \end{aligned}
\end{equation*}
Because $\bff_\eps=\nabla_{\bfd}F_\eps$ and $\mdiv\bfzeta=0$ with zero boundary trace,
\[
 \int_\Omega (\bfzeta\cdot\nabla)\bfd\cdot\bff_\eps(\bfd) \srd \bfx
 =\int_\Omega\bfzeta\cdot\nabla F_\eps(\bfd) \srd \bfx=0.
\]
Since $\bfM+n\varphi(\bfd-\bb)=\Delta\bfd-\bff_\eps(\bfd)$, the last term on the left of the preceding identity therefore reduces to
$\int_0^T\int_\Omega(\bfzeta\cdot\nabla)\bfd\cdot\Delta\bfd\srd\bfx\srd t$. Moving it to the right-hand side gives the elastic term with the minus sign in \eqref{eq:penalized fluid}.
Because $\bigcup_{m\ge1}\cV_m$ is dense in $X_u=\cV(\Omega)\cap\rH^3(\Omega)$ and the spectral projections satisfy \eqref{eq:velocity spectral projection}, the preceding identity is continuous with respect to the norm
 $\|\bfzeta\|_{\rH^1(0,T;\rL^2(\Omega))}+
 \|\bfzeta\|_{\rL^2(0,T;X_u)}.$
Time mollification followed by the $\rL^2$-orthogonal projection $\cP_k$ therefore converges in this norm (the projection is contractive both in $\rL^2(\Omega)$ and in $X_u$) and extends the identity to every test function $\bfzeta$ admitted in \autoref{def:weak sol penalized}. This proves \eqref{eq:penalized fluid}.

For the director equation, let
$\bfpsi\in\rH^1(0,T;\rL^2(\Omega;\R^3))\cap\rL^2(0,T;\rH_0^1(\Omega;\R^3))$ with $\bfpsi(T)=0$.  The Galerkin director equation reads
\begin{equation}\label{eq:k director integrated}
\begin{aligned}
 &\int_0^T\int_\Omega\Bigl[
 -\bfd^k\cdot\partial_t\bfpsi
 +(\bfu^k\cdot\nabla)\bfd^k\cdot\bfpsi
 +\nabla\bfd^k:\nabla\bfpsi
 +\bff_\eps(\bfd^k)\cdot\bfpsi\\
 &\hspace{4.3cm}
 +n\varphi^k(\bfd^k-\bb^k)\cdot\bfpsi
 \Bigr]\srd\bfx\srd t
 =\int_\Omega\bfd_0\cdot\bfpsi(0)\srd\bfx.
\end{aligned}
\end{equation}
 By \eqref{eq:conv dk}, we can pass the limit for the first and diffusion terms. Since $\eps>0$ is fixed and $|\bfd^k|\le1$, the map $\bff_\eps$ is Lipschitz on the range of the sequence with a constant bounded by $C\eps^{-2}$, uniformly with respect to~$k$. 
Hence \eqref{eq:strong conv dk} implies that $\bff_\eps(\bfd^k)\to\bff_\eps(\bfd)$ strongly in $\rL^2((0,T)\times\Omega).$
Similarly, by \eqref{eq:phi strong conv}, \eqref{eq:b strong conv}, and \eqref{eq:strong conv dk} we obtain $\varphi^k(\bfd^k-\bb^k)\to
 \varphi(\bfd-\bb)$ strongly in $\rL^2((0,T)\times\Omega).$
For the transport term, choose any $s>3/2$ in \eqref{eq:strong conv uk}, then we have $\bfu^k\to\bfu$ strongly in $\rL^2_t\rL^\infty_{\bfx}$ and $\nabla\bfd^k$ is bounded in $\rL^\infty_t\rL^2_{\bfx}$. Hence, $  (\bfu^k\cdot\nabla)\bfd^k
  \quad\text{is bounded in }\rL^2(0,T;\rL^2(\Omega)).$
For an arbitrary $\bfpsi\in \rL^2(0,T;\rH_0^1(\Omega))$, write the difference of the transport pairings as
\[
 \int_0^T\int_{\Omega} ((\bfu^k-\bfu)\cdot\nabla)\bfd^k\cdot\bfpsi \srd \bfx \srd t
 +\int_0^T\int_{\Omega} (\bfu\cdot\nabla)(\bfd^k-\bfd)\cdot\bfpsi \srd \bfx \srd t.
\]
The first term tends to zero by the strong $\rL^2_t\rL^\infty_{\bfx}$-convergence of
$\bfu^k$, the uniform $\rL^\infty_t\rL^2_{\bfx}$-bound for~$\nabla\bfd^k$, and
$\bfpsi\in \rL^2_t\rL^2_{\bfx}$.  For the second term, $\bfu\,\bfpsi\in
\rL^1_t\rL^2_{\bfx}$, so the weak-* convergence of $\nabla\bfd^k$ in $\rL^\infty_t\rL^2_{\bfx}$ implies the desired limit.  Thus the transport term passes directly
for the full stated energy test space. Passing to the limit in
\eqref{eq:k director integrated} we yield \eqref{eq:penalized d}.

Finally, the constitutive identity
\[
 \bfM^{n,k}=\Delta\bfd^k-\bff_\eps(\bfd^k)-n\varphi^k(\bfd^k-\bb^k)
\]
passes in distributions, using the weak $\rL^2_t\rH^2_{\bfx}$-convergence of $\bfd^k$ and the two strong convergences just proved. Hence
\[
 \bfM=\Delta\bfd-\bff_\eps(\bfd)-n\varphi(\bfd-\bb)
 \quad\text{a.e.\ in }(0,T)\times\Omega,
\]
which also identifies the molecular field appearing in the limiting momentum equation.

It remains to pass to the energy inequality and to preserve separation from the wall.  From
\eqref{eq:strong conv uk} and \eqref{eq:strong conv dk}, after one subsequence
\[
 \bfu^k(t)\to\bfu(t)\quad\text{strongly in }\rL^2(\Omega),
 \qquad
 \bfd^k(t)\to\bfd(t)\quad\text{strongly in }\rH^1(\Omega)
\]
for a.e.\ $t$.  At those times the kinetic, director-gradient, Ginzburg--Landau, and director-penalty
parts of the energy pass by strong convergence (using also \eqref{eq:phi strong conv}--\eqref{eq:b strong conv}).
The viscous, artificial-viscosity, and molecular-field dissipations pass by weak lower semicontinuity.
Taking the liminf in \eqref{eq:k-energy-identity} therefore gives \eqref{eq:n-energy-inequality} for a.e.\ $t$.
Finally, the separation estimate \eqref{eq:lipschitz flow separation} holds for every $k$ by \autoref{prop:k-weak sols} and passes to the limit by the uniform convergence \eqref{eq:flow convergence k} of the flows. This completes the proof.
\end{proof}

\subsection{Rigidity limit \texorpdfstring{$n\to\infty$}{n to infinity}}\label{ssec:limit n to infinity}

We pass to the rigidity limit $n\to\infty$ with $\delta$ and $\eps$ fixed. The director is compact on every spacetime region that remains a positive distance from the limiting solid. For the momentum equation this local compactness is sufficient. Admissible test fields are rigid in a neighborhood of the body, and the symmetric Ericksen tensor contracts only with the symmetric gradient of the test field, hence contributes nothing in that neighborhood. The velocity penalization gives $\varphi\bD(\bfu)=0$ in the limit, from which the limiting transported solid is a genuine rigid motion.

Fix $T>0$. By \autoref{thm:ex of sol to penalized probl}, the penalized problem has a solution $(\bfu^n,\bfd^n,\varphi^n,\bb^n)$ on $[0,T]$ for every $n$, and, by \eqref{eq:lipschitz flow separation}, all approximate solids satisfy $\dist(\cS(\varphi^n(t)),\partial\Omega)\ge\rho_\delta(T)>0$ for $0\le t\le T$. The energy inequality~\eqref{eq:n-energy-inequality} provides bounds on $\bfu^n$ and $\bfd^n$ that are uniform in $n$.
Applying the Banach--Alaoglu theorem, we can extract a subsequence (not relabeled) such that
\begin{align}
    \bfu^n &\rightarrow \bfu \quad\text{weakly in } \rL^2(0,T;\rH^3(\Omega;\R^3)) \;\text{and weak-* in } \rL^\infty(0,T;\cH(\Omega)),\label{con:un}\\
    \bfd^n & \rightarrow\bfd \quad\text{weak-* in } \rL^\infty(0,T;\rH^1(\Omega;\R^3)),\label{con:dn}\\
    \bfM^n
    &\rightarrow \bfM \text{ weakly in }\rL^2(0,T; \rL^2(\Omega;\R^3)), \label{eq:conv Mn}\\
    \varphi^n & \rightarrow \varphi \quad \text{weak-* in } \rL^\infty(0,T;\rL^\infty(\Omega)),\label{con:phin}
\end{align}
for some limit functions
\begin{gather*}
    \bfu\in\rL^2(0,T;\rH^3(\Omega;\R^3))\cap\rL^\infty(0,T;\cH(\Omega)),\\
    \bfd\in\rL^\infty(0,T;\rH^1(\Omega;\R^3)),\qquad \bfM\in \rL^2((0,T)\times\Omega;\R^3),\qquad 0 \le \varphi \le 1.
\end{gather*}

We first show that the flows of the approximate velocities converge, identify the limits of the transported indicator and lifting, and deduce that the limiting solid moves rigidly.
For fixed $\delta$, by the energy estimate and the embedding $\rH^3(\Omega)\hookrightarrow\rW^{1,\infty}(\Omega)\cap \rW^{2,6}(\Omega)$, one has uniform bi-Lipschitz flow bounds. Differentiating the flow equation twice in space, we have
\[
 \frac{\rd}{\rd t}\nabla^2\bfX^n_{0,t}
 = (\nabla^2\bfu^n)\circ\bfX^n_{0,t}
      [\nabla\bfX^n_{0,t},\nabla\bfX^n_{0,t}]
   +(\nabla\bfu^n)\circ\bfX^n_{0,t}\,\nabla^2\bfX^n_{0,t},
\]
and using volume preservation and Gr\"onwall's inequality one obtains
\[
 \sup_{n}\sup_{0\le t\le T}
 \left(
  \|\bfX^n_{0,t}\|_{\rW^{2,6}(\Omega)}
  +\|\bfY^n(t)\|_{\rW^{2,6}(\Omega)}
 \right)\le C_\delta.
\]
The inverse estimate follows from the derivative formulas for inverse maps. The maps and their inverses are also uniformly $1/2$-H\"older in time with values in $\rC(\overline\Omega)$. By the Arzel\`a--Ascoli theorem, we get, up to a subsequence,
\begin{equation}\label{eq:compactness flows n}
    \bfX^n_{0,t}\to\bfX_{0,t},
    \qquad
    \bfY^n(t,\cdot)\to\bfY(t,\cdot)\quad\text{uniformly on }[0,T]\times\overline\Omega.
\end{equation}
These are the flows of the weak limit $\bfu$. Indeed, for fixed $t$ and $\bfy$, the error in the integral flow equation splits into
\[
\int_0^t \bigl[\bfu^n(s,\bfX_{0,s}^n(\bfy)) -\bfu^n(s,\bfX_{0,s}(\bfy)) \bigr] \srd s +\int_0^t (\bfu^n -\bfu)(s, \bfX_{0,s}(\bfy))\srd s.
\]
The first term tends to $0$ by the uniform convergence of the maps and the $\rL^1_t\rW^{1,\infty}_{\bfx}$-bound. The second term tends to $0$ by weak $\rL_t^2\rH^3_{\bfx}$-convergence, since evaluation along the
fixed continuous path is a bounded element of
$\rL^2(0,T;(\rH^3(\Omega))')$. Uniqueness of the Lipschitz flow
identifies the limits. They preserve volume and fix $\partial\Omega$.

The indicator and lifting convergence now follows exactly as in
\eqref{eq:phi strong conv} and \eqref{eq:b strong conv}, namely
\begin{equation}\label{eq:phin strong before compactness}
 \begin{aligned}
 \varphi^n&\to\varphi:=\varphi_0\circ\bfY
 \quad\text{in }\rC([0,T];\rL^p(\Omega)),\quad p<\infty,\\
 \bb^n&\to\bb:=\bfd_0\circ\bfY
 \quad\text{in }\rC([0,T];\rL^\infty(\Omega)).
 \end{aligned}
\end{equation}
Here the indicator estimate uses $|\partial\cS_0|=0$,
and the lifting estimate uses the H\"older continuity of~$\bfd_0\in\rH^2(\Omega)$.
In particular, with $\cS(t):=\bfX_{0,t}(\cS_0)$,
\begin{equation}\label{eq:conv of Hausdorff dist n}
 \sup_{t\in[0,T]}d_{\rm H}(\cS(\varphi^n(t)),\cS(t))\to0.
\end{equation}
Thus both transport identities pass to the limit, and the following
inclusions hold.
\begin{lem}\label{lem:rigid-set Hausdorff inclusions}
For every $\sigma>0$ and all sufficiently large $n$,
\[
\cS(\varphi^n(t))\subset\cS(t)_\sigma,\qquad \cS(t)\subset\cS(\varphi^n(t))_\sigma
 \quad(0\le t\le T).
\]
\end{lem}
The two penalty bounds give
\[
\|\varphi^n\bD(\bfu^n)\|_{\rL^2_{t,\bfx}}\le Cn^{-1/2},
 \qquad
 \|\varphi^n(\bfd^n-\bb^n)\|_{\rL^\infty_t\rL^2_{\bfx}}
 \le Cn^{-1/2}.
\]
Combining them with \eqref{eq:phin strong before compactness} and
the weak convergences yields
$\varphi\bD(\bfu)=0$ and $\varphi(\bfd-\bb)=0$.
Set $\bfh(t)=|\cS_0|^{-1}\int_{\cS(t)}\bfx\srd\bfx$.
The connected transported set therefore has the rigid velocity~${\bfell(t)+\bfomega(t)\times(\bfx-\bfh(t))}$.
Differentiating its center of mass along the flow gives
$\bfh'=\bfell$. The relative trajectories solve
$\partial_t(\bfX_{0,t}(\bfy)-\bfh)
 =\bfomega\times(\bfX_{0,t}(\bfy)-\bfh)$, whence
\[
 \bfX_{0,t}(\bfy)=\bfh(t)+\bfO(t)\bfy\quad(\bfy\in\cS_0),
 \qquad \bfO'=\bfm\bfO,\quad \bfO(0)=\I_3.
\]
Consequently $\cS(t)=\bfh(t)+\bfO(t)\cS_0$ and
\begin{equation}\label{eq:d equals b solid} \bfd=\bb=\bfd_0(\bfO(t)^\top(\bfx-\bfh(t)))
 \quad\text{a.e.\ in }\cS(t).
\end{equation}

Next, we prove that the director converges strongly away from the solid and that the velocity converges strongly in the whole container.
By \eqref{eq:d equals b solid}, the lifted director field $\bfw:=\bfd-\bb$ vanishes almost everywhere in $\cS(t)$.

The director sequence is locally strongly compact as follows.
The proof is given in \autoref{ssec:aux ests for rigidity limit}.

\begin{lem}\label{lem:strong conv d}
For every $\sigma>0$, the sequence $\{\bfd^n\}$ converges strongly to $\bfd$ in $\rL^2(0,T;\rH^1(\widetilde{\cF}_\sigma(t);\R^3))$, where
\[
   \widetilde{\cF}_\sigma(t)
   \coloneqq \Omega\setminus\overline{\cS(\varphi(t))_\sigma}
   =\{\bfx\in\Omega:\dist(\bfx,\cS(\varphi(t)))>\sigma\}.
\]
Thus $\widetilde{\cF}_\sigma(t)$ is an open set compactly separated from the solid whenever $\sigma$ is below the tubular radius.
\end{lem}

The velocity compactness is based on the following criterion.
Its proof can be found in \autoref{ssec:aux ests for rigidity limit}.

\begin{lem}\label{lem:penalization velocity compactness}
Let $S(t)=\bfh(t)+\bfO(t)S_0$ be a rigid trajectory of a connected $\rC^2$-set $S_0$.  Let~$\chi_n=\dsone_{S_n(t)}$ and $\chi=\dsone_{S(t)}$, where $S_n(t)=\Phi_n(t,S_0)$.
Assume that every $\Phi_n(t)$ is volume preserving on a fixed neighborhood of $S_0$ and that
\begin{equation}\label{eq:uniform bilipschitz penal compactness}
 \sup_{n,t}\Bigl(
 \|\Phi_n(t)\|_{\rW^{1,\infty}}
 +\|\Phi_n(t)^{-1}\|_{\rW^{1,\infty}}
 +\|\Phi_n(t)\|_{\rW^{2,6}}
 +\|\Phi_n(t)^{-1}\|_{\rW^{2,6}}
 \Bigr)<\infty.
\end{equation}
Assume also that $S_n(t)\to S(t)$ uniformly in Hausdorff distance, the sets remain uniformly separated from $\partial\Omega$, and
$\chi_n\to\chi$ strongly in $\rC([0,T];\rL^p(\Omega))$ for every finite $p$. Suppose
\[
    \bfv_n \stackrel{*}{\rightharpoonup}\bfv \enspace \text{in }\rL^\infty(0,T;\cH(\Omega)), \quad \bfv_n \rightharpoonup\bfv \enspace \text{in } \rL^2(0,T;\cV(\Omega)),
\]
as well as 
\begin{equation}\label{eq:L2H3-bound}
    \sup_n\|\bfv_n\|_{\rL^2(0,T;\rH^3(\Omega))} <\infty, \tand
\end{equation}
\begin{equation}\label{eq:penalty est}
    \begin{aligned}
        n\int_0^T\int_\Omega\chi_n|\bD(\bfv_n)|^2\srd\bfx\srd t \le C.
    \end{aligned}
\end{equation}
Assume in addition that, for every sufficiently small $\sigma>0$, $[0,T]$ admits a finite covering by intervals~$I_j$ and fixed connected $\rC^2$-sets $S_j\Subset\Omega$ satisfying
\[
 S_n(t)\subset S_j\subset S(t)_\sigma\qquad(t\in I_j,
 \ n\ge n_\sigma),
\]
and there is $r_*>0$, independent of small $\sigma$ and of $j$, such that
\[
 D_j:=(S_j)_{r_*}\setminus\overline{S_j}\Subset\Omega
\]
are $\rC^2$-collars with uniformly bounded $\rC^2$-character.  
Assume finally that on each $I_j$ there is a constant $C_{\sigma,j}$, independent of $n$, such that
\begin{equation*}
 \left|\int_{I_j}\int_\Omega \bfv_n\cdot\partial_t\bfzeta\srd\bfx\srd t\right|
 \le C_{\sigma,j}\Bigl(
 \|\bfzeta\|_{\rL^\infty(I_j;\rH^3(\Omega))}
 +\|\bfzeta\|_{\rL^2(I_j;\rH^3(\Omega))}\Bigr)
\end{equation*}
for every smooth, compactly time-supported, divergence-free $\bfzeta$ with
$\bD(\bfzeta)=0$ on $S_j$. Then
\[
 \bfv_n\to \bfv\qquad\text{strongly in }\rL^2(0,T;\rL^2(\Omega)).
\]
\end{lem}

The application of the preceding criterion to the present sequence is stated next. Its proof is given in \autoref{ssec:aux ests for rigidity limit}.

\begin{lem}\label{lem:strong conv u}
The sequence $\{\bfu^n\}$ in~\eqref{con:un} converges strongly to $\bfu$ in $\rL^2(0,T;\rL^2(\Omega;\R^3))$.
\end{lem}

Observe that the strong convergence of $\bfu^n \to \bfu$ in $\rL^2(0,T;\rL^2(\Omega;\R^3))$ together with the uniform boundedness of $\bfu^n$ in $\rL^2(0,T;\rH^3(\Omega;\R^3))$ and interpolation implies that $\bfu^n \to \bfu$ strongly in~$\rL^2(0,T;\rH^s(\Omega;\R^3))$ for all $s < 3$.
In particular, choosing $s \in (\frac{5}{2},3)$, we find that
\[
    \bfu^n \to \bfu \quad \text{strongly in } \rL^1(0,T;\rW^{1,\infty}(\Omega;\R^3)).
\]
Consequently, the limits $\widehat{\bfX}$ and $\widehat{\bfY}$ obtained in \eqref{eq:compactness flows n} indeed coincide with the flow map $\bfX$ and inverse flow map $\bfY$ associated with the limit fluid velocity field $\bfu$, i.e., $\widehat{\bfX} = \bfX$ and $\widehat{\bfY} = \bfY$.

We now pass to the limit in the director equation, first for test functions supported in the fluid region and then, by density, for the full class of energy tests.
Recall the fluid spacetime $\cQ_{\cF}$ from \autoref{ssec:main result}, now associated with the limiting rigid trajectory. Let $\bfpsi\in \rC_c^\infty([0,T)\times\Omega;\R^3)$ with~$\mathrm{supp}\,\bfpsi\subset\cQ_{\cF}\cup(\{0\}\times\cF_0)$.
By compactness of the support and continuity of the rigid body motion, there
exists $\sigma>0$ such that
\begin{equation*}
    \bfpsi(t,\cdot)=0
    \quad\text{on }\cS(t)_\sigma
    \quad\tforall t\in[0,T].
\end{equation*}
From \autoref{lem:rigid-set Hausdorff inclusions}, for $n$ sufficiently large, we deduce that $\cS(\varphi^n(t)) \subset \cS(\varphi(t))_{\sigma/2}$ for all $t \in [0,T]$ and therefore $\varphi^n\bfpsi=0$ a.e.\ in $(0,T)\times\Omega$.
Using $\bfpsi$ in \eqref{eq:penalized d}, the director penalization term $n \varphi^n (\bfd^n - \bb^n) \cdot \bfpsi$ vanishes identically for all sufficiently large $n$, and we obtain
\begin{equation}\label{eq:director n localized}
    \begin{aligned}
    \int_0^T\int_\Omega
    \Bigl[
        -\bfd^n\cdot\del_t\bfpsi
        +(\bfu^n\cdot\nabla)\bfd^n\cdot\bfpsi
        +\nabla\bfd^n:\nabla\bfpsi
        +\bff_\eps(\bfd^n)\cdot\bfpsi
    \Bigr]\srd\bfx\srd t
    =
    \int_{\cF_0}\bfd_0\cdot\bfpsi(0)\srd\bfx.
    \end{aligned}
\end{equation}
Since $\bfpsi$ vanishes on $\cS(\varphi(t))_\sigma$, its support is contained in
$\{\bfx:\dist(\bfx,\cS(\varphi(t)))\ge\sigma\}$ and hence is compactly contained in the open set $\widetilde{\cF}_{\sigma/2}(t)$. The local compactness result, \autoref{lem:strong conv d}, therefore gives $\bfd^n\to\bfd$ strongly in $\rL^2(0,T;\rH^1(\widetilde{\cF}_{\sigma/2}(t);\R^3))$.
Together with the strong convergence of $\bfu^n$ from \autoref{lem:strong conv u}, this allows us to pass to
the limit in \eqref{eq:director n localized}.
By the properties of the support of $\bfpsi$, we find
\begin{equation}\label{eq:director limit localized}
    \begin{aligned}
    \int_0^T\int_{\cF(t)}
    \Bigl[
        -\bfd\cdot\del_t\bfpsi
        +(\bfu\cdot\nabla)\bfd\cdot\bfpsi
        +\nabla\bfd:\nabla\bfpsi
        +\bff_\eps(\bfd)\cdot\bfpsi
    \Bigr]\srd\bfx\srd t
    =
    \int_{\cF_0}\bfd_0\cdot\bfpsi(0)\srd\bfx.
    \end{aligned}
\end{equation}
Consequently,
\[
    \del_t\bfd+(\bfu\cdot\nabla)\bfd
    =
    \Delta\bfd-\bff_\eps(\bfd)
    \quad\text{in }\cD'(\cQ_{\cF}).
\]
To identify the weak limit $\bfM$ on the fluid region, let $\bfxi\in \rC_c^\infty(\cQ_{\cF};\R^3)$.
For all sufficiently large $n$, one has
$\varphi^n\bfxi=0$. Hence the constitutive relation
\[
    \bfM^n
    =
    \Delta\bfd^n-\bff_\eps(\bfd^n)
    -n\varphi^n(\bfd^n-\bb^n)
\]
implies
\[
    \int_0^T\int_\Omega
    \left[
        \bfM^n\cdot\bfxi
        +\nabla\bfd^n:\nabla\bfxi
        +\bff_\eps(\bfd^n)\cdot\bfxi
    \right]\srd\bfx\srd t
    =0.
\]
Using the weak convergence of $\bfM^n$ from \eqref{eq:conv Mn}, the weak convergence of
$\nabla\bfd^n$ from \eqref{con:dn}, and the strong local convergence of $\bfd^n$ as established in \autoref{lem:strong conv d}, we pass to
the limit and obtain
\[
    \int_0^T\int_\Omega
    \left[
        \bfM\cdot\bfxi
        +\nabla\bfd:\nabla\bfxi
        +\bff_\eps(\bfd)\cdot\bfxi
    \right]\srd\bfx\srd t
    =0.
\]
Therefore,
\[
    \bfM=\Delta\bfd-\bff_\eps(\bfd)
    \quad\text{in }\cD'(\cQ_{\cF}).
\]
Since both sides belong locally to $\rL^2(\cQ_{\cF})$, the equality holds
almost everywhere in $\cQ_{\cF}$.

We extend \eqref{eq:director limit localized} to the moving energy test class $\cX_{\bfd}(\varphi)$ from \eqref{eq:director energy tests}. Here and below, a function defined on $\cF(t)$ is identified with its extension by zero to $\cS(t)$. We use the following density result, proved in \autoref{ssec:aux ests for rigidity limit}.

\begin{lem}
\label{lem:density director tests}
Let $\cD_{\cF}:=
    \left\{
    \bfpsi\in\rC_c^\infty([0,T)\times\Omega;\R^3):
    \supp\bfpsi
    \subset
    \cQ_{\cF}\cup(\{0\}\times\cF_0)
    \right\}$.
Then $\cD_{\cF}$ is dense in the space $\cX_{\bfd}(\varphi)$ defined in \eqref{eq:director energy tests} with respect to the norm
\[
    \|\bfpsi\|_{\rH^1(0,T;\rL^2(\Omega))}
    +
    \|\bfpsi\|_{\rL^2(0,T;\rH^1(\Omega))}.
\]
Moreover, the approximating sequence may be chosen such that $\bfpsi_m(0)\to\bfpsi(0)$ in $\rL^2(\Omega)$.
\end{lem}

By \autoref{lem:density director tests}, every $\bfpsi\in\cX_{\bfd}(\varphi)$ is the limit, in the norm of that lemma, of a sequence $(\bfpsi_m)\subset\cD_{\cF}$ with $\bfpsi_m(0)\to\bfpsi(0)$ in $\rL^2(\Omega)$, and \eqref{eq:director limit localized} holds for each $\bfpsi_m$. Every term in \eqref{eq:director limit localized} is continuous with respect to this convergence, because
\begin{gather*}
    \bfd\in\rL^\infty(0,T;\rH^1(\Omega;\R^3)),\qquad \bfu\in\rL^2(0,T;\rL^\infty(\Omega;\R^3)),\\
    \bff_\eps(\bfd)\in\rL^\infty((0,T)\times\Omega;\R^3),
\end{gather*}
where the last property follows from \autoref{lem:relaxed dir cond} and the fixed value of $\eps$. Consequently, for every~$\bfpsi\in\cX_{\bfd}(\varphi)$,
\begin{equation*}
    \int_0^T\int_{\cF(t)}
    \Bigl[
        -\bfd\cdot\del_t\bfpsi
        +(\bfu\cdot\nabla)\bfd\cdot\bfpsi
        +\nabla\bfd:\nabla\bfpsi
        +\bff_\eps(\bfd)\cdot\bfpsi
    \Bigr]\srd\bfx\srd t
    =
    \int_{\cF_0}
    \bfd_0\cdot\bfpsi(0)\srd\bfx.
\end{equation*}

We turn to the momentum equation, in which the local compactness of the director suffices because the test fields are rigid near the solid.
Let $\bfzeta\in\rC_c^1([0,T);\cV(\Omega)\cap\rH^3(\Omega;\R^3))$ be rigid in a collar of the limiting solid, i.e., for some $\sigma>0$,
\begin{equation}\label{eq:zeta rigid neighborhood}
    \bD(\bfzeta)=0
    \quad\text{a.e.\ on }\cS(t)_\sigma
    \quad\tforall t\in[0,T].
\end{equation}
By the convergence of the Hausdorff distance as established in \eqref{eq:conv of Hausdorff dist n}, for all sufficiently large
$n$, we get
\[
    \cS(\varphi^n(t))\subset\cS(t)_{\sigma/2}
    \quad\tforall t\in[0,T].
\]
In particular, $\varphi^n\bD(\bfzeta)=0$, and hence
\[
    n\int_0^T\int_\Omega
    \varphi^n\bD(\bfu^n):\bD(\bfzeta)\srd\bfx\srd t=0.
\]
For the director contribution, the identity $\bfM^n+n\varphi^n(\bfd^n-\bb^n)
    =
    \Delta\bfd^n-\bff_\eps(\bfd^n)$
gives
\[
\begin{aligned}
&-\int_\Omega
(\bfzeta\cdot\nabla)\bfd^n\cdot
\bigl(\bfM^n+n\varphi^n(\bfd^n-\bb^n)\bigr)\srd\bfx
=
-\int_\Omega
(\bfzeta\cdot\nabla)\bfd^n\cdot
\bigl(\Delta\bfd^n-\bff_\eps(\bfd^n)\bigr)\srd\bfx.
\end{aligned}
\]
Since $\bff_\eps=\nabla_{\bfd}F_\eps$, $\mdiv\bfzeta=0$, and
$\bfzeta=0$ on $\del\Omega$,
\[
    \int_\Omega
    (\bfzeta\cdot\nabla)\bfd^n\cdot\bff_\eps(\bfd^n)\srd\bfx
    =
    \int_\Omega\bfzeta\cdot\nabla F_\eps(\bfd^n)\srd\bfx
    =0.
\]
Integration by parts therefore yields
\begin{equation*}
\begin{aligned}
&-\int_\Omega
(\bfzeta\cdot\nabla)\bfd^n\cdot
\bigl(\bfM^n+n\varphi^n(\bfd^n-\bb^n)\bigr)\srd\bfx\\
&\qquad=
\int_\Omega
(\nabla\bfd^n)^\top\nabla\bfd^n:\nabla\bfzeta\srd\bfx
=
\int_\Omega
(\nabla\bfd^n)^\top\nabla\bfd^n:\bD(\bfzeta)\srd\bfx.
\end{aligned}
\end{equation*}
By \eqref{eq:zeta rigid neighborhood}, the last integrand vanishes on
$\cS(t)_\sigma$. Hence,
\[
\begin{aligned}
\int_0^T\int_\Omega
(\nabla\bfd^n)^\top\nabla\bfd^n:\nabla\bfzeta
\srd\bfx\srd t
=
\int_0^T\int_{\widetilde{\cF}_\sigma(t)}
(\nabla\bfd^n)^\top\nabla\bfd^n:\nabla\bfzeta
\srd\bfx\srd t.
\end{aligned}
\]
On the open set $\widetilde{\cF}_\sigma(t)$, the local strong compactness established in \autoref{lem:strong conv d} gives
\[
    \nabla\bfd^n\to\nabla\bfd
    \quad\text{strongly in }\rL^2(0,T;\rL^2(\widetilde{\cF}_\sigma(t);\R^{3\times 3})).
\]
Consequently, $(\nabla\bfd^n)^\top\nabla\bfd^n \to (\nabla\bfd)^\top\nabla\bfd$ strongly in $\rL^1(0,T;\rL^1(\widetilde{\cF}_\sigma(t);\R^{3\times3}))$ on this region, and therefore
\begin{equation*}
\begin{aligned}
&-\int_0^T\int_\Omega
(\bfzeta\cdot\nabla)\bfd^n\cdot
\bigl(\bfM^n+n\varphi^n(\bfd^n-\bb^n)\bigr)
\srd\bfx\srd t
\\
&\qquad\rightarrow
\int_0^T\int_{\widetilde{\cF}_\sigma(t)}
(\nabla\bfd)^\top\nabla\bfd:\nabla\bfzeta
\srd\bfx\srd t=\int_0^T\int_{\cF(t)}
(\nabla\bfd)^\top\nabla\bfd:\nabla\bfzeta
\srd\bfx\srd t,
\end{aligned}
\end{equation*}
where we additionally used that $\bD(\bfzeta) = 0$ on $\cS(t)_\sigma \supset \cF(t) \cap \cS(t)_\sigma$.

The remaining terms in \eqref{eq:penalized fluid} pass to the limit by \autoref{lem:strong conv u} and \eqref{con:un}. Hence
\begin{equation*}
\begin{aligned}
 &\int_0^T\int_\Omega\Bigl[-\bfu\cdot\partial_t\bfzeta-(\bfu\otimes\bfu):\nabla\bfzeta
 +2\bD(\bfu):\bD(\bfzeta)\Bigr]\srd\bfx\srd t+\delta\int_0^T a_u(\bfu,\bfzeta)\srd t\\
 &\qquad=\int_\Omega\bfu_0\cdot\bfzeta(0)\srd\bfx
 +\int_0^T\int_{\cF(t)}(\nabla\bfd)^\top\nabla\bfd:\nabla\bfzeta\srd\bfx\srd t
\end{aligned}
\end{equation*}
for every $\bfzeta\in\rC_c^1([0,T);\cV(\Omega)\cap\rH^3(\Omega;\R^3))$ that is rigid in a collar of the limiting solid.

It remains to derive the energy inequality for the limit and to show that the limiting solid stays away from the container wall.
To this end, we pass to the limit in the approximate energy inequality~\eqref{eq:n-energy-inequality}. Besides the weak convergences established
above, we use the global compactness
\begin{equation}\label{eq:strong conv dn in CL2}
    \bfd^n\to\bfd
    \quad\text{strongly in }
    \rC([0,T];\rL^2(\Omega;\R^3)).
\end{equation}
Indeed, by the director equation,
\[
    \del_t\bfd^n
    =
    \bfM^n-(\bfu^n\cdot\nabla)\bfd^n,
\]
and the right-hand side is uniformly bounded in
$\rL^2(0,T;\rL^2(\Omega;\R^3))$. Together with the uniform bound in
$\rL^\infty(0,T;\rH^1(\Omega;\R^3))$, this yields the asserted compactness.
By interpolation, $\bfd^n\to\bfd$ strongly in $\rC([0,T];\rL^q(\Omega;\R^3))$ for all $q\in[2,6)$.
In particular,
\[
    F_\eps(\bfd^n)\to F_\eps(\bfd)
    \quad\text{strongly in }
    \rC([0,T];\rL^1(\Omega)).
\]

The strong velocity compactness from \autoref{lem:strong conv u} and the strong
$\rC([0,T];\rL^2(\Omega;\R^3))$-convergence of the directors from \eqref{eq:strong conv dn in CL2} allow one common subsequence such that
$\bfu^n(t)\to\bfu(t)$ and $\bfd^n(t)\to\bfd(t)$ strongly in $\rL^2(\Omega)$ for a.e.\ $t$.
At any such time, the uniform $\rH^1$-bound on $\bfd^n(t)$ identifies every weak~$\rH^1$-cluster point with $\bfd(t)$, so the director gradient is lower semicontinuous.
The kinetic energy converges strongly, and the potential converges by the preceding $\rC_t\rL^1_{\bfx}$-convergence.  
The viscous and
artificial-viscosity dissipations and $\bfM^n\rightharpoonup\bfM$ pass by weak lower semicontinuity. 
We obtain, initially for almost every $t\in(0,T)$,
\[
\begin{aligned}
&\frac12\int_\Omega|\bfu(t)|^2\srd\bfx
+\int_\Omega
\left(
    \frac12|\nabla\bfd(t)|^2+F_\eps(\bfd(t))
\right)\srd\bfx
+2\int_0^t\int_\Omega|\bD(\bfu)|^2\srd\bfx\srd s
+\delta\int_0^t\|\bfu(s)\|_{X_u}^2\srd s\\
&+\int_0^t\int_{\cF(s)}|\bfM|^2\srd\bfx\srd s
\leq
\frac12\int_\Omega|\bfu_0|^2\srd\bfx
+\int_\Omega
\left(
    \frac12|\nabla\bfd_0|^2
    +F_\eps(\bfd_0)
\right)\srd\bfx,
\end{aligned}
\]
where we used the identity $\bfM = \Delta\bfd+\frac1{\eps^2}(1-|\bfd|^2)\bfd$ a.e.\ in $\cQ_{\cF}$.

To remove the director energy contained in the solid, set
\[
    C_{\cS}
    :=
    \int_{\cS_0}
    \left(
        \frac12|\nabla\bfd_0|^2+F_\eps(\bfd_0)
    \right)\srd\bfx
    =
    \frac12\int_{\cS_0}|\nabla\bfd_0|^2\srd\bfx.
\]
Since $|\bfd_0|=1$ a.e.\ in $\Omega$, one has $F_\eps(\bfd_0)=0$. The initial director energy therefore splits as
\[
    \int_\Omega
    \left(
        \frac12|\nabla\bfd_0|^2
        +F_\eps(\bfd_0)
    \right)\srd\bfx
    =
    \int_{\cF_0}
    \left(
        \frac12|\nabla\bfd_0|^2
        +F_\eps(\bfd_0)
    \right)\srd\bfx
    +C_{\cS}.
\]
Furthermore, $\bfd=\bb$ in $\cS(t)$, and the restriction of the limiting
flow to the rigid body is a rigid transformation. Since $\bb(t,\bfx)=\bfd_0(\bfY(t,\bfx))$ and $|\bfd_0|=1$, we obtain
\[
    \int_{\cS(t)}
    \left(
        \frac12|\nabla\bfd(t)|^2+F_\eps(\bfd(t))
    \right)\srd\bfx
    =
    C_{\cS}.
\]
Subtracting $C_{\cS}$ from both sides gives
\begin{equation}\label{eq:energy est after n lim}
\begin{aligned}
E^{\delta,\eps}(t)
&+
2\int_0^t\int_\Omega|\bD(\bfu)|^2\srd\bfx\srd s
+
\delta\int_0^t\|\bfu(s)\|_{X_u}^2\srd s
\\
&+
\int_0^t\int_{\cF(s)}
\left|
    \Delta\bfd
    +\frac1{\eps^2}(1-|\bfd|^2)\bfd
\right|^2
\srd\bfx\srd s
\leq
E^{\delta,\eps}(0)=E_0,
\end{aligned}
\end{equation}
where
\[
    E^{\delta,\eps}(t)
    :=
    \frac12\int_\Omega|\bfu(t)|^2\srd\bfx
    +
    \int_{\cF(t)}
    \left(
        \frac12|\nabla\bfd(t)|^2
        +F_\eps(\bfd(t))
    \right)\srd\bfx,
\]
and, because $|\bfd_0|=1$, the initial value $E_0$ is independent of both $\delta$ and $\eps$.

Finally, \eqref{eq:lipschitz flow separation} and the uniform convergence of the forward flow maps give
\[
 \dist(\cS(\varphi(t)),\partial\Omega)\ge\rho_\delta(T)>0
 \qquad(0\le t\le T).
\]
Since $T>0$ was arbitrary, the rigidity limit can be carried out on every bounded time interval. At fixed $\delta$, the body never reaches the container wall.

\subsection{Vanishing artificial-viscosity limit \texorpdfstring{$\delta\to0$}{delta to zero}}\label{ssec:delta limit}
Fix $\eps>0$ and assume $E_0>0$. If $E_0=0$, the time $T_*$ below may be replaced by any finite time. For each $\delta>0$, let $(\bfu_\delta,\bfd_\delta,\bfM_\delta,\varphi_\delta)$ be the post-rigidity solution on $[0,T_*]$ provided by \autoref{ssec:limit n to infinity}. We denote
\[
 \cS_\delta(t)=\cS(\varphi_\delta(t)),\qquad
 \cF_\delta(t)=\Omega\setminus\overline{\cS_\delta(t)}.
\]
In the fluid region,
\begin{equation*}
 \bfM_\delta:=\Delta\bfd_\delta-\bff_\eps(\bfd_\delta)
 =\Delta\bfd_\delta+\frac1{\eps^2}(1-|\bfd_\delta|^2)\bfd_\delta.
\end{equation*}

The post-rigidity solution satisfies
\begin{align*}
 &\bfu_\delta\in \rL^\infty(0,T_*;\cH(\Omega))\cap\rL^2(0,T_*;\cK(\varphi_\delta)\cap \rH^3(\Omega;\R^3)),\\
 &\bfd_\delta|_{\cF_\delta(t)}\in \rL^\infty(0,T_*;\rH^1(\cF_\delta(t);\R^3)),\qquad
 \bfM_\delta\in \rL^2(\cQ_{\cF_\delta};\R^3),\\
 &\varphi_\delta\in \Char([0,T_*]\times\Omega)\cap \rC^{0,1/p}([0,T_*];\rL^p(\Omega))\quad(1\le p<\infty).
\end{align*}
Its momentum equation reads
\begin{equation}\label{eq:fluid delta}
\begin{aligned}
 &\int_0^{T_*}\int_\Omega\Bigl[-\bfu_\delta\cdot\partial_t\bfzeta
 -(\bfu_\delta\otimes\bfu_\delta):\nabla\bfzeta
 +2\bD(\bfu_\delta):\bD(\bfzeta)\Bigr] \srd \bfx \srd t\\
 &\qquad+\delta\int_0^{T_*}a_u(\bfu_\delta,\bfzeta)\srd t
 =\int_\Omega\bfu_0\cdot\bfzeta(0) \srd \bfx
 +\int_0^{T_*}\int_{\cF_\delta(t)}(\nabla\bfd_\delta)^\top\nabla\bfd_\delta:\nabla\bfzeta \srd \bfx \srd t,
\end{aligned}
\end{equation}
for every
\[
 \bfzeta\in \rH^1(0,T_*;\rH^1(\Omega;\R^3))\cap
 \rL^2(0,T_*;\cV(\Omega)\cap \rH^3(\Omega;\R^3)),\qquad \bfzeta(T_*)=0,
\]
for which there exists $r>0$ such that
$\bD(\bfzeta)=0$ on $(\cS_\delta(t))_r$ for a.e.\ $t\in(0,T_*)$. Its director equation reads
\begin{equation*}
 \int_0^{T_*}\int_{\cF_\delta(t)}
 \Bigl[-\bfd_\delta\cdot\partial_t\bfpsi
 +(\bfu_\delta\cdot\nabla)\bfd_\delta\cdot\bfpsi
 +\nabla\bfd_\delta:\nabla\bfpsi+\bff_\eps(\bfd_\delta)\cdot\bfpsi\Bigr] \srd \bfx \srd t
 =\int_{\cF_0}\bfd_0\cdot\bfpsi(0) \srd \bfx,
\end{equation*}
for every $\bfpsi\in\cX_{\bfd}(\varphi_\delta)$, see \eqref{eq:director energy tests}, and its transport equation reads
\begin{equation*}
 \int_0^{T_*}\int_\Omega\varphi_\delta
 [\partial_t\eta+\bfu_\delta\cdot\nabla\eta] \srd \bfx \srd t
 =-\int_\Omega\varphi_0\eta(0) \srd \bfx
\end{equation*}
for every $\eta\in\rC^1([0,T_*]\times\overline\Omega)$ with $\eta(T_*)=0$.
The energy inequality \eqref{eq:energy est after n lim} gives, with $C$ independent of $\delta$,
\begin{equation*}
 \|\bfu_\delta\|_{\rL^\infty_t\rL^2_{\bfx}}
 +\|\bfu_\delta\|_{\rL^2_t\rH^1_{0,\bfx}}
 +\|\nabla\bfd_\delta\|_{\rL^\infty_t\rL^2(\cF_\delta(t))}
 +\|\bfM_\delta\|_{\rL^2(\cQ_{\cF_\delta})}
 +\sqrt\delta\,\|\bfu_\delta\|_{\rL^2(0,T_*;X_u)}\le C.
\end{equation*}
The maximum principle of \autoref{lem:relaxed dir cond}, applied at the approximation level, also gives
\begin{equation}\label{eq:delta director convex constraint}
 |\bfd_\delta|\le1,\qquad (\bfd_\delta)_3\ge0
 \quad\text{a.e.\ in }\cQ_{\cF_\delta}.
\end{equation}

We first show that the solids stay uniformly away from the container wall and that their rigid motions are compact.
Since $\bfu_\delta$ is rigid in $\cS_\delta(t)$, there are~$\bfh_\delta$,~$\bfell_\delta=\bfh_\delta'$,~$\bfomega_\delta$, and~$\bfO_\delta\in \SO(3)$ such that
\[
 \cS_\delta(t)=\bfh_\delta(t)+\bfO_\delta(t)\cS_0,
 \qquad
 \bfu_\delta=\bfell_\delta+\bfomega_\delta\times(\bfx-\bfh_\delta)
 \quad\text{in }\cS_\delta(t).
\]
Combining the identity $\int_{\cS_\delta(t)}|\bfu_\delta|^2\srd\bfx=m|\bfell_\delta|^2+\bJ_\delta(t)\bfomega_\delta\cdot\bfomega_\delta$ and the energy inequality \eqref{eq:energy est after n lim} we have
\[
 \|\bfell_\delta\|_{\rL^\infty(0,T_*)}+\|\bfomega_\delta\|_{\rL^\infty(0,T_*)}\le C.
\]
Moreover, \eqref{eq:energy est after n lim} contains the kinetic energy and the viscous dissipation with initial energy $E_0$, so that \autoref{cor:uniform separation} yields, for every $\delta>0$,
\begin{equation}\label{eq:uniform noncontact delta}
 \dist(\cS_\delta(t),\partial\Omega)\ge\frac{\delta_0}{2}
 \qquad(0\le t\le T_*).
\end{equation}

After extraction, $\bfh_\delta\to\bfh$ and $\bfO_\delta\to\bfO$ uniformly on $[0,T_*]$.
With $\cS(t)=\bfh(t)+\bfO(t)\cS_0$ and~$\varphi=\dsone_{\cS(t)}$, this implies
\begin{equation}\label{eq:strong phi delta}
    \sup_{t\in[0,T_*]}d_{\rm H}(\cS_\delta(t),\cS(t))\to0, \qquad \varphi_\delta\to\varphi
 \quad\text{strongly in }\rC([0,T_*];\rL^p(\Omega)),\qquad 1\le p<\infty,
\end{equation}
and, for every $\sigma>0$ and all small $\delta$,
\begin{equation}\label{eq:solid inclusions delta}
 \cS_\delta(t)\subset\cS(t)_\sigma,
 \qquad \cS(t)\subset\cS_\delta(t)_\sigma
 \quad(0\le t\le T_*).
\end{equation}
After a further extraction, as $\delta \to 0$, we get
\[
 \bfell_\delta\stackrel{*}{\rightharpoonup}\bfell,
 \qquad \bfomega_\delta\stackrel{*}{\rightharpoonup}\bfomega
 \quad\text{in }\rL^\infty(0,T_*),
\]
whereas $\bfh_\delta\to\bfh$ and $\bfO_\delta\to\bfO$
uniformly as recorded above. 
Passing to the integral rigid-motion equations gives
\[
 \bfh'=\bfell,\qquad
 \bfO'=\bfm\bfO,\qquad \bfm\bfx=\bfomega\times\bfx.
\]
In addition,
\[
 \varphi_\delta\bfu_\delta
 =\varphi_\delta\bigl[
   \bfell_\delta+\bfomega_\delta\times(\bfx-\bfh_\delta)
   \bigr].
\]
The strong convergence of $\varphi_\delta$ and $\bfh_\delta$,
together with the weak velocity convergence and the weak-* convergence
of the rigid coefficients, passes this identity in distributions.
Consequently
\begin{equation}\label{eq:rigid body vel after delta lim}
    \bfu(t,\bfx)=\bfell(t)+\bfomega(t)\times(\bfx-\bfh(t))
 \quad\text{for a.e.\ }\bfx\in\cS(t).
\end{equation}
Thus the limiting rigid trajectory is generated by the rigid
restriction of the limiting velocity, including its orientation.

Next, we extend the director through the solid in order to obtain compactness in the whole container, and then strengthen the convergence locally in the fluid region.
To this end, define the whole-container material extension
\begin{equation*}
 \overline\bfd_\delta(t,\bfx):=
 \begin{cases}
   \bfd_\delta(t,\bfx),&\bfx\in\cF_\delta(t),\\[2pt]
   \bfd_0\bigl(\bfO_\delta(t)^\top(\bfx-\bfh_\delta(t))\bigr),&\bfx\in\cS_\delta(t).
 \end{cases}
\end{equation*}
The two traces agree on $\partial\cS_\delta(t)$, so $\overline\bfd_\delta(t)\in \rH^1(\Omega;\R^3)$ for a.e.\ $t$.  The solid contribution to its Dirichlet energy is invariant under the rigid motion.  Hence
\begin{equation}\label{eq:global H1 director delta}
 \sup_{\delta>0}\|\overline\bfd_\delta\|_{\rL^\infty(0,T_*;\rH^1(\Omega))}\le C.
\end{equation}
Extend $\bfM_\delta$ by zero to the solid and denote the extension by $\widetilde\bfM_\delta$.  Applying \autoref{lem:material gluing} to the fluid equation and the rigidly transported solid branch gives
\begin{equation}\label{eq:global material director delta}
 \partial_t\overline\bfd_\delta+
 (\bfu_\delta\cdot\nabla)\overline\bfd_\delta
 =\widetilde\bfM_\delta
 \quad\text{in }\cD'((0,T_*)\times\Omega).
\end{equation}
Because $\bfu_\delta\in \rL^2_t\rL^6_{\bfx}$ and
$\nabla\overline\bfd_\delta\in \rL^\infty_t\rL^2_{\bfx}$, the transport term is bounded in $\rL^2_t\rL^{3/2}_{\bfx}\hookrightarrow \rL^2_t\rH^{-1}_{\bfx}$.  
Thus, $\partial_t\overline\bfd_\delta$ is bounded in $\rL^2(0,T_*;\rH^{-1}(\Omega;\R^3))$.
By the Aubin--Lions lemma, after extracting a subsequence we have
\begin{equation}\label{eq:strong director delta L2 global}
 \overline\bfd_\delta\to\overline\bfd
 \quad\text{strongly in }\rL^2((0,T_*)\times\Omega)
\end{equation}
and weak-* in $\rL^\infty_t\rH^1_{\bfx}$.  In view of \eqref{eq:delta director convex constraint} and the unit-length solid branch, the limit satisfies $|\overline\bfd|\le1$ and $(\overline\bfd)_3\ge0$.

We strengthen the convergence on compact subsets of the fluid region.  Let $K\Subset\cQ_\cF$.  By \eqref{eq:strong phi delta}, a finite cylinder covering of $K$ provides $V_j\Subset U_j\Subset\Omega$ and time intervals $I_j$ such that $U_j\subset\cF_\delta(t)\cap\cF(t)$ for $t\in I_j$ and all sufficiently small $\delta$.  On $I_j\times U_j$,
\[
 \Delta\bfd_\delta=\bfM_\delta+\bff_\eps(\bfd_\delta).
\]
For fixed $\eps$, \eqref{eq:delta director convex constraint} makes $\bff_\eps(\bfd_\delta)$ uniformly bounded.  Interior elliptic regularity therefore gives a uniform $\rL^2(I_j;\rH^2(V_j))$-bound.  Equation \eqref{eq:global material director delta} gives the corresponding $\rL^2(I_j;\rH^{-1}(V_j))$-bound on the time derivative.  Applying the Aubin--Lions lemma, we yield
\begin{equation}\label{eq:strong director delta local}
 \bfd_\delta\to\bfd
 \quad\text{strongly in }\rL^2(I_j;\rH^1(V_j)),
\end{equation}
where $\bfd=\overline\bfd|_{\cQ_\cF}$.  Consequently, for every $K\Subset\cQ_\cF$,
\begin{equation}\label{eq:strong stress delta local}
 (\nabla\bfd_\delta)^\top\nabla\bfd_\delta
 \to(\nabla\bfd)^\top\nabla\bfd
 \quad\text{strongly in }\rL^1(K).
\end{equation}
Moreover, since $\bff_\eps$ is Lipschitz on the closed unit ball, \eqref{eq:strong director delta L2 global} gives
\begin{equation}\label{eq:strong f director delta}
 \bff_\eps(\overline\bfd_\delta)\to\bff_\eps(\overline\bfd)
 \quad\text{strongly in }\rL^2((0,T_*)\times\Omega).
\end{equation}
After extracting a subsequence,
\begin{equation}\label{eq:weak M delta}
 \widetilde\bfM_\delta\rightharpoonup\widetilde\bfM
 \quad\text{weakly in }\rL^2((0,T_*)\times\Omega).
\end{equation}
Testing on compact subsets of $\cQ_\cF$ and using \eqref{eq:strong director delta local} we get
\[
 \widetilde\bfM=\Delta\bfd-\bff_\eps(\bfd)
 \quad\text{in }\cD'(\cQ_\cF).
\]
By construction of the weak extension, $\widetilde\bfM=0$ a.e.\ in the spacetime interior of the solid. 

The anchoring condition can be passed using only local collars, so no global moving-domain
diffeomorphism is needed.  Choose $r>0$ smaller than one quarter of the uniform separation in
\eqref{eq:uniform noncontact delta} and smaller than the tubular radii of $\partial\Omega$ and
$\partial\cS_0$.  The fixed outer collar
\[
 U_{\rm out}:=\{\bfx\in\Omega:\dist(\bfx,\partial\Omega)<r\}
\]
is contained in $\cF_\delta(t)\cap\cF(t)$ for every $t$ and every small $\delta$.  By
\eqref{eq:global H1 director delta} and \eqref{eq:strong director delta L2 global}, the restrictions~$\overline\bfd_\delta|_{U_{\rm out}}$ converge weakly in $\rL^2_t\rH^1_{\bfx}$ and strongly in
$\rL^2_{t,\bfx}$ to $\overline\bfd|_{U_{\rm out}}$.  Their traces on $\partial\Omega$ are all
$\bfd_0$, so continuity of the trace operator gives
$\operatorname{Tr}_{\partial\Omega}\bfd=\bfd_0$ for a.e.\ time.

For the inner boundary, let
\[
 U_{\rm in}^0:=\{\bfy\in\cF_0:\dist(\bfy,\partial\cS_0)<r\}
\]
and use only the rigid maps
$\bfZ_\delta(t,\bfy)=\bfh_\delta(t)+\bfO_\delta(t)\bfy$ and
$\bfZ(t,\bfy)=\bfh(t)+\bfO(t)\bfy$.  The pullbacks
\[
 \widehat\bfd_\delta(t,\bfy):=
 \overline\bfd_\delta(t,\bfZ_\delta(t,\bfy)),\qquad \bfy\in U_{\rm in}^0,
\]
are uniformly bounded in $\rL^\infty_t\rH^1(U_{\rm in}^0)$ because the maps are rigid
isometries.  The strong convergence~\eqref{eq:strong director delta L2 global} together with the
uniform convergence of $\bfh_\delta$ and $\bfO_\delta$ implies
$\widehat\bfd_\delta\to\overline\bfd\circ \bfZ$ strongly in~$\rL^2((0,T_*)\times U_{\rm in}^0)$. For the varying composition, we first establish this for smooth
approximation of $\overline\bfd$ and then by $\rL^2$-approximation using preservation of volume.
Thus the weak $\rL^2_t\rH^1_{\bfy}$-limit is uniquely $\overline\bfd\circ \bfZ$.  On
$\partial\cS_0$ every pulled-back trace equals $\bfd_0$.  The fixed-domain trace
operator is continuous, so we have 
\[
 \operatorname{Tr}_{\partial\cS(t)}\bfd(t,\bfx)
 =\bfd_0(\bfO(t)^\top(\bfx-\bfh(t))), \quad \text{ for a.e. } t.
\]
The proof of the lemma below is deferred to \autoref{ssec:aux ests for rigidity limit}.

\begin{lem}\label{lem:FSI compactness criterion}
Let $S_j(t)=\bfh_j(t)+\bfO_j(t)\cS_0\Subset\Omega$ be rigid trajectories with $\bfh_j, \bfO_j \in \rC([0,T])$ satisfying
\[
 \inf_{j,t}\dist(S_j(t),\partial\Omega)>0,
 \qquad \bfh_j\to \bfh,\quad \bfO_j\to \bfO\quad\text{uniformly on }[0,T].
\]
Let $\bfv_j$ satisfy
\[
 \bfv_j\in \rL^\infty(0,T;\cH(\Omega))\cap \rL^2(0,T;\cV(\Omega)),
 \qquad \bD(\bfv_j)=0\quad\text{a.e.\ on }S_j(t),
\]
with a uniform bound in these spaces. Assume moreover that for every $\sigma>0$ and every interval~${I\Subset(0,T)}$ there are $j_\sigma$ and $C_{\sigma,I}$ such that, for all $j\ge j_\sigma$,
\begin{equation}\label{eq:FSI temporal functional bound}
 \left|\int_I\int_\Omega \bfv_j\cdot\partial_t\bfzeta\srd\bfx\srd t\right|
 \le C_{\sigma,I}\Bigl(
 \|\bfzeta\|_{\rL^\infty(I;\rH^3(\Omega))}
 +\|\bfzeta\|_{\rL^2(I;\rH^3(\Omega))}
 \Bigr)
\end{equation}
for every $\bfzeta\in\rC_c^\infty(I\times\Omega;\R^3)$ satisfying
$\mdiv\bfzeta=0$ and $\bD(\bfzeta(t))=0$ on $S(t)_\sigma$
for every $t\in I$. Then, after extraction,
\[
 \bfv_j\to \bfv\qquad\text{strongly in }\rL^2(0,T;\rL^2(\Omega;\R^3)).
\]
\end{lem}

By the energy inequality \eqref{eq:energy est after n lim}, we imply that as $\delta \to 0$,
\begin{align}
 \bfu_\delta&\stackrel{*}{\rightharpoonup}\bfu
\text{ in }\rL^\infty(0,T_*;\cH(\Omega)), \quad \tand\quad  
 \bfu_\delta\rightharpoonup\bfu
 \text{ in }\rL^2(0,T_*;\cV(\Omega)).\label{eq:weak u delta}
\end{align}
The exact rigidity and \eqref{eq:strong phi delta} imply
$\varphi\bD(\bfu)=0$, so $\bfu\in \rL^2(0,T_*;\cK(\varphi))$.
For a smooth test rigid on~$\cS(t)_\sigma$, \eqref{eq:solid inclusions delta} makes it admissible in \eqref{eq:fluid delta} for all small $\delta$.  The convection and viscous terms have their standard energy bounds, and by \eqref{eq:global H1 director delta},
\[
\left|\int_\Omega(\nabla\overline\bfd_\delta)^\top\nabla\overline\bfd_\delta:\nabla\bfzeta\srd \bfx\right|
 \le C\|\bfzeta\|_{\rH^3(\Omega)}.
\]
Finally, we have
\begin{equation}\label{eq:artificial term compactness estimate}
 \left|\delta\int_0^{T_*}a_u(\bfu_\delta,\bfzeta)\srd t\right|
 \le \sqrt\delta\,\|\sqrt\delta\,\bfu_\delta\|_{\rL^2(0,T_*;X_u)}
 \|\bfzeta\|_{\rL^2(0,T_*;X_u)}
 \le C\sqrt\delta\,\|\bfzeta\|_{\rL^2(0,T_*;X_u)}.
\end{equation}
Thus \eqref{eq:FSI temporal functional bound} holds uniformly in $\delta$, and by \autoref{lem:FSI compactness criterion} we have
\begin{equation}\label{eq:strong u delta}
 \bfu_\delta\to\bfu
 \quad\text{strongly in }\rL^2(0,T_*;\rL^2(\Omega;\R^3)) \tas \delta \to 0.
\end{equation}
In particular, $\bfu_\delta\otimes\bfu_\delta\to\bfu\otimes\bfu
 \quad\text{strongly in }\rL^1((0,T_*)\times\Omega) \tas \delta \to 0.$

We can now pass to the limit for the equations.  From
\eqref{eq:strong phi delta} and \eqref{eq:strong u delta}, we obtain the transport identity
\begin{equation*}
 \int_0^{T_*}\int_\Omega\varphi
 [\partial_t\eta+\bfu\cdot\nabla\eta] \srd \bfx \srd t
 =-\int_\Omega\varphi_0\eta(0) \srd \bfx.
\end{equation*}
For a director test function compactly supported in $\cQ_\cF\cup(\{0\}\times\cF_0)$, the support lies in $\cF_\delta(t)$ for all small $\delta$.  Combining \eqref{eq:strong director delta local}, \eqref{eq:strong f director delta}, and \eqref{eq:strong u delta} we have
\begin{equation*}
 \int_0^{T_*}\int_{\cF(t)}
 \Bigl[-\bfd\cdot\partial_t\bfpsi
 +(\bfu\cdot\nabla)\bfd\cdot\bfpsi
 +\nabla\bfd:\nabla\bfpsi+\bff_\eps(\bfd)\cdot\bfpsi\Bigr] \srd \bfx \srd t
 =\int_{\cF_0}\bfd_0\cdot\bfpsi(0) \srd \bfx. 
\end{equation*}
By \autoref{lem:density director tests}, this identity can be extended to the natural moving energy test class for fixed $\eps$.

For the momentum equation, let $\bfzeta$ be a smooth test function that is rigid in $\cS(t)_\sigma$.  The artificial-viscosity term vanishes by \eqref{eq:artificial term compactness estimate}.  We can pass to the limit for all velocity terms by \eqref{eq:weak u delta} and \eqref{eq:strong u delta}.  The director term needs only \eqref{eq:strong stress delta local}. Indeed, since the Ericksen tensor is symmetric,
\[
 (\nabla\bfd_\delta)^\top\nabla\bfd_\delta:\nabla\bfzeta
 =(\nabla\bfd_\delta)^\top\nabla\bfd_\delta:\bD(\bfzeta)=0
 \quad\text{on }\cS(t)_\sigma,
\]
and the test function is compactly supported outside the collar.
The initial and terminal time layers are controlled directly by the uniform director energy.  Hence
\begin{equation*}
\begin{aligned}
 &\int_0^{T_*}\int_\Omega
 \Bigl[-\bfu\cdot\partial_t\bfzeta-(\bfu\otimes\bfu):\nabla\bfzeta
 +2\bD(\bfu):\bD(\bfzeta)\Bigr]\srd\bfx\srd t\\
 &\qquad=\int_\Omega\bfu_0\cdot\bfzeta(0)\srd\bfx
 +\int_0^{T_*}\int_{\cF(t)}(\nabla\bfd)^\top\nabla\bfd:\nabla\bfzeta\srd\bfx\srd t.
\end{aligned}
\end{equation*}
The global $\rH^1$ material extension makes the limiting momentum functional continuous under \eqref{eq:collar density convergence}.  Hence, by \autoref{prop:rigid test equivalence} this identity can be extended to the standard test class $\cT_{\bfu}$.  In the subsequent singular limit we use the collar subclass only as an approximation device.

To pass to the limit for the energy inequality, we extract one common subsequence.  By strong
convergences \eqref{eq:strong director delta L2 global} and \eqref{eq:strong u delta}, we can
 take, up to a subsequence,
\[
 \overline\bfd_\delta(t)\to\overline\bfd(t)
 \quad\text{and}\quad
 \bfu_\delta(t)\to\bfu(t)
 \qquad\text{strongly in }\rL^2(\Omega) \tas \delta\to 0,
\]
for a.e.\ $t\in(0,T_*)$.  Fix such a time that is also a Lebesgue time for all energy terms.  The
uniform~$\rH^1$-bound \eqref{eq:global H1 director delta} implies that every weak $\rH^1$-cluster point
of $\overline\bfd_\delta(t)$ equals $\overline\bfd(t)$, so
\[
 \|\nabla\overline\bfd(t)\|_{\rL^2(\Omega)}^2
 \le\liminf_{\delta\to0}
 \|\nabla\overline\bfd_\delta(t)\|_{\rL^2(\Omega)}^2.
\]
For fixed $\eps$, by the maximum principle we have $|\overline\bfd_\delta|\le1$, and $F_\eps$ is Lipschitz
on the closed unit ball. Hence, the strong convergence $\overline\bfd_\delta(t)\to\overline\bfd(t)$ in $\rL^2(\Omega)$ at the fixed time $t$ implies the strong $\rL^1$-convergence of the
Ginzburg--Landau potential.  The kinetic energy converges along this common subsequence.  For the fluid
Dirichlet energy we subtract the rigid-motion invariant solid contribution
\[
 \frac12\int_{\cS_\delta(t)}|\nabla\overline\bfd_\delta|^2\srd\bfx
 =\frac12\int_{\cS_0}|\nabla\bfd_0|^2\srd\bfx
 =\frac12\int_{\cS(t)}|\nabla\overline\bfd|^2\srd\bfx.
\]
By the weak lower semicontinuity of the viscous and molecular-field dissipations, and by
\eqref{eq:weak M delta}, for a.e.\ $t\in(0,T_*)$ we have
\begin{equation*}
 E_\eps(t)+2\int_0^t\int_\Omega|\bD(\bfu)|^2\srd\bfx\srd s
 +\int_0^t\int_{\cF(s)}|\Delta\bfd-\bff_\eps(\bfd)|^2\srd\bfx\srd s
 \le E_\eps(0),
\end{equation*}
where
\[
 E_\eps(t)=\frac12\|\bfu(t)\|_{\rL^2(\Omega)}^2
 +\int_{\cF(t)}\left(\frac12|\nabla\bfd(t)|^2+F_\eps(\bfd(t))\right)\srd\bfx.
\]
We can pass to the limit in the separation estimate \eqref{eq:uniform noncontact delta}. This completes the limit $\delta\to0$ and the proof of \autoref{thm:weak sol GL}.

\subsection{Fixed-domain Ginzburg--Landau theory on arbitrary open sets}\label{ssec:fixed eps continuation}
Since the post-impact fluid region need not be Lipschitz, we establish fixed-$\eps$ existence and compactness for the Ginzburg--Landau system on an arbitrary bounded open set.  This construction is performed after the limiting contact configuration has been identified.

For a bounded open set $G\subset\R^3$, let $\cD_G$, $\cH_G$, and $\cV_G$ be defined in \eqref{eq:arbitrary domain solenoidal spaces}.  For $\bfv\in\cD_G$, integrating by parts we have
\[
\|\nabla\bfv\|_{\rL^2(G)}^2=2\|\bD(\bfv)\|_{\rL^2(G)}^2,
\]
and by the density argument, we can extend this identity to every $\bfv\in\cV_G$.  Thus the viscous dissipation controls the full $\rH_0^1$-norm on $\cV_G$ without any assumption on $\partial G$.

\begin{lem}\label{lem:fixed domain velocity compactness}
Let $G\subset\R^3$ be bounded and open.  Let $(\bfv_m)_{m\ge1}\subset\cD_G$ be $\rL^2(G)$-orthonormal, with nested spans
$X_N:=\operatorname{span}\{\bfv_1,\ldots,\bfv_N\}$ whose union is dense in both $\cH_G$ and $\cV_G$.  Suppose that
$(\bfw_j)$ is bounded in $\rL^\infty(0,T;\cH_G)\cap \rL^2(0,T;\cV_G)$
and that, for every fixed $m$, the scalar functions $ a_{j,m}(t):=(\bfw_j(t),\bfv_m)_{\rL^2(G)}$
have uniformly bounded variation on $(0,T)$ for all sufficiently large $j$.  Then $(\bfw_j)$ is relatively compact in $\rL^2((0,T)\times G)$.
\end{lem}

\begin{proof}
Let $P_N$ be the $\rL^2(G)$-orthogonal projection onto $X_N$.  Since $\rH_0^1(G)\Subset \rL^2(G)$ for every bounded open set by the Rellich--Kondrachov theorem with zero extension to a containing ball, the embedding $\cV_G\Subset\cH_G$ is compact.  Hence
\begin{equation}\label{eq:fixed domain high mode compactness}
 \alpha_N:=\sup_{\substack{\bfz\in\cV_G\\ \|\bfz\|_{\cV_G}\le1}}
 \|(I-P_N)\bfz\|_{\rL^2(G)}\rightarrow0.
\end{equation}
Indeed, the image of the unit ball of $\cV_G$ is relatively compact in $\cH_G$, and $P_N\bfz\to \bfz$ in $\cH_G$ for every~$\bfz\in\cH_G$.
The convergence is therefore uniform on that compact set.

For fixed $N$, each coefficient of $P_N\bfw_j$ is bounded in $\rL^\infty(0,T)$ and in $\mathrm{BV}(0,T)$. Applying Helly's compactness theorem, these coefficients are relatively compact in $\rL^1(0,T)$. The uniform $\rL^\infty(0,T)$-bound implies that they are relatively compact in $\rL^2(0,T)$.  Thus, $(P_N\bfw_j)_j$ is relatively compact in $\rL^2(0,T;\cH_G)$.  On the other hand, by \eqref{eq:fixed domain high mode compactness},
\[
 \|(I-P_N)\bfw_j\|_{\rL^2(0,T;\cH_G)}
 \le \alpha_N\|\bfw_j\|_{\rL^2(0,T;\cV_G)},
\]
uniformly in $j$.  Letting $N\to\infty$, the claim is derived.
\end{proof}

\begin{lem}\label{lem:arbitrary open maximum principle}
Let $G\subset\R^3$ be bounded and open, let
$\bfu\in \rL^2(0,T;\rH_0^1(G;\R^3))$ be solenoidal, and let~$\bb\in \rH^1(G;K_+)$. Suppose
\[
 \bfd\in\rL^\infty(0,T;\rH^1(G;\R^3)), \qquad \bfd-\bb\in \rL^2(0,T;\rH_0^1(G;\R^3)),\qquad
 \partial_t\bfd\in \rL^2(0,T;\rH^{-1}(G;\R^3)),
\]
$\bfd(0)\in K_+$ a.e., and
\begin{equation}\label{eq:arbitrary open parabolic}
 \partial_t\bfd+(\bfu\cdot\nabla)\bfd
 =\Delta\bfd+\frac1{\eps^2}(1-|\bfd|^2)\bfd
 \quad\text{in }\rL^2(0,T;\rH^{-1}(G)).
\end{equation}
Then $|\bfd|\le1$ and $d_3\ge0$ a.e.\ on $(0,T)\times G$.
\end{lem}

\begin{proof}
Let us write $\bfd=\bb+\bfw$ with $\bfw\in\rL^2(0,T;\rH_0^1(G;\R^3))$. Since $\partial_t\bfw=\partial_t\bfd\in\rL^2(0,T;\rH^{-1}(G;\R^3))$, we have $\bfd\in\rC([0,T];\rL^2(G;\R^3))$. 

We first show that if $\Psi\colon\R^3\to\R^k$ is Lipschitz with $\Psi=0$ on $K_+$, then $\Psi(\bfd(t))\in\rH_0^1(G;\R^k)$ for a.e. $t$. Indeed, fix $t$ and choose $\bfw_j\in\rC_c^\infty(G;\R^3)$ with $\bfw_j\to\bfw(t)$ in $\rH^1(G)$. Since $\bb$ takes values in $K_+$, the function $\Psi(\bb+\bfw_j)$ vanishes outside $\supp\bfw_j$. Hence, it has compact support in $G$ and belongs to $\rH_0^1(G;\R^k)$. By the chain rule for Lipschitz maps, these functions $\Psi(\bb+\bfw_j)$ are bounded in $\rH^1(G)$, and they converge to $\Psi(\bfd(t))$ in $\rL^2(G)$. Since $\rH_0^1(G)$ is weakly closed in $\rH^1(G)$, we conclude that $\Psi(\bfd(t))\in\rH_0^1(G;\R^k)$.

\smallskip
Next we prove that $|\bfd|\le1$. Let $\Pi$ be the nearest-point projection onto the closed unit ball $B\subset\R^3$, and define
\[
 \mathbf F(\bfz):=\bfz-\Pi(\bfz)=\Bigl(1-\frac1{|\bfz|}\Bigr)_+\bfz,
 \qquad
 \Phi(\bfz):=\frac12|\mathbf F(\bfz)|^2=\frac12\bigl(|\bfz|-1\bigr)_+^2 .
\]
The function $\Phi$ is convex and of class $\rC^{1,1}$ with $\nabla\Phi=\mathbf F$, and $\mathbf F$ is Lipschitz and vanishes on $B\supset K_+$. Hence $\mathbf F(\bfd)\in\rL^2(0,T;\rH_0^1(G;\R^3))$ by the above arguments.

Now test equation \eqref{eq:arbitrary open parabolic} with $\mathbf F(\bfd)$. Due to the convexity of $\Phi$, for all $t,s \in (0,T)$ we have
\[
 \mathbf F(\bfd(s))\cdot\bigl(\bfd(t)-\bfd(s)\bigr)
 \le\Phi(\bfd(t))-\Phi(\bfd(s))
 \le\mathbf F(\bfd(t))\cdot\bigl(\bfd(t)-\bfd(s)\bigr).
\]
Integrating these inequalities over $G$ and using the Steklov averages in time, we imply that the integral of $\Phi(\bfd(t))$ over $G$ is absolutely continuous as a function of $t$ and
\begin{equation}\label{eq:rel Phi and F}
    \frac{\rd}{\rd t}\int_G\Phi(\bfd)\srd\bfx=\langle\partial_t\bfd,\mathbf F(\bfd)\rangle .
\end{equation}
By definition of $\mathbf F$, we have on $\{|\bfd|>1\}$, $\partial_i\mathbf F(\bfd)=\Bigl(1-\frac1{|\bfd|}\Bigr)\partial_i\bfd+\frac{\bfd\cdot\partial_i\bfd}{|\bfd|^3}\,\bfd,$
and $\nabla\mathbf F(\bfd)=0$ a.e.\ on $\{|\bfd|\le1\}$. Consequently,
\begin{equation}\label{e:aux1}
  \nabla\bfd:\nabla\mathbf F(\bfd)
 =\Bigl(1-\frac1{|\bfd|}\Bigr)_+|\nabla\bfd|^2
 +\dsone_{\{|\bfd|>1\}}\frac{1}{|\bfd|^3}\sum_{i=1}^3(\bfd\cdot\partial_i\bfd)^2\ \ge0 .   
\end{equation}
Similarly, we have $\mathbf F(\bfd)\cdot\partial_i\bfd=\mathbf F(\bfd)\cdot\partial_i\mathbf F(\bfd)=\partial_i\Phi(\bfd)$, hence, $(\bfu\cdot\nabla)\bfd\cdot\mathbf F(\bfd)=\bfu\cdot\nabla\Phi(\bfd).$
Choose $\bfphi_j\in\rC_c^\infty(G;\R^3)$ with $\bfphi_j\to\mathbf F(\bfd(t))$ in $\rH^1(G)$. Since $\mdiv\bfu=0$, we have $\int_G\bfu\cdot\nabla\tfrac12|\bfphi_j|^2\srd\bfx=0$, and since $\bfphi_j\cdot\nabla\bfphi_j\to\mathbf F(\bfd)\cdot\nabla\mathbf F(\bfd)$ in $\rL^{3/2}(G)$ and $\bfu(t)\in\rL^6(G)$, we have
\begin{equation}\label{e:aux2}
\int_G(\bfu\cdot\nabla)\bfd\cdot\mathbf F(\bfd)\srd \bfx = 0 \quad \text{ for a.e. } t. 
\end{equation}
Note that $\bfd\cdot\mathbf F(\bfd)=|\bfd|(|\bfd|-1)_+$, so  $(1-|\bfd|^2)\,\bfd\cdot\mathbf F(\bfd)\le0$. Combining this inequality with \eqref{e:aux1}, \eqref{e:aux2}, we obtain
\[
\frac{\rd}{\rd t}\int_G\Phi(\bfd)\srd\bfx = - \int_G(\bfu\cdot\nabla)\bfd\cdot\mathbf F(\bfd)\srd \bfx + \int_G  \nabla\bfd:\nabla\mathbf F(\bfd)\srd \bfx + \dfrac{1}{\eps^2}\int_G (1-|\bfd|^2)\,\bfd\cdot\mathbf F(\bfd)\srd\bfx \le0.
\]
As $\Phi\ge0$ and $\Phi(\bfd(0))=0$, it implies $\Phi(\bfd(t))=0$, that is, $|\bfd|\le1$ a.e.\ in $(0,T)\times G$.

\smallskip
Lastly, we show that $d_3\ge0$. Let us denote $q:=d_3$. By our assumption \eqref{eq:arbitrary open parabolic}, it satisfies
\begin{equation}\label{e:aux3}
 \partial_tq+\bfu\cdot\nabla q-\Delta q= \frac{1}{\eps^2}\bigl(1-|\bfd|^2\bigr)q
 \quad\text{in }\rL^2(0,T;\rH^{-1}(G)),    
\end{equation}
and $0\le \frac{1}{\eps^2}(1-|\bfd|^2)\le\eps^{-2}$ by the above verification that $|\bfd| \le 1$. The negative part $q_-=\max\{-q,0\}$ equals $\Psi(\bfd)$ for the Lipschitz function $\Psi(\bfz):=(z_3)_-$, which vanishes on $K_+$. Therefore, we have $q_-\in\rL^2(0,T;\rH_0^1(G))$. Moreover, we have the identities
\[
 \langle\partial_tq,-q_-\rangle=\frac12\frac{\rd}{\rd t}\|q_-\|_{\rL^2(G)}^2,
 \qquad \nabla q\cdot\nabla(-q_-)=|\nabla q_-|^2,
 \qquad q\,(-q_-)=q_-^2,
\]
where the first identity can be derived as in \eqref{eq:rel Phi and F}. Testing \eqref{e:aux3} with $-q_-$, using these identities, and noting that the transport term vanishes as in the above computation of $\langle\partial_t\bfd,\mathbf F(\bfd)\rangle$, we get
\[
 \frac12\frac{\rd}{\rd t}\|q_-\|_{\rL^2(G)}^2+\|\nabla q_-\|_{\rL^2(G)}^2
 =\int_G\frac{1}{\eps^2}(1-|\bfd|^2)\,q_-^2\srd\bfx
 \le\frac{1}{\eps^2}\|q_-\|_{\rL^2(G)}^2 .
\]
Since $q_-(0)=0$, by Gr\"onwall's inequality we have $q_-=0$ meaning that $d_3\ge0$ a.e.\ in $(0,T)\times G$.
\end{proof}

We first solve the Ginzburg--Landau relaxation of \eqref{eq:simplified EL-fixed dom} on an arbitrary bounded open set $G$,
\begin{equation}\label{eq:GL fixed dom}
    \left\{
    \begin{aligned}
        \partial_t \bfu + (\bfu \cdot \nabla)\bfu - \Delta\bfu + \nabla\pi &= -\mdiv\bigl((\nabla\bfd)^\top \nabla\bfd\bigr), \quad \mdiv\bfu = 0, &&\text{in } (0,T) \times G,\\
        \partial_t \bfd + (\bfu \cdot \nabla)\bfd - \Delta\bfd &= -\bff_\eps(\bfd), &&\text{in } (0,T) \times G,
    \end{aligned}
    \right.
\end{equation}
with $\bfu=0$ and $\bfd=\bfd_*$ on $\partial G$ in the variational sense, and with initial data $(\bfu_*,\bfd_*)$.

\begin{prop}
\label{prop:postcontact GL arbitrary domain}
Let $G\subset\R^3$ be a bounded open set, let $\bfu_*\in\cH_G$, and let $\bfd_*\in \rH^1(G;K_+)$.  Let $\eps>0$ be fixed.  Then there is a global trajectory
$(\bfu,\bfd,\bfM)$ on $[0,\infty)$ such that, for every $T>0$,
\begin{align*}
 &\bfu\in \rL^\infty(0,T;\cH_G)\cap \rL^2(0,T;\cV_G),\\
 &\bfd\in \rL^\infty(0,T;\rH^1(G;K_+)),\qquad \bfd-\bfd_*\in \rL^\infty(0,T;\rH_0^1(G;\R^3)),\\
 &\bfM:=\Delta \bfd-\bff_\eps(\bfd)\in \rL^2((0,T)\times G;\R^3),
\end{align*}
where the molecular-field identity is understood variationally, that is, for a.e. $t$
\begin{equation*}
 \int_G \bfM(t)\cdot\bfxi\srd\bfx
 =-\int_G\nabla \bfd(t):\nabla\bfxi\srd\bfx
   -\int_G\bff_\eps(\bfd(t))\cdot\bfxi\srd\bfx
 \qquad\forall\bfxi\in \rH_0^1(G;\R^3).
\end{equation*} 
The pair $(\bfu,\bfd)$ solves \eqref{eq:GL fixed dom} in the following weak form:
for every
$\bfzeta\in \rC_c^1([0,T);\cD_G)$ and every~$\bfpsi\in \rC_c^1([0,T);\rC_c^\infty(G;\R^3))$,
\begin{align*}
 &\int_0^T\int_G\Bigl[-\bfu\cdot\partial_t\bfzeta-(\bfu\otimes \bfu):\nabla\bfzeta   +2\bD(\bfu):\bD(\bfzeta)\Bigr]\srd\bfx\srd t\\
&\qquad=\int_G\bfu_*\cdot\bfzeta(0)\srd\bfx
       +\int_0^T\int_G(\nabla \bfd)^\top\nabla \bfd:\nabla\bfzeta\srd\bfx\srd t,\\
 &\int_0^T\int_G\Bigl[-\bfd\cdot\partial_t\bfpsi+(\bfu\cdot\nabla)\bfd\cdot\bfpsi
       +\nabla \bfd:\nabla\bfpsi+\bff_\eps(\bfd)\cdot\bfpsi\Bigr]\srd\bfx\srd t
   =\int_G\bfd_*\cdot\bfpsi(0)\srd\bfx.
\end{align*}
For a.e.\ $t\in(0,T)$ we have
\begin{equation}\label{eq:postcontact energy GL}
\begin{aligned}
 &\frac12\|\bfu(t)\|_{\rL^2(G)}^2
 +\frac12\|\nabla \bfd(t)\|_{\rL^2(G)}^2+\int_G F_\eps(\bfd(t))\srd\bfx +2\int_0^t\int_G|\bD(\bfu)|^2\srd\bfx\srd s\\
 &\quad
 +\int_0^t\int_G|\bfM|^2\srd\bfx\srd s
 \le
 \frac12\|\bfu_*\|_{\rL^2(G)}^2
 +\frac12\|\nabla \bfd_*\|_{\rL^2(G)}^2+\int_G F_\eps(\bfd_*)\srd\bfx.
\end{aligned}
\end{equation}
\end{prop}

\begin{proof}
Note that $G$ may fail to be Lipschitz, we will utilize compactness properties that do not require boundary regularity.  Let $\bb_*$ be
the componentwise harmonic representative of the affine class $\bfd_*+\rH_0^1(G)$, i.e.,
\[
 \bb_*-\bfd_*\in \rH_0^1(G),\qquad
 \int_G\nabla \bb_*:\nabla\bfxi\srd\bfx=0
 \quad\forall\bfxi\in \rH_0^1(G;\R^3).
\]
We claim that $\bb_*$ takes values in $K_+$. 
By definition \eqref{e:Kplus}, we can write
\[
K_+ = \{\bfz:z_3\ge0\} \bigcap_{\bfe\in\bS^2}  \{\bfz:\bfe\cdot\bfz\le1\}.
\]
Let $\ell(\bfz)=\bfe\cdot\bfz-1$ or $\ell(\bfz)=-z_3$, so that $\ell(\bfd_*)\le0$. 
We claim that the scalar function $h:=\ell(\bb_*)$ is weakly harmonic, and $h_+= \max\{h,0\}\in\rH_0^1(G)$. Indeed, choose $\bfq_m\in\rC_c^\infty(G;\R^3)$ with $\bfq_m\to\bb_*-\bfd_*$ in $\rH_0^1(G)$.
Then $\ell(\bfd_*+\bfq_m)_+$ has compact support in $G$ because $\ell(\bfd_*)\le0$, is bounded in $\rH^1(G)$, and converges to $h_+$ in $\rL^2(G)$. By the weak closedness of $\rH_0^1(G)$, $h_+\in\rH_0^1(G)$. 
The harmonicity of $h$ also implies that $\|\nabla h_+\|_{\rL^2(G)}^2=0$, hence $h_+=0$ by Poincar\'e's inequality. 
Applying this argument to $-z_3$ and to a countable dense set of $\bfe\in\bS^2$, we imply that $\bb_*\in K_+$ almost everywhere.

Choose a sequence $(\bfv_j)_{j\ge1}\subset\cD_G$ whose span is dense both in $\cH_G$ with respect to
$\rL^2(G)$ and in $\overline{\cD_G}^{\rH_0^1(G)}$ with respect to $\rH_0^1(G)$, and orthonormalize it in
$\rL^2(G)$.  For fixed $N$, we consider velocities of the form $\bfu^N(t)=\sum_{j=1}^Na_j(t)\bfv_j.$
If
$a\in \rC([0,T];\R^N)$ is prescribed, the vector field $\bfu^N$ is smooth in space, divergence free, and has
zero trace.  We solve the equation
\[
 \partial_t \bfd+(\bfu^N\cdot\nabla)\bfd=\Delta \bfd-\bff_\eps(\bfd),
 \qquad \bfd(0)=\bfd_*,\qquad \bfd-\bb_*\in \rL^2(0,T;\rH_0^1(G)),
\]
by the eigenfunction Galerkin scheme for the Dirichlet Laplacian associated with $\rH_0^1(G)$.  The energy estimate and Banach--Alaoglu give weak and weak-* compactness at the energy level.  After the auxiliary Galerkin index is removed, \autoref{lem:arbitrary open maximum principle} applies with boundary representative $\bb_*$ and shows that the resulting parabolic solution satisfies $\bfd\in K_+$ almost everywhere.  To show the dependence on the given velocity, let
$\bfd_1,\bfd_2$ correspond to coefficient paths $a_1,a_2$ in a bounded subset of
$\rC([0,T];\R^N)$ and set~$\bfr=\bfd_1-\bfd_2$.  Testing the difference equation by $\bfr$, using
$\mdiv\bfu[a_1]=0$, the Lipschitz continuity of $\bff_\eps$ on $K_+$, and the finite-dimensional bound
$\|\bfu[a_1]-\bfu[a_2]\|_{\rL^\infty(G)}\le C_N|a_1-a_2|$, we have
\[
 \frac12\frac{\rd}{\rd t}\|\bfr\|_{\rL^2(G)}^2+\|\nabla \bfr\|_{\rL^2(G)}^2
 \le C_{N,\eps,R}\|\bfr\|_{\rL^2(G)}^2
 +C_{N,R}|a_1-a_2|^2\|\nabla\bfd_2\|_{\rL^2(G)}^2.
\]
Applying Gr\"onwall's lemma and using the uniform director-energy bound we imply that $\bfd[a]$ depends continuously on $a$, as a map from bounded subsets of $\rC([0,T];\R^N)$ into
$\rC([0,T];\rL^2(G))\cap \rL^2(0,T;\rH^1(G))$.

The finite-dimensional coupling can be closed without any regularity of $\partial G$.
Given a coefficient path $a$, let $\bfd[a]$ be the director constructed as above and let
$\mathcal T_N(a)=b$ where $b$ solves the following equation
\[
 b_i'(t)+2\sum_{j=1}^N (\bD(\bfv_j),\bD(\bfv_i))_{\rL^2(G)}b_j(t)
 =\int_G(\bfu[a]\otimes \bfu[a]):\nabla\bfv_i\srd\bfx+\int_G(\nabla \bfd[a])^\top\nabla \bfd[a]:\nabla\bfv_i\srd\bfx,
\]
with initial datum to be the $\rL^2(G)$-orthogonal projection of $\bfu_*$.  Here
$\bfu[a]=\sum_{j=1}^Na_j\bfv_j$.  Let us fix $R>0$ and take $a$ in the closed ball of radius $R$ in
$\rC([0,T_1];\R^N)$ for some $T_1$ to be determined.  Then, on this ball the parabolic energy estimate for the given smooth drift $\bfu[a]$ implies a bound for $\bfd[a]$ in $\rL^\infty(0,T_1;\rH^1(G))$ depending only on $R$, $N$, $\eps$, and the initial
energy.  Then, using the integral form of the above ODE we have
\[
 \|\mathcal T_N(a)(t)-\mathcal T_N(a)(s)\|_{\R^N}
 \le C|t-s|^{1/2},
\]
where $C>0$ depending on $N, R,\eps$ and initial energy. Choose $R$ larger than the norm of the projected initial coefficients and then choose
$T_1$ sufficiently small, $\mathcal T_N$ maps the above ball into itself.  The preceding difference estimate
for the director implies that $\mathcal T_N$ is continuous, and the above modulus implies the relative compactness
by the Arzel\`a--Ascoli theorem. Thus, we can apply Schauder's theorem to conclude that there exists a fixed point on $[0,T_1]$.

After removing the auxiliary director Galerkin index, we derive the a priori estimate at the fixed point.  For the 
above velocity $\bfu^N$, we first remove the auxiliary Dirichlet Galerkin index in the director
equation.  The resulting parabolic solution has the full molecular field
\[
   \bfM^N=\Delta \bfd^N-\bff_\eps(\bfd^N)\in \rL^2((0,T_1)\times G)
\]
and satisfies the director energy identity.    Since
$\Delta\bfd^N = \bfM^N+\bff_\eps(\bfd^N) \in \rL^2(G)$, and since $\bfu^N$ has compact support in $G$, the interior elliptic regularity implies $\bfd^N\in \rH^2$ on a neighborhood of $\supp\bfu^N$.  Then, since $\mdiv \bfu^N=0$, by integrating by parts we have
\[
 \int_G (\bfu^N\cdot\nabla)\bfd^N\cdot \bfM^N\srd\bfx
 =-\int_G (\nabla \bfd^N)^\top\nabla \bfd^N:\nabla \bfu^N\srd\bfx.
\]
Indeed, the contribution of $-\bff_\eps(\bfd^N)=-\nabla_{\bfd}F_\eps(\bfd^N)$ is the integral of
$-\bfu^N\cdot\nabla F_\eps(\bfd^N)$ and vanishes, while the Laplacian term yields the Ericksen stress.
Testing the finite-dimensional momentum equation by $\bfu^N$ and adding the director energy identity we get
\eqref{eq:postcontact energy GL}.  It implies that the semi-Galerkin energy is non-increasing.  In particular, the
coefficient vector and the director energy remain bounded on every finite interval, so the local Schauder
construction can be continued with a lifespan depending only on the current bounded quantities. There is
no finite-time accumulation of restart times.  Hence $(\bfu^N,\bfd^N)$ exists on every finite interval.

The energy inequality \eqref{eq:postcontact energy GL} for $(\bfu^N,\bfd^N,\bfM^N)$ implies the uniform bounds
\begin{equation}\label{eq:postcontact Galerkin bounds}
 \bfu^N\text{ bounded in }\rL^\infty_t\rL^2_{\bfx}\cap \rL^2_t\rH^1_{0,\bfx},
 \quad \bfd^N\text{ bounded in }\rL^\infty_t\rH^1_{\bfx},
 \quad \bfM^N\text{ bounded in }\rL^2_{t,\bfx}.
\end{equation}
Since $\bfd^N\in K_+$ and $\eps$ is fixed, $\bff_\eps(\bfd^N)$ is bounded in $\rL^\infty_{t,\bfx}$.  Set
$\bfw^N:=\bfd^N-\bb_*\in \rH_0^1(G;\R^3)$. Since $\bb_*$ is weakly harmonic, by the variational
molecular-field equation we have
\[
  -\Delta \bfw^N=-\bfM^N-\bff_\eps(\bfd^N)
\]
in $\rH^{-1}(G)$. The right-hand side belongs to $\rL^2(G)$, hence $\bfw^N$ is bounded in
$\rL^2(0,T;D(A_D))$, where $A_D$ is the Dirichlet Laplacian on the arbitrary open set $G$.  Because $G$ is bounded, by the Rellich--Kondrachov theorem with  zero extension to a containing ball we have
$\rH_0^1(G)\Subset \rL^2(G)$, hence $A_D$ has compact resolvent.  If $A_D\bfe_j=\lambda_j \bfe_j$ is its spectral resolution, then, with
$c_j=(\bfw,\bfe_j)_{\rL^2(G)}$,
\[
\|\nabla\bfw\|_{\rL^2(G)}^2=\sum_j\lambda_j|c_j|^2,
 \qquad
 \|A_D\bfw\|_{\rL^2(G)}^2=\sum_j\lambda_j^2|c_j|^2.
\]
The first identity is the quadratic-form identity for the variational
Dirichlet Laplacian. On the bounded set $G$, its left-hand side is an equivalent norm
on $\rH_0^1(G)$ by Poincar\'e's inequality. Since $\lambda_j\to\infty$ as $j\to \infty$, the
tails of the first series are uniformly small on bounded subsets of
$D(A_D)$.  Thus, without any boundary regularity assumption,
\[
  D(A_D)\Subset \rH_0^1(G)\Subset \rL^2(G).
\]
Moreover
\[
  \partial_t \bfw^N=\bfM^N-\bfu^N\cdot\nabla \bfd^N
\]
is bounded in $\rL^2(0,T;\rH^{-1}(G))$. Therefore, by the Aubin--Lions lemma, after extracting a subsequence we have
\begin{equation}\label{eq:postcontact strong H1}
   \bfd^N\to \bfd\qquad\text{strongly in }\rL^2(0,T;\rH^1(G)).
\end{equation}
For the velocity we use \autoref{lem:fixed domain velocity compactness}.  The basis $(\bfv_m)$ was chosen above so that its spans are dense in $\cH_G$ and $\cV_G$ as introduced in \autoref{ssec:main result}.  For every fixed $m$ and every $N\ge m$, by the Galerkin equation we have
\[
 \frac{\rd}{\rd t}(\bfu^N,\bfv_m)_{\rL^2(G)}
 =\int_G(\bfu^N\otimes\bfu^N):\nabla\bfv_m\srd\bfx
 -2\int_G\bD(\bfu^N):\bD(\bfv_m)\srd\bfx
 +\int_G(\nabla\bfd^N)^\top\nabla\bfd^N:\nabla\bfv_m\srd\bfx.
\]
The uniform bounds \eqref{eq:postcontact Galerkin bounds} therefore yield
\[
  \int_0^T\left|\frac{\rd}{\rd t}(\bfu^N,\bfv_m)_{\rL^2(G)}\right|\srd t
 \le C_m\left(
 T\|\bfu^N\|_{\rL^\infty_t\rL^2_{\bfx}}^2
 +T^{1/2}\|\bD(\bfu^N)\|_{\rL^2_{t,\bfx}}
 +T\|\nabla\bfd^N\|_{\rL^\infty_t\rL^2_{\bfx}}^2\right)\le C_{m,T}.
\]
Thus the functions $t\mapsto(\bfu^N(t),\bfv_m)_{\rL^2(G)}$, $N\ge m$, have uniformly bounded variation on $(0,T)$, and by \autoref{lem:fixed domain velocity compactness} we have, after extracting a subsequence, $\bfu^N\to\bfu$ strongly in $\rL^2((0,T)\times G)$. By \eqref{eq:postcontact Galerkin bounds}, we may also assume that $\bfM^N\rightharpoonup\bfM$ weakly in $\rL^2((0,T)\times G)$. Together with \eqref{eq:postcontact strong H1}, these convergences allow us to pass to the limit in the Galerkin momentum equation, in the director equation for $\bfd^N$, and in the relation $\bfM^N=\Delta\bfd^N-\bff_\eps(\bfd^N)$. This implies the two weak equations and the molecular-field identity of the statement. Note that $(\nabla\bfd^N)^\top\nabla\bfd^N\to(\nabla\bfd)^\top\nabla\bfd$ in $\rL^1((0,T)\times G)$, which suffices for test fields in $\cD_G$. 

Since $\bfd^N\in K_+$ by \autoref{lem:arbitrary open maximum principle}, the strong convergence implies $\bfd\in K_+$ a.e., and since $\bfd^N-\bb_*$ and $\bb_*-\bfd_*$ belong to $\rH_0^1(G;\R^3)$, the weak closedness of this space implies $\bfd-\bfd_*\in\rH_0^1(G;\R^3)$. 

Finally, we recall that $(\bfu^N,\bfd^N,\bfM^N)$ satisfies the energy inequality \eqref{eq:postcontact energy GL}. After a further extraction, $\bfu^N(t)\to\bfu(t)$ in $\rL^2(G)$ and $\bfd^N(t)\to\bfd(t)$ in $\rH^1(G)$ for almost everywhere $t$, so the first three terms on its left-hand side converge, and the two dissipation integrals are weakly lower semicontinuous. Hence \eqref{eq:postcontact energy GL} holds for the limit. Since the approximations exist on every bounded time interval with the same bounds, a diagonal argument over $T=1,2,\dots$ yields a single trajectory on $[0,\infty)$.
\end{proof}

\section{The singular limit $\eps\to0$ on collision-free intervals}\label{sec:eps to 0}

In this section, we carry out the limit passage $\eps\to0$ in the Ginzburg--Landau regularization, see \autoref{thm:weak sol GL} for the existence result. 
We refer to the resulting solutions in which the constraint $|\bfd|=1$ holds as \emph{constrained weak solutions} as opposed to the solutions with regularization parameter~$\eps > 0$.
The following proposition records the existence of such constrained weak solutions.

\begin{prop}\label{prop:constrained local existence}
Under the assumptions of \autoref{exist:main}, there exists a constrained weak solution~$(\bfu,\bfd,\varphi)$ on $[0,T_*)$, that is, a weak solution on $[0,T]$ in the sense of \autoref{def:weak sol} for every~$T<T_*$, where $T_*$ is given by \eqref{eq:T star main theorem}. It satisfies $\dist(\cS(t),\partial\Omega)\ge\delta_0/2$ for $0\le t<T_*$.
\end{prop}

The rest of this section is devoted to the proof of \autoref{prop:constrained local existence}.

In this section and in \autoref{ssec:test transport auxiliary}, we also use the collar-rigid subclass of the standard momentum test class \eqref{eq:standard momentum tests},
\begin{equation}\label{eq:collar momentum tests}
\cT_{\bfu}^{\rm col}:=\bigl\{\bfzeta\in\cT_{\bfu}:\exists r>0\text{ such that }
\bD(\bfzeta(t))=0\text{ on }\cS(t)_r\text{ for every }t\in[0,T]\bigr\},
\end{equation}
which serves as an approximation device. The final weak formulation uses \eqref{eq:standard momentum tests}.

We return to the datum $\bfd_0\in \rH^1(\Omega;\bS_+^2)$ from \autoref{exist:main}.  By \autoref{lem:hemisphere preparation}, choose a family
\begin{equation}\label{eq:prepared anchors eps}
 \bfd_0^\eps\in \rH^2(\Omega;\bS_+^2),\qquad
 \bfd_0^\eps\to\bfd_0\quad\text{strongly in }\rH^1(\Omega)
 \quad(\eps\to0).
\end{equation}
Set
\[
 E_0^\eps:=\frac12\|\bfu_0\|_{\rL^2(\Omega)}^2+
 \frac12\|\nabla\bfd_0^\eps\|_{\rL^2(\cF_0)}^2,
 \qquad E_0^\eps\to E_0.
\]
Let $T_*$ be the time in \eqref{eq:T star main theorem}.
In the case $E_0=0$, we get $\bfu_0=0$ and $\nabla\bfd_0=0$ a.e.\ in $\cF_0$, so the stationary body together with $\bfd(t)=\bfd_0|_{\cF_0}$ is a constrained solution.  
We thus assume $E_0>0$ in the following.

Fix $T<T_*$. 
For every sufficiently small $\eps$, \autoref{thm:weak sol GL} applied with anchor $\bfd_0^\eps$ gives a solution $(\bfu_\eps,\bfd_\eps,\varphi_\eps)$ on $[0,T]$. 
Due to the convergence of the energy $E_0^\eps\to E_0$, \autoref{lem:dissipative rigid displacement} gives a number $\rho_T>0$, independent of small $\eps$, such that
\begin{equation}\label{eq:eps uniform separation}
 \dist(\cS_\eps(t),\partial\Omega)\ge\rho_T
 \qquad(0\le t\le T).
\end{equation}
Define $\bftau_\eps:=\Delta\bfd_\eps+\frac1{\eps^2}(1-|\bfd_\eps|^2)\bfd_\eps$.
The energy inequality \eqref{eq:energy-main GL} for $(\bfu_\eps,\bfd_\eps,\varphi_\eps)$, whose right-hand side equals $E_0^\eps$ because $|\bfd_0^\eps|=1$, together with \eqref{eq:prepared anchors eps} and Korn's inequality, implies
\begin{equation}\label{eq:eps uniform energy}
\begin{aligned}
 \sup_{0<\eps<\eps_T}\Bigg[&\|\bfu_\eps\|_{\rL^\infty(0,T;\rL^2(\Omega))}^2
 +\|\bfu_\eps\|_{\rL^2(0,T;\rH^1_0(\Omega))}^2
 +\|\nabla\bfd_\eps\|_{\rL^\infty(0,T;\rL^2(\cF_\eps(t)))}^2\\
 &\qquad+\int_0^T\int_{\cF_\eps(t)}\frac{(1-|\bfd_\eps|^2)^2}{\eps^2}\srd\bfx\srd t
 +\|\bftau_\eps\|_{\rL^2(\cQ_{\cF_\eps})}^2\Bigg]\le C_T.
\end{aligned}
\end{equation}

\subsection{Rigid motion and velocity compactness}

\begin{lem}\label{lem:eps director bounds}
For every sufficiently small $\eps$,
\[
 |\bfd_\eps|\le1,\qquad (\bfd_\eps)_3\ge0
 \quad\text{a.e.\ in }\cQ_{\cF_\eps}.
\]
\end{lem}
\begin{proof}
By construction, see \autoref{thm:weak sol GL}, the pair $(\bfu_\eps,\bfd_\eps)$ is obtained from solutions of the semi-Galerkin system of \autoref{ssec:semi-Galerkin approx} with anchor $\bfd_0^\eps$ through the limits $k\to\infty$, $n\to\infty$, and $\delta\to0$. The Galerkin director equation \eqref{eq:k-director} is of the form \eqref{eq:penalized maximum principle equation} with
\[
 \bfu=\bfu^k,\qquad a=n\varphi^k\ge0,\qquad \bb=\bb^k=\bfd_0^\eps\circ\bfY^k,
\]
and its initial and boundary datum is $\bfd_0^\eps$. Since $\bfd_0^\eps$ takes values in $\bS_+^2$, \autoref{lem:relaxed dir cond} gives $|\bfd^k|\le1$ and $(\bfd^k)_3\ge0$ a.e.\ in $(0,T)\times\Omega$. These closed convex constraints are preserved by the strong~$\rL^2$-convergence of the directors in each of the three limits, see \eqref{eq:delta director convex constraint} for the post-rigidity level.
\end{proof}

Let $\bfh_\eps$ and $\bfO_\eps$ be the rigid motions. 
From \autoref{lem:rigid velocity energy control} and \eqref{eq:eps uniform energy}, we infer uniform $\rW^{1,\infty}$-bounds.  
After extraction, as $\eps \to 0$, we have
\begin{equation}\label{eq:rigid motion eps convergence}
 \bfh_\eps\to\bfh,\qquad \bfO_\eps\to\bfO
 \quad\text{uniformly on }[0,T],
\end{equation}
and, with $\varphi=\dsone_{\bfh(t)+\bfO(t)\cS_0}$ capturing the position of the rigid body, we get the convergence
\begin{equation}\label{eq:phi eps strong}
    \varphi_\eps\to\varphi \tin \rC([0,T];\rL^p(\Omega)),\tfor 1\le p<\infty.
\end{equation}
This convergence is due to the fact that corresponding material points of $\cS_\eps(t)$ and $\cS(t)$ have distance at most
\[
 \eta_\eps:=\|\bfh_\eps-\bfh\|_{\rL^\infty(0,T)}+\operatorname{diam}(\cS_0)\|\bfO_\eps-\bfO\|_{\rL^\infty(0,T)}\to0.
\]
Hence the solids converge uniformly in Hausdorff distance, and their symmetric difference is contained in the layer $\{\bfx:\dist(\bfx,\partial\cS(t))\le2\eta_\eps\}$ around $\partial\cS(t)$, a rigid image of the $\rC^2$-surface $\partial\cS_0$, whose measure is at most $C\eta_\eps$. Consequently,
\[
 \sup_{t\in[0,T]}\|\varphi_\eps(t)-\varphi(t)\|_{\rL^p(\Omega)}^p=\sup_{t\in[0,T]}|\cS_\eps(t)\triangle\cS(t)|\le C\eta_\eps.
\]

We apply \autoref{lem:FSI compactness criterion} to $\bfu_\eps$ and the trajectories $\cS_\eps(t)$. The separation, convergence, and boundedness hypotheses hold by \eqref{eq:eps uniform separation}, \eqref{eq:rigid motion eps convergence}, and \eqref{eq:eps uniform energy}. To verify \eqref{eq:FSI temporal functional bound}, let $\bfzeta$ be a test field as in the lemma, rigid on $\cS(t)_\sigma$. By the uniform Hausdorff convergence of the solids, $\cS_\eps(t)\subset\cS(t)_\sigma$ for all small $\eps$, so $\bfzeta$ is admissible in the momentum equation of \autoref{def:weak sol GL}. By the Sobolev embedding of~$\rH^3(\Omega)$ into $\rW^{1,\infty}(\Omega)$, the convective term and the Ericksen stress are bounded by~$C_T\|\bfzeta\|_{\rL^\infty(I;\rH^3(\Omega))}$ and the viscous term by $C_T\|\bfzeta\|_{\rL^2(I;\rH^3(\Omega))}$, uniformly in $\eps$ by \eqref{eq:eps uniform energy}. 

Consequently, after extraction,
\begin{equation}\label{eq:u eps compact}
\begin{aligned}
 \bfu_\eps&\to\bfu &&\text{strongly in }\rL^2((0,T)\times\Omega),\\
 \bfu_\eps&\rightharpoonup\bfu &&\text{weakly in }\rL^2(0,T;\rH^1_0(\Omega)),\\
 \bfu_\eps&\stackrel{*}{\rightharpoonup}\bfu &&\text{weak-* in }\rL^\infty(0,T;\rL^2(\Omega)).
\end{aligned}
\end{equation}
In the same way as in the limit passage $\delta \to 0$ in \eqref{eq:rigid body vel after delta lim}, for the limit $\eps \to 0$, we can show that
\[
 \bfu(t,\bfx)=\bfell(t)+\bfomega(t)\times(\bfx-\bfh(t))
 \quad\text{for a.e.\ }\bfx\in\cS(t),
\]
so the limiting rigid trajectory is generated by the rigid restriction of the limiting velocity, including its orientation.

The material extension through the rigid phase also allows us to deduce global compactness of the director.  
Define
\begin{equation*}
    \overline\bfd_\eps(t,\bfx) \coloneqq
    \left\{
    \begin{aligned}
         &\bfd_\eps(t,\bfx),&&\tfor \bfx\in\cF_\eps(t),\\
         &\bfd_0^\eps\left(\bfO_\eps(t)^\top(\bfx-\bfh_\eps(t))\right),&&\tfor \bfx\in\cS_\eps(t).
    \end{aligned}
     \right.
\end{equation*}
The anchoring condition in \autoref{def:weak sol GL} ensures that the two traces agree on the moving interface, so \autoref{lem:material gluing} implies that $\overline\bfd_\eps(t)\in\rH^1(\Omega;\R^3)$.
The fact that the change of variables $\bfx=\bfh_\eps(t)+\bfO_\eps(t)\bfy$ is an isometry of $\cS_0$ onto $\cS_\eps(t)$ leads to
\[
    \|\nabla\overline\bfd_\eps(t)\|_{\rL^2(\Omega)}^2 =\|\nabla\bfd_\eps(t)\|_{\rL^2(\cF_\eps(t))}^2+\|\nabla\bfd_0^\eps\|_{\rL^2(\cS_0)}^2 .
\]
The first term is bounded by \eqref{eq:eps uniform energy} and the second by \eqref{eq:prepared anchors eps}. 
Joint with $|\overline\bfd_\eps|\le1$, this yields
\[
    \sup_{\eps<\eps_T}\|\overline\bfd_\eps\|_{\rL^\infty(0,T;\rH^1(\Omega))}\le C_T.
\]
Let $\widetilde\bftau_\eps$ be the zero extension of $\bftau_\eps$ to the solid. 
An application of \autoref{lem:material gluing} with anchor $\bfd_0^\eps$ gives
\[
    \partial_t\overline\bfd_\eps+(\bfu_\eps\cdot\nabla)\overline\bfd_\eps=\widetilde\bftau_\eps \tin \cD'((0,T)\times\Omega).
\]
By $\bfu_\eps\in\rL^2_t\rL^6_{\bfx}$ and~$\nabla\overline\bfd_\eps\in\rL^\infty_t\rL^2_{\bfx}$, $(\bfu_\eps\cdot\nabla)\overline\bfd_\eps$ is bounded in $\rL^2(0,T;\rL^{3/2}(\Omega))\hookrightarrow \rL^2(0,T;\rH^{-1}(\Omega))$.
On the other hand, $\widetilde\bftau_\eps$ is bounded in $\rL^2_{t,\bfx}$.  
Hence, by the Aubin--Lions lemma, after extracting a subsequence we have
\begin{equation}\label{eq:eps global material strong L2}
 \overline\bfd_\eps\to\overline\bfd
 \quad\text{strongly in }\rL^2((0,T)\times\Omega),
 \qquad
 \overline\bfd_\eps\stackrel{*}{\rightharpoonup}\overline\bfd
 \quad\text{in }\rL^\infty(0,T;\rH^1(\Omega)).
\end{equation}
On every compact subset of the limiting fluid spacetime, the restriction of $\overline\bfd$ is the director
identified below by the local Lin--Wang argument. On the solid, it is the rigidly transported limit of
$\bfd_0^\eps$, namely~$\bfd_0(\bfO(t)^\top(\bfx-\bfh(t)))$.

\subsection{A quantitative Lin--Wang compactness principle}

The singular director limit is localized to compact subsets of the moving fluid region. The required static compactness statement is the following.

\begin{lem}\label{lem:LW compactness input}
Let $U\Subset\R^3$ be smooth and bounded and let $V\Subset U$.  Suppose $e_j\to0$ and
$\bfv_j\in \rH^1(U;\R^3)$ satisfies
\begin{equation}\label{eq:LW static equation}
 \Delta\bfv_j+\frac1{e_j^2}(1-|\bfv_j|^2)\bfv_j=\bftau_j
 \quad\text{in }U,
\end{equation}
with
\[
 |\bfv_j|\le1,\qquad (\bfv_j)_3\ge0,
\]
\begin{equation}\label{eq:LW static bounds}
 \sup_j\int_U\left(\frac12|\nabla\bfv_j|^2+
 \frac{(1-|\bfv_j|^2)^2}{4e_j^2}\right)\srd\bfx<\infty,
 \qquad
 \sup_j\|\bftau_j\|_{\rL^2(U)}<\infty.
\end{equation}
Then $(\bfv_j)$ is precompact in $\rH^1(V)$.  Every limit belongs to
$\rH^1(V;\bS^2_+)$.
\end{lem}

\begin{proof}
This is a direct local application of \cite[Theorem~6.1]{LW:16}.  In the
notation of that theorem,~\mbox{\eqref{eq:LW static equation}--\eqref{eq:LW static bounds}}
place $\bfv_j$ in the class $\mathbf X(L_1,L_2,a;U)$ with $a=1$. Indeed,
$|\bfv_j|\le1$, and $(\bfv_j)_3\ge0$ is the condition $(\bfv_j)_3\ge-1+a$ for $a=1$,
while the two remaining defining bounds are the uniform bound for the Ginzburg--Landau energy of $\bfv_j$ and the uniform $\rL^2(U)$-bound for its tension $\bftau_j$, both of which are assumed in \eqref{eq:LW static bounds}.
The cited theorem states that $\mathbf X(L_1,L_2,a;U)$ is precompact
in~$\rH^1_{\rm loc}(U;\R^3)$.  Hence $(\bfv_j)$ is precompact in
$\rH^1(V)$ for every $V\Subset U$.

If $\bfv_j\to\bfv$ strongly in $\rH^1(V)$ along a convergent subsequence, then
\[
 \|1-|\bfv_j|^2\|_{\rL^2(V)}^2
 \le 4e_j^2\sup_i\int_U
 \left(\frac12|\nabla\bfv_i|^2+
 \frac{(1-|\bfv_i|^2)^2}{4e_i^2}\right)\srd\bfx\rightarrow0.
\]
Thus $|\bfv|=1$ a.e.\ on $V$.  The closed half-space condition
$(\bfv_j)_3\ge0$ passes under the strong~$\rL^2$-convergence, so
$\bfv\in \rH^1(V;\bS^2_+)$.
\end{proof}

\begin{rem}
Note that \autoref{lem:LW compactness input} only provides compactness, and it does not identify the weak limit of the tensions $\bftau_j$. 
In fact, the limiting tension is identified separately, see the cross-product identity and the tangency relation \eqref{eq:tension tangent limit} in the proof of \autoref{prop:eps local strong d}.
\end{rem}

The following lemma establishes a uniform local bound of director fields in $\rH^1$.

\begin{lem}\label{lem:LW compactness modulus}
Consider smooth bounded $V\Subset U\Subset\R^3$ as well as numbers $L,\Lambda>0$.  
For every $\eta>0$, there exists~$r=r(U,V,L,\Lambda,\eta)>0$ such that if $0<e_1,e_2\le1$ and $\bfv_i\in \rH^1(U;\R^3)$ satisfy
\[
 \Delta\bfv_i+\frac1{e_i^2}(1-|\bfv_i|^2)\bfv_i=\bftau_i,
 \quad |\bfv_i|\le1,\quad (\bfv_i)_3\ge0,
\]
\[
 \int_U\left(\frac12|\nabla\bfv_i|^2+
 \frac{(1-|\bfv_i|^2)^2}{4e_i^2}\right)\srd\bfx\le L,
 \qquad \|\bftau_i\|_{\rL^2(U)}\le\Lambda,
\]
and $\|\bfv_1-\bfv_2\|_{\rL^2(U)}\le r$, then $\|\bfv_1-\bfv_2\|_{\rH^1(V)}\le\eta$.
\end{lem}
\begin{proof}
Suppose by contradiction that the assertion fails. Then there are pairs $(\bfv_{1,j},\bfv_{2,j})$ with parameters $e_{1,j},e_{2,j}\in(0,1]$ such that both maps satisfy the above bounds,
$\|\bfv_{1,j}-\bfv_{2,j}\|_{\rL^2(U)}\to0$, but
$\|\bfv_{1,j}-\bfv_{2,j}\|_{\rH^1(V)}>\eta$ for every $j$.
upon passing to a subsequence, for $i=1,2$, we may assume $e_{i,j}\rightarrow \bar e_i\in[0,1]$.
Next, we choose an intermediate smooth domain $W$ such that $V\Subset W\Subset U$. 
We claim that, for each $i$, the sequence $\{\bfv_{i,j}\}_j$ is precompact in $\rH^1(V)$.

If $\bar e_i=0$, then $(\bfv_{i,j})_j$ satisfies the hypotheses of \autoref{lem:LW compactness input} with $e_j=e_{i,j}\to0$ and $\bftau_j=\bftau_{i,j}$. 
In fact, \eqref{eq:LW static equation} and the constraints $|\bfv_{i,j}|\le1$, $(\bfv_{i,j})_3\ge0$ hold by assumption, and the bounds \eqref{eq:LW static bounds} hold with the constants $L$ and $\Lambda$, showing precompactness of the sequence in $\rH^1(V)$.

On the other hand, if $\bar e_i>0$, then $e_{i,j}\ge \bar e_i/2$ for all sufficiently large $j$. 
In
\[
 \Delta\bfv_{i,j}
 =\bftau_{i,j}-\frac1{e_{i,j}^2}(1-|\bfv_{i,j}|^2)\bfv_{i,j}
\]
the bound $|\bfv_{i,j}|\le1$ yields that the right-hand side is uniformly bounded in $\rL^2(U)$.
By the standard interior elliptic estimates we get a uniform $\rH^2(W)$-bound. 
From $\rH^2(W)\Subset\rH^1(V)$, it follows that the sequence is precompact in $\rH^1(V)$.

In summary, after passing to a common subsequence,
$\bfv_{i,j}\to\bfv_i$ strongly in $\rH^1(V)$ for $i=1,2$. 
The $\rL^2(U)$-distance tending to zero forces $\bfv_1=\bfv_2$ a.e.\ on $V$, and hence the $\rH^1(V)$-distance also tends to zero, a contradiction.
\end{proof}

\begin{lem}\label{lem:LW spacetime Cauchy}
Let $I$ be a bounded interval and $V\Subset U$ as in \autoref{lem:LW compactness modulus}. Let $e_j\in(0,1]$, and let~$\bfd_j\in\rL^2(I;\rH^1(U;\R^3))$ and $\bftau_j\in\rL^2(I\times U;\R^3)$ satisfy, for a.e.\ $t\in I$,
\begin{equation}\label{eq:LW slice equation}
 \Delta\bfd_j(t)+\frac1{e_j^2}\bigl(1-|\bfd_j(t)|^2\bigr)\bfd_j(t)=\bftau_j(t)\quad\text{in }U,
 \qquad |\bfd_j(t)|\le1,\qquad (\bfd_j(t))_3\ge0,
\end{equation}
\begin{equation}\label{eq:LW slice bounds}
 \int_U\left(\frac12|\nabla\bfd_j(t)|^2+\frac{(1-|\bfd_j(t)|^2)^2}{4e_j^2}\right)\srd\bfx\le L,
 \qquad
 \int_I\|\bftau_j(t)\|_{\rL^2(U)}^2\srd t\le C.
\end{equation}
If $(\bfd_j)$ is a Cauchy sequence in $\rL^2(I\times U)$, then it is a Cauchy sequence in $\rL^2(I;\rH^1(V))$.
\end{lem}
\begin{proof}
Fix $\eta>0$ and $\Lambda>0$, and let $r$ be given by \autoref{lem:LW compactness modulus}.  Define
\[
 B_{j,\Lambda}:=\{t\in I:\|\bftau_j(t)\|_{\rL^2(U)}>\Lambda\},
 \qquad
 A_{ij,r}:=\{t\in I:\|\bfd_i(t)-\bfd_j(t)\|_{\rL^2(U)}>r\}.
\]
By Chebyshev's inequality we have
\[
 |B_{j,\Lambda}|\le C\Lambda^{-2},\qquad
 |A_{ij,r}|\le r^{-2}\|\bfd_i-\bfd_j\|_{\rL^2(I\times U)}^2.
\]
We now address the set $I\setminus(B_{i,\Lambda}\cup B_{j,\Lambda}\cup A_{ij,r})$.
In fact, for a.e.\ $t\in I\setminus(B_{i,\Lambda}\cup B_{j,\Lambda}\cup A_{ij,r})$, \eqref{eq:LW slice equation} and \eqref{eq:LW slice bounds} imply that $\bfd_i(t)$ and $\bfd_j(t)$ satisfy the hypotheses of \autoref{lem:LW compactness modulus}.
Thus, $\|\bfd_i(t)-\bfd_j(t)\|_{\rH^1(V)}\le\eta$.  
The energy bound in~\eqref{eq:LW slice bounds} gives $ \int_V|\nabla\bfd_i-\nabla\bfd_j|^2\srd\bfx\le8L$ and $|\bfd_i-\bfd_j|\le2$.
In total, we get
\[
    \limsup_{i,j\to\infty}\|\bfd_i-\bfd_j\|_{\rL^2(I;\rH^1(V))}^2 \le |I|\eta^2+C_{U,L}\Lambda^{-2}.
\]
The proof is completed upon letting $i,j\to\infty$, then $\Lambda\to\infty$, and finally $\eta\to0$.
\end{proof}

\begin{prop}\label{prop:eps local strong d}
There exist a director field $\bfd$ and a tension field $\bftau\in \rL^2(\cQ_\cF;\R^3)$ such that, after extraction, for every compact~$K\Subset\cQ_\cF$, we get
\begin{equation}\label{eq:local strong d eps}
    \bfd_\eps\to\bfd\quad\text{strongly in }\rL^2(K), \qquad \nabla\bfd_\eps\to\nabla\bfd\quad\text{strongly in }\rL^2(K),
\end{equation}
and $\bftau_\eps\rightharpoonup\bftau$ weakly in $\rL^2(K)$.  
Moreover, we have $|\bfd|=1$, $d_3\ge0$, and $\bftau=\Delta\bfd+|\nabla\bfd|^2\bfd$ in $\cD'(\cQ_\cF)$.
\end{prop}
\begin{proof}
Let $K\Subset\cQ_\cF$.  A finite covering gives cylinders $I_j\times V_j$ and smooth domains $V_j\Subset U_j$ such that $U_j\subset\cF_\eps(t)\cap\cF(t)$ for $t\in I_j$ and all small $\eps$.  On $I_j\times U_j$,
\[
 \partial_t\bfd_\eps=\bftau_\eps-(\bfu_\eps\cdot\nabla)\bfd_\eps.
\]
By \eqref{eq:eps uniform energy}, $\bftau_\eps$ is bounded in $\rL^2(I_j\times U_j)$, $\bfu_\eps$ is bounded in $\rL^2(I_j;\rL^6(U_j))$ by the Sobolev embedding $\rH_0^1(\Omega)\hookrightarrow\rL^6(\Omega)$, and $\nabla\bfd_\eps$ is bounded in $\rL^\infty(I_j;\rL^2(U_j))$.
H\"older's inequality then yields that $(\bfu_\eps\cdot\nabla)\bfd_\eps$ is bounded in~$\rL^2(I_j;\rL^{3/2}(U_j))$, so in total, the right-hand side is bounded in $\rL^2(I_j;\rW^{-1,3/2}(U_j))$.
Moreover, due to $|\bfd_\eps|\le1$, $\bfd_\eps$ is bounded in $\rL^2(I_j;\rH^1(U_j))$.
By $\rH^1(U_j)\Subset \rL^2(U_j)\hookrightarrow \rW^{-1,3/2}(U_j)$, applying the Aubin--Lions lemma, after extracting a subsequence we have
\[
 \bfd_\eps\to\bfd\quad\text{strongly in }\rL^2(I_j\times U_j).
\]
For a.e.\ $t\in I_j$, $\bfd_\eps(t)$ satisfies \eqref{eq:LW slice equation} on $U_j$ with $e_j=\eps$ and tension $\bftau_\eps(t)$, by the definition of $\bftau_\eps$ and the bounds $|\bfd_\eps|\le1$, $(\bfd_\eps)_3\ge0$.
On the other hand, the bounds \eqref{eq:LW slice bounds} follow from \eqref{eq:eps uniform energy}.  
\autoref{lem:LW spacetime Cauchy} then upgrades the previous convergence to
\[
    \bfd_\eps\to\bfd\quad\text{strongly in }\rL^2(I_j;\rH^1(V_j)).
\]
Invoking a finite covering, we conclude \eqref{eq:local strong d eps}. 
The estimate $\int_{I_j}\int_{U_j}(1-|\bfd_\eps|^2)^2\srd\bfx\srd t\le C_T\eps^2$, which follows from \eqref{eq:eps uniform energy}, together with the strong convergence gives $|\bfd|=1$ a.e., while $(\bfd_\eps)_3\ge0$ yields $d_3\ge0$.

To identify the tension, we take the cross product of the elliptic equation with $\bfd_\eps$ and obtain
\[
 \bfd_\eps\times\bftau_\eps
 =\bfd_\eps\times\Delta\bfd_\eps
 =\operatorname{div}(\bfd_\eps\times\nabla\bfd_\eps).
\]
By \eqref{eq:eps uniform energy}, the tensions $\bftau_\eps$ are bounded in $\rL^2(K)$ for every compact $K\Subset\cQ_\cF$, uniformly in $K$. 
Hence, after extraction, we infer the existence of $\bftau\in\rL^2(\cQ_\cF;\R^3)$ such that $\bftau_\eps\rightharpoonup\bftau$ weakly in $\rL^2(K)$ for every such~$K$. Testing the above identity against smooth compactly supported fields and using the strong convergence \eqref{eq:local strong d eps} together with this weak convergence gives, in the limit,
\begin{equation}\label{eq:cross product tension limit}
    \bfd\times\bftau=\operatorname{div}(\bfd\times\nabla\bfd) \tin \cD'(\cQ_\cF).
\end{equation}
We now fix a cylinder $I\times U$ with $\overline I\times\overline U\Subset\cQ_{\cF}$. 
Concatenating $|\bfd|=1$, $\bftau\in\rL^2(\cQ_\cF)$, and $\nabla\bfd\in\rL^2(I\times U)$ by~\eqref{eq:local strong d eps}, we find that $\bfd\times\bftau$ and $\bfd\times\nabla\bfd$ belong to $\rL^2(I\times U)$.
As \eqref{eq:cross product tension limit} contains no time derivative, it extends by density to
\begin{equation}\label{eq:cross product weak}
 \int_I\int_U(\bfd\times\bftau)\cdot\boldsymbol\chi\srd\bfx\srd t
 =-\sum_{i=1}^3\int_I\int_U(\bfd\times\partial_i\bfd)\cdot\partial_i\boldsymbol\chi\srd\bfx\srd t
 \tfor \boldsymbol\chi\in\rL^2(I;\rH_0^1(U;\R^3)).
\end{equation}

Next, we show that the tension is tangential. 
The equation $\partial_t\bfd_\eps+(\bfu_\eps\cdot\nabla)\bfd_\eps=\bftau_\eps$ gives
\[
 \partial_t\frac{|\bfd_\eps|^2}{2}
 +\bfu_\eps\cdot\nabla\frac{|\bfd_\eps|^2}{2}
 =\bfd_\eps\cdot\bftau_\eps.
\]
Let $\phi\in \rC_c^\infty(\cQ_\cF)$. Since $\supp\phi$ is a compact subset of $\cQ_\cF$ and the solids $\cS_\eps(t)$ converge to $\cS(t)$ uniformly in Hausdorff distance, see \eqref{eq:rigid motion eps convergence}, we have $\supp\phi\subset\cQ_{\cF_\eps}$ for all sufficiently small $\eps$, and integration by parts yields
\[
 \int_0^T\int_\Omega \phi\,\bfd_\eps\cdot\bftau_\eps \srd \bfx \srd t
 =-\frac12\int_0^T\int_\Omega |\bfd_\eps|^2\bigl(\partial_t\phi+\bfu_\eps\cdot\nabla\phi\bigr) \srd \bfx \srd t,
\]
where we used $\mdiv\bfu_\eps=0$. The local strong $\rL^2$-convergence of $\bfd_\eps$, the bound $|\bfd_\eps|\le1$, and~$|\bfd|=1$ imply $|\bfd_\eps|^2\to1$ strongly in $\rL^p(\supp\phi)$ for every $p<\infty$. Together with the strong local $\rL^2$-convergence of $\bfu_\eps$ and the weak $\rL^2$-convergence of $\bftau_\eps$, this yields
\[
 \int_0^T\int_\Omega \phi\,\bfd\cdot\bftau \srd \bfx \srd t
 =-\frac12\int_0^T\int_\Omega\bigl(\partial_t\phi+\bfu\cdot\nabla\phi\bigr) \srd \bfx \srd t=0,
\]
because $\phi$ has compact support in spacetime and $\mdiv\bfu=0$. 
This implies
\begin{equation}\label{eq:tension tangent limit}
 \bfd\cdot\bftau=0\qquad\text{a.e.\ in }\cQ_\cF.
\end{equation}

Finally, let $\bfphi\in \rC_c^\infty(I\times U;\R^3)$, and set $\bfphi_T:=\bfphi-(\bfd\cdot\bfphi)\bfd$ and $\boldsymbol\chi:=\bfd\times\bfphi_T=\bfd\times\bfphi$. With regard to $\bfd\in\rL^2(I;\rH^1(U))$ and $|\bfd|=1$, both fields belong to $\rL^2(I;\rH_0^1(U;\R^3))$, and no time regularity is needed.
As $|\bfd|=1$, $\bfd\cdot\partial_i\bfd=0$, $\bfd\cdot\bfphi_T=0$, and \eqref{eq:tension tangent limit} holds, we infer that a.e.
\begin{equation*}
    \begin{aligned}
        \bftau\cdot\bfphi=\bftau\cdot\bfphi_T
        &=(\bfd\times\bftau)\cdot(\bfd\times\bfphi_T)=(\bfd\times\bftau)\cdot\boldsymbol\chi, \tand\\
        (\bfd\times\partial_i\bfd)\cdot\partial_i\boldsymbol\chi
        &=(\bfd\times\partial_i\bfd)\cdot(\partial_i\bfd\times\bfphi_T) +(\bfd\times\partial_i\bfd)\cdot(\bfd\times\partial_i\bfphi_T)=\partial_i\bfd\cdot\partial_i\bfphi_T .
    \end{aligned}
\end{equation*}
Inserting $\boldsymbol\chi$ into \eqref{eq:cross product weak} therefore gives
\[
 \int_I\int_U\bftau\cdot\bfphi \srd \bfx \srd t
 =-\int_I\int_U\nabla\bfd:\nabla\bfphi_T \srd \bfx \srd t.
\]
Since $\bfd\cdot\partial_i\bfd=0$, we have $\nabla\bfd:\nabla\bigl((\bfd\cdot\bfphi)\bfd\bigr)=|\nabla\bfd|^2(\bfd\cdot\bfphi)$, and consequently
\begin{equation}\label{eq:id of tensor field}
    \int_I\int_U\bftau\cdot\bfphi \srd \bfx \srd t=-\int_I\int_U\nabla\bfd:\nabla\bfphi \srd \bfx \srd t+\int_I\int_U|\nabla\bfd|^2\bfd\cdot\bfphi \srd \bfx \srd t.
\end{equation}
As the cylinder $I\times U$ was arbitrary, $\bftau=\Delta\bfd+|\nabla\bfd|^2\bfd$ in $\cD'(\cQ_\cF)$.
\end{proof}

Let $\widetilde\bftau_\eps$ denote the extension of $\bftau_\eps$ by zero to the solid. By \eqref{eq:eps uniform energy}, a subsequence converges weakly in $\rL^2((0,T)\times\Omega)$ to some $\widetilde\bftau$. By the uniform convergence \eqref{eq:rigid motion eps convergence} of the rigid motions, every compact subset of $\cQ_\cF$ lies in $\cQ_{\cF_\eps}$, and every compact subset of the open set $\{(t,\bfx):\bfx\in\cS(t)\}$ lies in~$\{(t,\bfx):\bfx\in\cS_\eps(t)\}$, for all sufficiently small $\eps$. Hence $\widetilde\bftau=\bftau$ in $\cQ_\cF$ by \autoref{prop:eps local strong d}, and $\widetilde\bftau=0$ a.e.\ in the solid region, since $\widetilde\bftau_\eps$ vanishes there. 
Due to the measure zero of the interface, $\widetilde\bftau$ is the zero extension of $\bftau$, and the weak lower semicontinuity of the norm gives the estimate
\begin{equation}\label{eq:tau global liminf}
    \int_0^T\int_{\cF(t)}|\bftau|^2 \srd \bfx \srd t \le\liminf_{\eps\to0}\int_0^T\int_{\cF_\eps(t)}|\bftau_\eps|^2 \srd \bfx \srd t.
\end{equation}

\subsection{Anchoring and passage to the weak equations}

In this subsection, we show that the limiting $\bfd$ satisfies the anchoring condition \eqref{eq:weak trace anchoring} with $\bfd_0$, and we pass to the limit in the director and fluid equations.

\begin{lem}\label{lem:anchoring eps limit}
The limiting director field $\bfd$ satisfies the anchoring condition \eqref{eq:weak trace anchoring} with the original $\rH^1$-datum $\bfd_0$.
\end{lem}
\begin{proof}
Let $r_0>0$ be such that the nearest-point projections onto $\partial\Omega$ and $\partial\cS_0$ are well defined and of class $\rC^1$ on the $r_0$-neighborhoods of these surfaces.
The existence of such an $r_0$ is guaranteed by the compactness and $\rC^2$-regularity of both surfaces.
We now choose $0<r<r_0$ smaller than one quarter of the uniform separation $\rho_T$ in~\eqref{eq:eps uniform separation}, establishing that $U_{\rm out}$ and $U_{\rm in}^0$ defined below are Lipschitz domains, and thereby showing that the trace operators onto $\partial\Omega$ and $\partial\cS_0$ are continuous on $\rH^1$ of $U_{\rm out}$ and $U_{\rm in}^0$.
Another consequence is that $U_{\rm out}$ does not meet any of the solids~$\cS_\eps(t)$.

On $U_{\rm out}=\{\bfx\in\Omega:\dist(\bfx,\partial\Omega)<r\}$, the material extensions $\overline\bfd_\eps$ coincide with the fluid directors.  
The convergences from \eqref{eq:eps global material strong L2} and the uniform $\rL^\infty_t\rH^1_{\bfx}$-bound imply that the restrictions of $\overline\bfd_\eps$ converge strongly in $\rL^2_{t,\bfx}$ and weakly in~$\rL^2_t\rH^1_{\bfx}$.  
The corresponding outer traces are $\bfd_0^\eps$, and they converge to~$\bfd_0$ in $\rH^{1/2}(\partial\Omega)$ by \eqref{eq:prepared anchors eps}. 
The outer anchoring condition then follows from the continuity of the trace operator.

For the inner boundary, we introduce $U_{\rm in}^0=\{\bfy\in\cF_0:\dist(\bfy,\partial\cS_0)<r\}$ along with the rigid maps $\bfZ_\eps(t,\bfy)=\bfh_\eps(t)+\bfO_\eps(t)\bfy$. 
The pullbacks $\widehat\bfd_\eps(t,\bfy)\coloneqq\overline\bfd_\eps(t,\bfZ_\eps(t,\bfy))$ are uniformly bounded in $\rL^\infty(0,T;\rH^1(U_{\rm in}^0))$.  
From \eqref{eq:eps global material strong L2} and \eqref{eq:rigid motion eps convergence}, we deduce that
$\widehat\bfd_\eps\to\overline\bfd(t,\bfh(t)+\bfO(t)\bfy)$ strongly in~$\rL^2((0,T)\times U_{\rm in}^0)$. 
The weak $\rL^2_t\rH^1_{\bfy}$-limit is thus uniquely the indicated pullback of $\overline\bfd$.  
On $\partial\cS_0$ the pulled-back traces equal $\bfd_0^\eps$, so their weak trace limit is $\bfd_0$. 
Invoking the push-forward of this identity by the limiting rigid map, we obtain the inner anchoring in
\eqref{eq:weak trace anchoring}.
\end{proof}

Let now $\bfzeta\in\cT_{\bfu}^{\rm col}$, where we recall the test space $\cT_{\bfu}^{\rm col}$ from \eqref{eq:collar momentum tests}.
Its symmetric gradient vanishes in a fixed collar of~$\cS(t)$, which contains $\cS_\eps(t)$ for all sufficiently small $\eps$.  
By symmetry of the Ericksen tensor, its contraction with~$\nabla\bfzeta$ vanishes in that collar.  
On the complement, local strong $\rH^1$-convergence from \autoref{prop:eps local strong d} allows for the limit passage in the stress. 
Note that this convergence does not cover a neighborhood of $t=0$. 
However, the bound for $\nabla\bfd_\eps$ in $\rL^\infty(0,T;\rL^2(\cF_\eps(t)))$ from \eqref{eq:eps uniform energy} yields that the contribution of the time layer $(0,t_0)$ to the stress term is at most $Ct_0\|\nabla\bfzeta\|_{\rL^\infty((0,T)\times\Omega)}$, uniformly in $\eps$.
The same bound holds for the limit. 
Thus, letting first $\eps\to0$ and then $t_0\to0$ gives the momentum identity for every collar-rigid test.
\autoref{prop:rigid test equivalence} then extends the identity to every $\bfzeta\in\cT_{\bfu}$, proving \eqref{eq:weak sol u}.

For the director, we first consider $\bfpsi\in\cT_{\bfd}$, with $\cT_{\bfd}$ as introduced in \eqref{eq:smooth loc test class dir}.  
For all sufficiently small $\eps$, it is supported in the approximate fluid region, and the molecular-field equation gives
\begin{equation*}
 \int_0^T\int_{\cF_\eps(t)}
 \left[-\bfd_\eps\cdot\partial_t\bfpsi
 +(\bfu_\eps\cdot\nabla)\bfd_\eps\cdot\bfpsi
 -\bftau_\eps\cdot\bfpsi\right] \srd\bfx\srd t
 =\int_{\cF_0}\bfd_0^\eps\cdot\bfpsi(0) \srd\bfx.
\end{equation*}
On every interval $[t_0,T]$, the limit follows from \eqref{eq:u eps compact} and \autoref{prop:eps local strong d}, while on $(0,t_0)$, we integrate the convection term by parts in space and obtain
\[
 \int_0^{t_0}\int_{\cF_\eps(t)} (\bfu_\eps\cdot\nabla)\bfd_\eps\cdot\bfpsi\srd\bfx\srd t
 =-\int_0^{t_0}\int_{\cF_\eps(t)} (\bfu_\eps\cdot\nabla)\bfpsi\cdot\bfd_\eps\srd\bfx\srd t.
\]
Thus, the bound $|\bfd_\eps|\le1$ from \autoref{lem:eps director bounds} and the bound for $\bfu_\eps$ in $\rL^\infty(0,T;\rL^2(\Omega))$ from \eqref{eq:eps uniform energy} make it~$O(t_0)$. 
By the Cauchy--Schwarz inequality in time and the bound for $\bftau_\eps$ in \eqref{eq:eps uniform energy}, the molecular-field term is $O(t_0^{1/2})$, and the term containing $\partial_t\bfpsi$ is $O(t_0)$ since $|\bfd_\eps|\le1$.
Letting first $\eps\to0$ and then~$t_0\to0$, and using \eqref{eq:prepared anchors eps}, yields
\begin{equation}\label{eq:director molecular limit eps}
 \int_0^T\int_{\cF(t)}
 \left[-\bfd\cdot\partial_t\bfpsi
 +(\bfu\cdot\nabla)\bfd\cdot\bfpsi
 -\bftau\cdot\bfpsi\right] \srd\bfx\srd t
 =\int_{\cF_0}\bfd_0\cdot\bfpsi(0) \srd\bfx.
\end{equation}
The functional in \eqref{eq:director molecular limit eps} is continuous on $\cX_{\bfd}(\varphi)$, since the only nontrivial term $(\bfu\cdot\nabla)\bfd\cdot\bfpsi$ is bounded by $\|\bfu\|_{\rL^2_t\rL^6_{\bfx}}\|\nabla\bfd\|_{\rL^\infty_t\rL^2_{\bfx}} \|\bfpsi\|_{\rL^2_t\rL^3_{\bfx}}$.
\autoref{lem:density director tests} then extends \eqref{eq:director molecular limit eps} to every $\bfpsi\in\cX_{\bfd}(\varphi)$, proving \eqref{eq:weak sol d} in its stronger tension form.  Finally, \eqref{eq:phi eps strong} and \eqref{eq:u eps compact} give \eqref{eq:weak sol solid}. Exact rigidity passes to the limit as before.

\subsection{Energy inequality and the quantitative interval}

We pass to the energy inequality along the same subsequence.  By
\eqref{eq:u eps compact} and \eqref{eq:eps global material strong L2}, after a further subsequence,
\[
 \bfu_\eps(t)\to\bfu(t)\quad\text{strongly in }\rL^2(\Omega),
 \qquad
 \overline\bfd_\eps(t)\to\overline\bfd(t)\quad\text{strongly in }\rL^2(\Omega)
\]
for a.e.\ $t\in(0,T)$.  
Fix such a time $t$ at which, in addition, the energy inequality \eqref{eq:energy-main GL} holds for every~$\eps$ in the sequence and $\|\overline\bfd_\eps(t)\|_{\rH^1(\Omega)}\le C_T$. Since the sequence is countable, a.e.\ $t\in(0,T)$ has these properties. Choose a subsequence $\eps_k\to0$, depending on $t$, such that
\[
 \lim_{k\to\infty}\int_\Omega|\nabla\overline\bfd_{\eps_k}(t)|^2\srd\bfx
 =\liminf_{\eps\to0}\int_\Omega|\nabla\overline\bfd_\eps(t)|^2\srd\bfx .
\]
Due to the boundedness of $(\overline\bfd_{\eps_k}(t))_k$ in $\rH^1(\Omega)$, a further subsequence converges weakly in $\rH^1(\Omega)$, and its limit is $\overline\bfd(t)$ thanks to the strong $\rL^2$-convergence. 
By the weak lower semicontinuity of the norm, we thus have
\[
    \int_\Omega|\nabla\overline\bfd(t)|^2\srd\bfx \le\liminf_{\eps\to0}\int_\Omega|\nabla\overline\bfd_\eps(t)|^2\srd\bfx.
\]
Invoking the change of variables $\bfx=\bfh_\eps(t)+\bfO_\eps(t)\bfy$ and \eqref{eq:prepared anchors eps}, we find that
\[
 \int_{\cS_\eps(t)}|\nabla\overline\bfd_\eps(t)|^2\srd\bfx
 =\int_{\cS_0}|\nabla\bfd_0^\eps|^2\srd\bfx
 \rightarrow
 \int_{\cS_0}|\nabla\bfd_0|^2\srd\bfx
 =\int_{\cS(t)}|\nabla\overline\bfd(t)|^2\srd\bfx.
\]
Subtracting this gives, for a.e.\ $t$,
\[
    \int_{\cF(t)}|\nabla\bfd(t)|^2\srd\bfx \le\liminf_{\eps\to0}\int_{\cF_\eps(t)}|\nabla\bfd_\eps(t)|^2\srd\bfx.
\]
The kinetic energy converges at the same a.e.\ time slices.  
Weak lower semicontinuity of the viscous dissipation and
\eqref{eq:tau global liminf}, together with the nonnegative Ginzburg--Landau potential and the convergence of the initial energy  $E_0^\eps\to E_0$, now yields \eqref{eq:energy-main} for a.e.\ $t\in(0,T)$.

The preceding arguments yield a constrained weak solution on $[0,T]$ for every $T<T_*$, obtained along a subsequence that may depend on $T$. 
In order to obtain a single solution on $[0,T_*)$, we use a diagonal argument. 
To this end, for each $\eps$, let $(\bfu_\eps,\bfd_\eps,\varphi_\eps)$ be the Ginzburg--Landau solution from \autoref{thm:weak sol GL} with anchor $\bfd_0^\eps$. 
This solution exists on $[0,T_*^\eps]$, where $T_*^\eps$ is given by \eqref{eq:T star main theorem} with~$E_0^\eps$ in place of $E_0$, and $T_*^\eps\to T_*$ because $E_0^\eps\to E_0$. 
Choose $T_m\uparrow T_*$ and a sequence $\eps_j\to0$. The arguments of this section, applied on $[0,T_1]$, give a subsequence along which the Ginzburg--Landau solutions converge to a constrained weak solution on $[0,T_1]$.
Applied to this subsequence on $[0,T_2]$, they give a further subsequence, and so on.
Along the diagonal subsequence, the Ginzburg--Landau solutions converge on $[0,T_m]$ for every~$m$, and the limits on different intervals coincide on the smaller one, since they are limits of the same sequence. 
Thus, they give rise to a single triple $(\bfu,\bfd,\varphi)$ on~$[0,T_*)$ that is a constrained weak solution on $[0,T]$ for every $T<T_*$.

By \eqref{eq:energy-main}, $(\bfu,\bfd,\varphi)$ satisfies \eqref{eq:kinetic energy dissipation hypothesis} on each such interval, so $\dist(\cS(t),\partial\Omega)\ge\delta_0/2$ for $0\le t<T_*$ by \autoref{cor:uniform separation}, completing the proof of \autoref{prop:constrained local existence}.

The following section is concerned with the continuation of this solution up to the first contact time in \autoref{ssec:maximal constrained} and then after contact in \autoref{ssec:postimpact constrained}.

\section{Continuation of the solution}\label{sec:continuation of sol}

In this section, we discuss the continuation of the solution after the limit passage in the Ginzburg--Landau regularization, first before possible contact with the outer fluid boundary, see \autoref{ssec:maximal constrained}, and then after contact, see \autoref{ssec:postimpact constrained}.

\subsection{Pre-contact continuation of the constrained solution}\label{ssec:maximal constrained}

We iterate the local constrained construction.  The lower bound for the next
collision-free interval depends only on the current separation and energy, not on an $\rH^2$-norm of the
current director.  Thus the loss of regularity in the limit $\eps\to0$ does not prevent restarting.

Restart of the constrained solution requires collision-free endpoint traces and a gluing property.

\begin{lem}\label{lem:constrained endpoint traces}
Let $(\bfu,\bfd,\varphi)$ be a constrained weak solution on $[0,T_0)$, that is, a weak solution on~$[0,T]$ in the sense of \autoref{def:weak sol} for every $T<T_0$, and assume
\[
  \inf_{t_1<t<T_0}\dist(\cS(t),\partial\Omega)>0
\]
for some $t_1<T_0$.  Then the rigid motion has a limit
$\cS_*:=\bfh_*+\bfO_*\cS_0\Subset\Omega$ as $t\to T_0$.  There exist
\[
  \bfu_*\in\cH_{\cS_*}(\Omega),
  \qquad
  \bfd_*\in \rH^1(\Omega\setminus\overline{\cS_*};\bS_+^2),
\]
satisfying the anchoring conditions
\begin{equation}\label{eq:endpoint anchoring}
 \operatorname{Tr}_{\partial\Omega}\bfd_*=\bfd_0,\qquad
 \operatorname{Tr}_{\partial\cS_*}\bfd_*(\bfx)=\bfd_0\bigl(\bfO_*^\top(\bfx-\bfh_*)\bigr),
\end{equation}
that is, \eqref{eq:weak trace anchoring} for the configuration $(\bfh_*,\bfO_*)$, and a sequence chosen from the common full-measure set of admissible times $t_j\to T_0$ such that
\begin{equation*}
  \bfu(t_j)\rightharpoonup\bfu_*\quad\text{in }\rL^2(\Omega),
  \qquad
  \overline{\bfd}(t)\to\overline{\bfd}_*
       \quad\text{strongly in }\rL^2(\Omega)\text{ as }t\to T_0,
\end{equation*}
where $\overline{\bfd}$ denotes the material extension through the solid and
$\overline{\bfd}_*|_{\Omega\setminus\overline{\cS_*}}=\bfd_*$.  Moreover
\begin{equation}\label{eq:constrained endpoint energy lsc}
 \frac12\|\bfu_*\|_{\rL^2(\Omega)}^2
 +\frac12\|\nabla\bfd_*\|_{\rL^2(\Omega\setminus\overline{\cS_*})}^2
 \le \liminf_{j\to\infty}E(t_j).
\end{equation}
For every $\bfzeta\in \rC_c^\infty(\Omega;\R^3)$ which is divergence free and rigid in a fixed positive
neighborhood of $\cS_*$,
\begin{equation}\label{eq:constrained endpoint momentum trace}
 \lim_{t\to T_0}\int_\Omega\bfu(t)\cdot\bfzeta\srd\bfx
 =\int_\Omega\bfu_*\cdot\bfzeta\srd\bfx.
\end{equation}
Finally, suppose that $(\bfu,\bfd,\varphi)$ is continued beyond $T_0$ by a constrained weak solution on $[T_0,T_1]$ with initial state $(\bfu_*,\bfd_*,\cS_*)$ and initial energy $E_*:=\frac12\|\bfu_*\|_{\rL^2(\Omega)}^2+\frac12\|\nabla\bfd_*\|_{\rL^2(\Omega\setminus\overline{\cS_*})}^2$. 

Then the concatenation is a constrained weak solution on $[0,T_1]$ that satisfies the energy inequality \eqref{eq:energy-main} with the original initial energy $E_0$.
\end{lem}

\begin{proof}
The energy identity for the rigid part and \autoref{lem:rigid velocity energy control} give uniform bounds for $\bfh'$ and $\bfomega$, so $\bfh(t)\to\bfh_*$ and $\bfO(t)\to\bfO_*$ as $t\to T_0$.
Next, we define the extension of $\bfd$ on the whole container by
\begin{equation*}
    \overline{\bfd}(t,\bfx):=
    \left\{
    \begin{aligned}
        &\bfd(t,\bfx),&&\tfor \bfx\in\cF(t),\\
        &\bfd_0\bigl(\bfO(t)^\top(\bfx-\bfh(t))\bigr),&&\tfor \bfx\in\cS(t).
    \end{aligned}
    \right.
\end{equation*}
Moreover, we denote by $\widetilde{\bftau}$ the zero extension of $\bftau\in \rL^2(\cQ_{\cF})$.  
From \autoref{lem:material gluing}, it follows that
\[
    \partial_t\overline{\bfd}+\bfu\cdot\nabla\overline{\bfd}=\widetilde{\bftau} \tin \cD'((t_1,T_0)\times\Omega).
\]
Now we set $\bfq:=\overline\bfd-\bfd_0$. 
The energy inequality \eqref{eq:energy-main} and the fact that the solid part of $\overline\bfd$ is a rigid image of~$\bfd_0|_{\cS_0}$ yield that the extension $\overline\bfd$ is bounded in $\rL^\infty(t_1,T_0;\rH^1(\Omega))$. 
Its trace on $\partial\Omega$ equals that of $\bfd_0$ by the outer anchoring condition in \eqref{eq:weak trace anchoring}. 
Furthermore, $\partial_t\overline\bfd\in\rL^2(t_1,T_0;\rH^{-1}(\Omega))$ by the last assertion of \autoref{lem:material gluing}.
Since $\bfd_0$ does not depend on time, we obtain
\[
 \bfq\in\rL^2(t_1,T_0;\rH_0^1(\Omega)),\qquad
 \partial_t\bfq\in\rL^2(t_1,T_0;\rH^{-1}(\Omega)).
\]
By the Lions--Magenes lemma for the triple $\rH_0^1(\Omega)\hookrightarrow\rL^2(\Omega)\hookrightarrow\rH^{-1}(\Omega)$, the function $\bfq$, and hence~$\overline\bfd$, belongs to $\rC([t_1,T_0];\rL^2(\Omega))$. In particular, $\overline\bfd(t)\to\overline\bfd_*$ strongly in $\rL^2(\Omega)$ as $t\to T_0$. 
By the boundedness of $\overline\bfd(t)$ in $\rH^1(\Omega)$, this convergence carries over to weak convergence in $\rH^1(\Omega)$, and $\overline\bfd_*\in\rH^1(\Omega)$.
Thanks to $\bfh(t)\to\bfh_*$ and $\bfO(t)\to\bfO_*$, the part of $\overline\bfd(t)$ defined on $\cS(t)$ converges as well, and a.e.\ in $\cS_*$, we get $\overline{\bfd}_*(\bfx)=\bfd_0\bigl(\bfO_*^\top(\bfx-\bfh_*)\bigr)$.
This reveals that $\bfd_*:=\overline\bfd_*|_{\Omega\setminus\overline{\cS_*}}$ satisfies \eqref{eq:endpoint anchoring}. 
The corresponding anchoring condition on the rigid body interface holds because $\overline\bfd_*\in\rH^1(\Omega)$ has the same trace on $\partial\cS_*$ from both sides, and the anchoring condition on the outer boundary holds because $\overline\bfd_*-\bfd_0\in\rH_0^1(\Omega)$.
The strong $\rL^2$-convergence allows us to infer that $|\bfd_*|=1$ and $(\bfd_*)_3\ge0$. 
Finally, the gradient norm of the part of $\overline{\bfd}$ defined on the solid domain $\cS(t)$ does not depend on time, i.e.,
\[
    \int_{\cS(t)}|\nabla\overline\bfd(t)|^2\srd\bfx=\int_{\cS_0}|\nabla\bfd_0|^2\srd\bfx=\int_{\cS_*}|\nabla\overline\bfd_*|^2\srd\bfx.
\]
Weak lower semicontinuity of the $\rH^1(\Omega)$-norm thus gives
\[
    \|\nabla\bfd_*\|_{\rL^2(\Omega\setminus\overline{\cS_*})}^2\le\liminf_{t\to T_0}\|\nabla\bfd(t)\|_{\rL^2(\cF(t))}^2 ,
\]
which is the director part of the energy inequality \eqref{eq:constrained endpoint energy lsc}.

Let $E(t)$ be the energy on the left-hand side of \eqref{eq:energy-main}. 
We then choose times $t_j\to T_0$, in the set with full measure on which the energy inequality holds and $\bfu(t_j)$ is rigid on $\cS(t_j)$, such that~$E(t_j)\to\operatorname*{ess\,liminf}_{t\to T_0}E(t)$.
Upon extracting, the boundedness of the sequence $\bfu(t_j)$ in $\cH(\Omega)$ allows us to derive its weak convergence in
$\rL^2(\Omega)$ to some $\bfu_*\in\cH(\Omega)$. 
On the set $\cS(t_j)$, we write
\[
  \bfu(t_j,\bfx)=\bfell_j+\bfomega_j\times(\bfx-\bfh(t_j)).
\]
The rigid kinetic-energy identity bounds the rigid velocities $(\bfell_j,\bfomega_j)$. 
A further extraction yields their convergence to~$(\bfell_*,\bfomega_*)$.  
Every compact subset of $\cS_*$ lies in $\cS(t_j)$ for all large $j$, so
weak convergence against tests supported in that compact set gives
\[
    \bfu_*(\bfx)=\bfell_*+\bfomega_*\times(\bfx-\bfh_*)
  \quad\text{for a.e.\ }\bfx\in\cS_*.
\]
This establishes $\bfu_*\in\cH_{\cS_*}(\Omega)$, and weak lower semicontinuity gives the kinetic part of
\eqref{eq:constrained endpoint energy lsc}.

Finally, fix a divergence-free $\bfzeta\in\rC_c^\infty(\Omega;\R^3)$ that is rigid on $(\cS_*)_r$ for some $r>0$. 
By the convergence $\cS(t)\to\cS_*$, there is $t_2<T_0$ such that $\bfzeta$ is rigid on $\cS(t)_{r/2}$ for all $t\in(t_2,T_0)$.
Multiplied by a cut-off function in time whose support is in $(t_2,T_0)$, it is thus an admissible test field in the class \eqref{eq:collar momentum tests}. 
The momentum equation shows that $t\mapsto(\bfu(t),\bfzeta)_{\rL^2(\Omega)}$ has a distributional derivative in $\rL^1(t_2,T_0)$, because the convective and Ericksen stresses are bounded in $\rL^\infty_t\rL^1_{\bfx}$, while the viscous contribution is in $\rL^2(t_2,T_0)$.
Hence the pairing has a unique one-sided limit.  
Evaluating it along $t_j$ proves \eqref{eq:constrained endpoint momentum trace}.

It remains to prove the gluing statement. 
For this purpose, let $\bfzeta$ be a test field on $[0,T_1)$ in the class $\cT_{\bfu}^{\rm col}$ from~\eqref{eq:collar momentum tests}. 
Testing the first piece with $\bfzeta$ multiplied by a cut-off function in time gives, for a.e.\ $s<T_0$, its momentum identity on $[0,s]$ with the additional term $\int_\Omega\bfu(s)\cdot\bfzeta(s)\srd\bfx$ on the left-hand side. 
By the uniform convergence $\bfzeta(s)\to\bfzeta(T_0)$ uniformly, this term also converges to $\int_\Omega\bfu_*\cdot\bfzeta(T_0)\srd\bfx$ as $s\to T_0$ by \eqref{eq:constrained endpoint momentum trace}, and this limit is the initial term in the momentum identity of the second piece. 
Adding the two identities gives the momentum identity on $[0,T_1]$ for tests $\bfzeta \in \cT_{\bfu}^{\rm col}$, and \autoref{prop:rigid test equivalence} extends it to all $\bfzeta\in\cT_{\bfu}$.
In the same way, the endpoint terms in the director equation cancel because~$\overline\bfd(s)\to\overline\bfd_*$ strongly in $\rL^2(\Omega)$, and the transport equation for $\varphi$ holds across~$T_0$ because the rigid motion is continuous.

For the energy inequality, let $\mathcal D(s,t)$ denote the two dissipation integrals in \eqref{eq:energy-main} taken over the time interval $(s,t)$. For a.e.\ $t\in(T_0,T_1)$, the energy inequality of the second piece, \eqref{eq:constrained endpoint energy lsc}, and the energy inequality of the first piece at the times $t_j$ give
\[
 E(t)+\mathcal D(0,t)\le E_*+\mathcal D(0,T_0)
 \le\liminf_{j\to\infty}\bigl(E(t_j)+\mathcal D(0,t_j)\bigr)\le E_0 .\qedhere
\]
\end{proof}

The constrained system can then be restarted from data in the energy class.

\begin{lem}\label{lem:constrained energy restart}
Let $\cS_* =\bfh_*+\bfO_*\cS_0$ satisfy $\rho_*:=\dist(\cS_*,\partial\Omega)>0$, consider $\bfu_*\in\cH_{\cS_*}(\Omega)$, and assume that $\bfd_*\in \rH^1(\Omega\setminus\overline{\cS_*};\bS_+^2)$ satisfies \eqref{eq:endpoint anchoring}. 
Define $E_*:=\frac12\|\bfu_*\|_{\rL^2(\Omega)}^2+
      \frac12\|\nabla \bfd_*\|_{\rL^2(\Omega\setminus\overline{\cS_*})}^2$.

Then there is a constrained weak solution starting from $(\bfu_*,\bfd_*,\cS_*)$ on a time interval of length
\begin{equation*}
    \theta_*:=\left\{
    \begin{aligned}
        &\dfrac{\rho_*}{2C_{\cS_0}\sqrt{2E_*}},&&\tfor E_*>0,\\
        &+\infty,&&\tfor E_*=0.
    \end{aligned}
    \right.
\end{equation*}
Moreover, the body remains at distance at least $\rho_*/2$ from $\partial\Omega$ on that interval. 
\end{lem}

\begin{proof}
If $E_*=0$, then $\bfu_*=0$ and $\nabla\bfd_*=0$.
This means that the body at rest at $\cS_*$, together with $\bfu=0$ and~$\bfd(t)=\bfd_*$, is a constrained weak solution for all times.

For a global definition, we extend $\bfd_*$ through the solid by the associated material datum, i.e.,
\begin{equation*}
    \overline\bfd_*(\bfx):= \left\{
    \begin{aligned}
        &\bfd_*(\bfx),&&\tfor \bfx\in\Omega\setminus\overline{\cS_*},\\
        &\bfd_0\bigl(\bfO_*^\top(\bfx-\bfh_*)\bigr),&&\tfor \bfx\in\cS_*.
    \end{aligned}
    \right.
\end{equation*}
From \eqref{eq:endpoint anchoring}, it follows that $\overline\bfd_*\in \rH^1(\Omega;\bS_+^2)$.  
This allows us to apply \autoref{lem:hemisphere preparation}, where the approximation index tied to the Ginzburg--Landau parameter $\eps$, and we choose
\begin{equation}\label{eq:restart prepared family}
 \overline\bfd_*^\eps\in \rH^2(\Omega;\bS_+^2),
 \qquad |\overline\bfd_*^\eps|=1,
 \qquad
 \overline\bfd_*^\eps\to\overline\bfd_*
 \quad\text{strongly in }\rH^1(\Omega)
\end{equation}
as $\eps\to0$.
Next, we apply the construction of \autoref{sec:GL-approx} and \autoref{sec:eps to 0} with the restart time as new initial time, with $\cS_*$, $\bfu_*$, and $\overline\bfd_*^\eps$ replacing $\cS_0$, $\bfu_0$, and $\bfd_0^\eps$.
In particular, $\varphi_0=\dsone_{\cS_*}$, and no change of coordinates is needed.
We keep the description of the rigid motion through the original reference body,~$\cS(t)=\bfh(t)+\bfO(t)\cS_0$ with $\bfh(0)=\bfh_*$ and $\bfO(0)=\bfO_*$, so that the points of $\cS_*$ move according to
\[
    \bfZ_*(t,\bfx)=\bfh(t)+\bfO(t)\bfO_*^\top(\bfx-\bfh_*),\qquad \bfx\in\cS_*.
\]
By definition, $\overline\bfd_*(\bfx)=\bfd_0(\bfO_*^\top(\bfx-\bfh_*))$ in $\cS_*$, so the anchoring datum transported by $\bfZ_*$ is $\bfd_0(\bfO(t)^\top(\bfx-\bfh(t)))$ on $\cS(t)$, as in \eqref{eq:weak trace anchoring}.
Let us observe that the solid mass and the eigenvalues of the inertia tensor are invariant under rigid placement, the Korn--Poincar\'e constant depends only on~$\Omega$, and the lifespan estimate depends on the initial gap to $\partial\Omega$.
As a result, the construction and the uniform separation estimates of \autoref{ssec:uniform precontact} and \autoref{sec:GL-approx} apply with $\cS_*$ as initial solid.

For each sufficiently small $\eps$, \autoref{thm:weak sol GL} thus gives a solution to the moving-domain problem with Ginzburg--Landau regularization starting from $(\bfu_*,\overline\bfd_*^\eps,\cS_*)$, and with initial energy
\[
    E_*^\eps=\frac12\|\bfu_*\|_{\rL^2(\Omega)}^2+\frac12\|\nabla\overline\bfd_*^\eps\|_{\rL^2(\Omega\setminus\overline{\cS_*})}^2 \rightarrow E_*.
\]
Here the initial potential energy vanishes due to the sphere-valuedness. 
By \autoref{thm:weak sol GL}, this solution then exists at least on $[0,\theta_\eps]$, where, in accordance with \eqref{eq:T star main theorem} with $\rho_*$ and $E_*^\eps$ in place of $\delta_0$ and $E_0$,
\[
    \theta_\eps:=\frac{\rho_*}{2C_{\cS_0}\sqrt{2E_*^\eps}}, \qquad \theta_\eps\rightarrow\theta_*.
\]
In particular, by \autoref{cor:uniform separation}, for every $\theta<\theta_*$ and all sufficiently small $\eps$,
\begin{equation}\label{eq:restart uniform gap}
 \dist(\cS_\eps(t),\partial\Omega)\ge\frac{\rho_*}{2}
 \qquad(0\le t\le\theta).
\end{equation}

It remains to pass to the limit with only $\rH^1$-regularity of $\bfd_*$.  
The rigid motions are compact in $\rC([0,\theta])$ due to the uniform energy inequality and~\eqref{eq:restart uniform gap}, the velocities are strongly compact in~${\rL^2((0,\theta)\times\Omega)}$ by \autoref{lem:FSI compactness criterion}, and the directors are strongly compact in $\rL^2_{\rm loc}(0,\theta;\rH^1_{\rm loc}(\cF(t)))$ by \autoref{lem:LW spacetime Cauchy}.
We recall that \eqref{eq:id of tensor field} from \autoref{prop:eps local strong d} identifies the limiting molecular field as~$\Delta\bfd+|\nabla\bfd|^2\bfd$. 
The argument of \autoref{lem:anchoring eps limit}, applied with $\overline\bfd_*^\eps$ from~\eqref{eq:restart prepared family} instead of~\eqref{eq:prepared anchors eps}, yields that the limiting director satisfies the anchoring condition \eqref{eq:weak trace anchoring} with the original datum $\bfd_0$.
The elastic stress first passes against tests $\bfzeta \in \cT_{\bfu}^{\rm col}$ and then against every standard rigid test $\bfzeta\in\cT_{\bfu}$ thanks to \autoref{prop:rigid test equivalence}.  

With regard to the energy inequality, we argue as at the end of \autoref{sec:eps to 0}, with the material extensions of the restart director fields.
The only change is that the solid part of the gradient norm now equals~$\|\nabla\overline\bfd_*^\eps\|_{\rL^2(\cS_*)}^2$, which converges to $\|\nabla\overline\bfd_*\|_{\rL^2(\cS_*)}^2$ by \eqref{eq:restart prepared family}. 
The convergence $E_*^\eps\to E_*$ shown above implies that the limiting energy inequality has initial energy $E_*$, so in total, the limit is a constrained weak solution with initial state $(\bfu_*,\bfd_*,\cS_*)$.

The preceding arguments give a constrained weak solution on $[0,\theta]$ for every $\theta<\theta_*$.
Repeating the diagonal argument at the end of \autoref{sec:eps to 0}, with a sequence $\theta_m\uparrow\theta_*$ instead of $T_m\uparrow T_*$, we obtain a single constrained weak solution on $[0,\theta_*)$.
Moreover, \autoref{cor:uniform separation} gives $\dist(\cS(t),\partial\Omega)\ge\rho_*/2$ for~$0\le t<\theta_*$. 
If $\theta_*<\infty$, \autoref{lem:constrained endpoint traces}, applied with $T_0=\theta_*$, provides the endpoint state at $\theta_*$.
\end{proof}

\begin{proof}[Proof of \autoref{exist:main}\,(i) and \autoref{cor:quantitative separation main}]
If $E_0=0$, then $\bfu_0=0$ and $\nabla\bfd_0=0$ a.e.\ in $\cF_0$.
Thus, the body at rest together with $\bfu=0$ and $\bfd(t)=\bfd_0|_{\cF_0}$ is a global constrained solution, and $T_{\rm col}=\infty$.
Therefore, we assume $E_0>0$ and write $C_R:=C_{\cS_0}$ in the remainder of the proof.

We start by constructing one branch of the solution.
\autoref{prop:constrained local existence} yields that there is a constrained weak solution on $[0,T_*)$.
We then fix $t_1\in(0,T_*)$ and restrict the solution to $[0,t_1]$.
By \autoref{lem:constrained endpoint traces}, applied with $T_0=t_1$, it has an endpoint state $(\bfu_*,\bfd_*,\cS_*)$ at $t_1$ that satisfies the hypotheses of \autoref{lem:constrained energy restart}.
We then extend this branch inductively without modifying pieces that have been constructed previously. 
Suppose that it has been constructed on $[0,t_j]$, with gap $\rho_j:=\dist(\cS(t_j),\partial\Omega)>0$ and endpoint energy at most $E_0$. By \autoref{lem:constrained energy restart}, the endpoint data at $t_j$ can be continued to a constrained weak solution for at least $\rho_j/(2C_R\sqrt{2E_0})$ units of time, during which the gap stays at least $\rho_j/2$. 
We retain the piece of length
\[
 t_{j+1}-t_j:=\frac{\rho_j}{4C_R\sqrt{2E_0}},
\]
whose terminal gap is positive, apply \autoref{lem:constrained endpoint traces} at $t_{j+1}$, and continue.
If the energy of an endpoint state vanishes, we continue by the solution at rest, as in the case $E_0=0$.
The last assertion of \autoref{lem:constrained endpoint traces} implies that every such concatenation is a constrained weak solution on the combined interval and satisfies the energy inequality \eqref{eq:energy-main} with the original initial energy $E_0$.

We define $T_{\rm col}:=\lim_{j\to\infty}t_j\in(0,\infty]$. 
The energy inequality \eqref{eq:energy-main} yields $\|\bfu(t)\|_{\rL^2(\Omega)}\le\sqrt{2E_0}$ for a.e.\ $t<T_{\rm col}$. 
By \autoref{lem:rigid velocity energy control}, the velocity of every point of the body is bounded by $C_{\cS_0}\sqrt{2E_0}$, showing that the rigid trajectory and the gap $t\mapsto\dist(\cS(t),\partial\Omega)$ are Lipschitz continuous on $[0,T_{\rm col})$.
If $T_{\rm col}<\infty$, the gap thus has a limit $\rho\ge0$ as $t\to T_{\rm col}$. 
The case $\rho>0$ implies $\rho_j\ge\rho/2$ for all large $j$, hence $t_{j+1}-t_j\ge\rho/(8C_R\sqrt{2E_0})$, which contradicts the convergence of $(t_j)$. 
Therefore, $\rho=0$, the rigid motion extends continuously to $T_{\rm col}$, and \eqref{eq:collision alternative main} holds. 
By positivity of every $\rho_j$, the time $T_{\rm col}$ is the first contact time of the constructed branch, proving \autoref{exist:main}\,(i) and the first assertion of \autoref{exist:main}\,(ii).

Finally, we establish a lower bound for the collision time $T_{\rm col}$.
In fact, the branch satisfies~\eqref{eq:energy-main} on every compact subinterval of $[0,T_{\rm col})$, so \autoref{lem:dissipative rigid displacement} applies, and \eqref{eq:coupled separation} is \eqref{eq:energy dissipation rigid displacement}. 
Recall from~\eqref{eq:T ed main theorem} that $T_{\rm ed}$ is the unique time at which the right-hand side of \eqref{eq:energy dissipation rigid displacement} vanishes, and that this right-hand side is positive for $t<T_{\rm ed}$. Hence, if $T_{\rm col}<T_{\rm ed}$, the gap would stay bounded away from zero as $t\to T_{\rm col}$, which contradicts $\rho=0$. Therefore $T_{\rm col}\ge T_{\rm ed}$.
This completes the proof of \autoref{cor:quantitative separation main}.
\end{proof}

\subsection{Post-impact constrained continuation}\label{ssec:postimpact constrained}

We identify the one-sided data at the first wall contact and then solve the constrained Ericksen--Leslie system on the resulting fixed fluid region. The post-impact solution is obtained by a Ginzburg--Landau limit on this fixed domain.

Let $G\subset\R^3$ be bounded and open, and let $\cD_G$, $\cH_G$, and $\cV_G$ be defined by \eqref{eq:arbitrary domain solenoidal spaces}.
We recall that the fixed-domain solution concept was stated in \autoref{def:postcontact constrained weak solution}.

\begin{prop}
\label{prop:postcontact constrained arbitrary domain}
For every bounded open $G\subset\R^3$, every $\bfu_*\in\cH_G$, and every
$\bfd_*\in \rH^1(G;\bS_+^2)$, there exists a post-impact constrained weak solution in the sense of
\autoref{def:postcontact constrained weak solution} on $[0,\infty)$.  Moreover, there are
$\eps_j\to0$ and corresponding global fixed-domain Ginzburg--Landau solutions
$(\bfu_j,\bfd_j,\bfM_j)$ with the same initial data such that, for every $T<\infty$,
\begin{align*}
 &\bfu_j\to\bfu &&\text{strongly in }\rL^2((0,T)\times G),
 \\
 &\bfd_j\to\bfd &&\text{strongly in }\rL^2((0,T)\times G),
 \quad \bfd_j\stackrel{*}{\rightharpoonup}\bfd\text{ in }\rL^\infty(0,T;\rH^1(G)),
 \\
 &\bfd_j\to\bfd &&\text{strongly in }\rL^2(0,T;\rH^1(V))
 \quad\text{for every }V\Subset G,
 \\
 &\bfM_j\rightharpoonup\bftau &&\text{weakly in }\rL^2((0,T)\times G).
\end{align*}
\end{prop}

\begin{proof}
We fix $T>0$, let $\eps_j\to0$, and apply
\autoref{prop:postcontact GL arbitrary domain} on $G$ with the same data $(\bfu_*,\bfd_*)$ for every~$j$.
The director condition $|\bfd_*|=1$ implies that the initial Ginzburg--Landau potential vanishes. 
From \eqref{eq:postcontact energy GL}, it then follows, uniformly in $j$,
\begin{equation}\label{eq:fixed domain eps uniform energy}
    \begin{aligned}
        \|\bfu_j\|_{\rL^\infty_t\rL^2_{\bfx}}+\|\bfu_j\|_{\rL^2_t\rH^1_{0,\bfx}}+\|\bfd_j\|_{\rL^\infty_t\rH^1_{\bfx}}+\|\bfM_j\|_{\rL^2_{t,\bfx}}
        &\le C_T,\\
        \sup_{0<t<T}\frac1{\eps_j^2}\int_G(1-|\bfd_j(t)|^2)^2\srd\bfx
        &\le C_T, \qquad |\bfd_j|\le1,\quad (\bfd_j)_3\ge0,
    \end{aligned}
\end{equation}
where we recall $\bfM_j=\Delta \bfd_j+\eps_j^{-2}(1-|\bfd_j|^2)\bfd_j$.
By $\bfu_j\in \rL^2_t\rL^6_{\bfx}$ and $\nabla \bfd_j\in \rL^\infty_t\rL^2_{\bfx}$, the director equation implies that $\partial_t \bfd_j=\bfM_j-\bfu_j\cdot\nabla \bfd_j$ is bounded in $\rL^2(0,T;\rH^{-1}(G))$.
For $\bfq_j:=\bfd_j-\bfd_*$, we find that $\bfq_j$ is bounded in $\rL^\infty(0,T;\rH_0^1(G))$, while
$\partial_t \bfq_j=\partial_t \bfd_j$ is bounded in $\rL^2(0,T;\rH^{-1}(G))$. 
The zero-extension map is continuous from $\rH_0^1(G)$ into $\rH_0^1(B)$ for any fixed ball $B\Supset\overline G$, so the Rellich--Kondrachov theorem on $B$ shows that $ \rH_0^1(G)\Subset \rL^2(G)$ for every bounded open $G$. 
Applying the Aubin--Lions lemma to $\bfq_j$ imply that $\bfd_j\to \bfd$ strongly in $\rL^2((0,T)\times G)$ and, after extraction, weak-* in $\rL^\infty_t\rH^1_{\bfx}$. 
The potential bound in
\eqref{eq:fixed domain eps uniform energy} therefore implies $|\bfd|=1$ a.e. 
Since $(\bfd_j)_3\ge0$ by \eqref{eq:fixed domain eps uniform energy}, the strong convergence gives $d_3\ge0$ a.e. 
The fact that $\bfd_j-\bfd_*\in\rH_0^1(G;\R^3)$ and the weak closedness of this space imply $\bfd-\bfd_*\in\rL^\infty(0,T;\rH_0^1(G;\R^3))$.

Let now $V\Subset U\Subset G$.
The sequence $(\bfd_j)$ is Cauchy in $\rL^2((0,T)\times U)$ by the preceding Aubin--Lions argument.
By \eqref{eq:fixed domain eps uniform energy}, the maps $\bfd_j(t)$ satisfy \eqref{eq:LW slice equation} on $U$ with $e_j=\eps_j$ and $\bftau_j=\bfM_j(t)$, as well as~\eqref{eq:LW slice bounds}. 
Thus, for every $V\Subset G$, \autoref{lem:LW spacetime Cauchy} yields
\begin{equation}\label{eq:postcontact local H1 constrained}
 \bfd_j\to \bfd\quad\text{strongly in }\rL^2(0,T;\rH^1(V)).
\end{equation}

After passing to the limit $\bfM_j\rightharpoonup\bftau$ in $\rL^2((0,T)\times G)$, arguing as in \eqref{eq:id of tensor field} in the proof of \autoref{prop:eps local strong d} gives $\bftau=\Delta \bfd+|\nabla \bfd|^2\bfd$ in $\cD'((0,T)\times G)$.

For the velocity, we fix the basis $(\bfv_m)$ from \autoref{lem:fixed domain velocity compactness}.  
Testing the weak momentum equation with~$\chi(t)\bfv_m$, $\chi\in \rC_c^1(0,T)$, shows that $t\mapsto(\bfu_j(t),\bfv_m)_{\rL^2(G)}$ belongs to $\rW^{1,1}(0,T)$, with distributional derivative given by the same expression as in the proof of \autoref{prop:postcontact GL arbitrary domain}. 
By \eqref{eq:fixed domain eps uniform energy}, we get
\[
    \int_0^T\left|\frac{\rd}{\rd t}(\bfu_j,\bfv_m)_{\rL^2(G)}\right|\srd t \le C_m\left(T\|\bfu_j\|_{\rL^\infty_t\rL^2_{\bfx}}^2+T^{1/2}\|\bD(\bfu_j)\|_{\rL^2_{t,\bfx}}+T\|\nabla\bfd_j\|_{\rL^\infty_t\rL^2_{\bfx}}^2\right)\le C_{m,T}.
\]
These functions thus have uniformly bounded variation on $(0,T)$, so \autoref{lem:fixed domain velocity compactness} implies
\begin{equation}\label{eq:postcontact eps strong velocity constrained}
    \bfu_j\to \bfu\qquad\text{strongly in }\rL^2((0,T)\times G).
\end{equation}
For test functions with compact support, \eqref{eq:postcontact local H1 constrained} and \eqref{eq:postcontact eps strong velocity constrained} allow for the limit passage in the two nonlinear terms in the momentum and director equations. 
The limiting director equation, written with the tension $\bftau$, depends continuously on $\bfpsi\in\cX_G$.
The transport term is bounded by~$\|\bfu\|_{\rL^2_t\rL^6_{\bfx}}\|\nabla\bfd\|_{\rL^\infty_t\rL^2_{\bfx}}\|\bfpsi\|_{\rL^2_t\rL^3_{\bfx}}$, and the tension term by~$\|\bftau\|_{\rL^2_{t,\bfx}}\|\bfpsi\|_{\rL^2_{t,\bfx}}$, where $\bftau\in\rL^2((0,T)\times G)$ as the weak limit of the fields $\bfM_j$, which are bounded by~\eqref{eq:fixed domain eps uniform energy}.
The density of $\rC_c^\infty([0,T)\times G)$ in~$\cX_G$ allows us to infer that the director equation holds for every~$\bfpsi\in\cX_G$.

For the energy inequality, we pass to a further subsequence such that
\[
 \bfu_j(t)\to\bfu(t),\qquad \bfd_j(t)\to\bfd(t)
 \quad\text{strongly in }\rL^2(G)
\]
for a.e.\ $t\in(0,T)$. Fix such a time $t$ at which \eqref{eq:postcontact energy GL} holds for every $j$. Since $(\bfd_j(t))_j$ is bounded in~$\rH^1(G)$ by \eqref{eq:fixed domain eps uniform energy}, every weak $\rH^1$-limit of a subsequence equals $\bfd(t)$.
Hence
\[
 \|\nabla\bfd(t)\|_{\rL^2(G)}^2
 \le\liminf_{j\to\infty}\|\nabla\bfd_j(t)\|_{\rL^2(G)}^2,
 \qquad
 \|\bfu_j(t)\|_{\rL^2(G)}^2\to\|\bfu(t)\|_{\rL^2(G)}^2.
\]
The Ginzburg--Landau potential $F_{\eps_j}$ is nonnegative, and the dissipation terms $2\int_0^t\int_G|\bD(\bfu)|^2\srd\bfx\srd s$ and $\int_0^t\int_G|\bfM|^2\srd\bfx\srd s$ are weakly lower semicontinuous, since $\bD(\bfu_j)\rightharpoonup\bD(\bfu)$ and $\bfM_j\rightharpoonup\bftau$ weakly in~$\rL^2((0,T)\times G)$.
Passing to the liminf in \eqref{eq:postcontact energy GL}, whose right-hand side is independent of $j$ thanks to $|\bfd_*|=1$, thus yields \eqref{eq:postcontact constrained energy} for a.e.\ $t$.

Finally, note that the trajectories provided by \autoref{prop:postcontact GL arbitrary domain} are defined on $[0,\infty)$. 
A diagonal argument over $T=1,2,\dots$, as at the end of \autoref{sec:eps to 0}, then yields a single solution on $[0,\infty)$. 
This completes the proof of this proposition and especially establishes \autoref{thm:arbitrary open constrained}.
\end{proof}

\begin{rem}\label{rem:arbitrary open GL limit}
The proof shows that the solution of \autoref{thm:arbitrary open constrained} is obtained as $\eps$-limit with $\bfu_\eps\to\bfu$ strongly in $\rL^2_{\rm loc}([0,\infty);\rL^2(G))$, $\bfd_\eps\to\bfd$ strongly in $\rL^2_{\rm loc}([0,\infty);\rL^2(G))$ and in $\rL^2(0,T;\rH^1(V))$ for every $T<\infty$ and every $V\Subset G$, while, for every $T<\infty$, the molecular fields $\bfM_\eps:=\Delta\bfd_\eps+\eps^{-2}(1-|\bfd_\eps|^2)\bfd_\eps$ satisfy $\bfM_\eps \rightharpoonup\bftau$ in $\rL^2((0,T)\times G)$.
\end{rem}

Next, we identify the state of the solution at the first contact time $T_c$. We write $\bfu_c^-$ for a pre-impact velocity trace, that is, a weak $\rL^2(\Omega)$-limit of $\bfu(t_j)$ along a sequence $t_j\to T_c$, which need not be unique, and $\bfu_c^+$ for the post-impact velocity on the residual fluid domain, which is uniquely determined.

\begin{lem}
\label{lem:constrained collision traces}
Let $(\bfu,\bfd,\varphi)$ be a constrained weak solution on $[0,T_c)$ in the sense of \autoref{def:weak sol}, with energy $E(t)$ as in \eqref{eq:energy-main}, such that
\[
 \dist(\cS(t),\partial\Omega)>0\quad(t<T_c),\qquad
 \dist(\cS(t),\partial\Omega)\to0\quad(t\to T_c).
\]
Then the following assertions hold.
\begin{enumerate}[(i)]
\item The rigid motion has a unique limit $\cS_c=\bfh_c+\bfO_c\cS_0$ as $t\to T_c$, and there exists $\bfd_c\in \rH^1(\Omega;\bS_+^2)$ such that the material extension $\overline{\bfd}$ of $\bfd$ through the solid satisfies
\begin{equation}\label{eq:constrained collision director trace}
 \overline{\bfd}(t)\to \bfd_c\quad\text{strongly in }\rL^2(\Omega)
 \text{ and weakly in }\rH^1(\Omega)\quad(t\to T_c).
\end{equation}
\item With $\Omega_c:=\Omega\setminus\overline{\cS_c}$, there exists a unique $\bfu_c^+\in\cH_{\Omega_c}$ such that, for every $\bfzeta\in\cD_{\Omega_c}$ with zero extension $\widetilde{\bfzeta}$ to $\Omega$,
\begin{equation*}
 (\bfu_c^+,\bfzeta)_{\rL^2(\Omega_c)}
 =\lim_{t\to T_c}(\bfu(t),\widetilde{\bfzeta})_{\rL^2(\Omega)}.
\end{equation*}
\item If $\bfu(t_j)\rightharpoonup\bfu_c^-$ in $\rL^2(\Omega)$ along a sequence $t_j\to T_c$ of times at which $\bfu(t_j)$ is rigid on $\cS(t_j)$, then $\bfu_c^-\in\cH_{\cS_c}(\Omega)$ and
$\bfu_c^+=P_{\cH_{\Omega_c}}(\bfu_c^-|_{\Omega_c})$.
\item We have
\begin{equation}\label{eq:constrained collision energy drop}
 \frac12\|\bfu_c^+\|_{\rL^2(\Omega_c)}^2
 +\frac12\|\nabla \bfd_c\|_{\rL^2(\Omega_c)}^2
 \le\operatorname*{ess\,liminf}_{t\to T_c}E(t).
\end{equation}
Moreover, there is a sequence $t_j\to T_c$ as in (iii), chosen from the full-measure set on which the energy inequality holds, such that
\[
 E(t_j)\to\operatorname*{ess\,liminf}_{t\to T_c}E(t),
 \qquad \bfu(t_j)\rightharpoonup\bfu_c^- \text{ in }\rL^2(\Omega),
\]
and every weak limit $\bfu_c^-$ obtained along such a sequence satisfies
\begin{equation}\label{eq:constrained collision energy drop strengthened}
\begin{aligned}
 &\frac12\|\bfu_c^+\|_{\rL^2(\Omega_c)}^2
 +\frac12\|\nabla \bfd_c\|_{\rL^2(\Omega_c)}^2\\
 &\quad+\frac12\|\bfu_c^-\|_{\rL^2(\cS_c)}^2
 +\frac12\|\bfu_c^-|_{\Omega_c}-\bfu_c^+\|_{\rL^2(\Omega_c)}^2
 \le\operatorname*{ess\,liminf}_{t\to T_c}E(t).
\end{aligned}
\end{equation}
\end{enumerate}
\end{lem}

\begin{proof}
By the energy inequality \eqref{eq:energy-main} and \autoref{lem:rigid velocity energy control}, the functions $\bfh$ and $\bfO$ are Lipschitz on~$[0,T_c)$.
Hence $\bfh(t)\to\bfh_c$ and $\bfO(t)\to\bfO_c$ as $t\to T_c$, and we set $\cS_c:=\bfh_c+\bfO_c\cS_0$.
Let $\overline\bfd$ be the material extension \eqref{eq:material extension notation main}, and let $\widetilde\bftau$ be the extension of the tension $\bftau=\Delta \bfd+|\nabla \bfd|^2\bfd$ by zero to the solid.
By \autoref{lem:material gluing}, we get $\partial_t\overline{\bfd}+\bfu\cdot\nabla\overline{\bfd}=\widetilde{\bftau}$ in $\cD'((0,T_c)\times\Omega)$.
The same argument as in the proof of \autoref{lem:constrained endpoint traces} shows that $\bfq:=\overline\bfd-\bfd_0$ satisfies $\bfq\in\rL^2(0,T_c;\rH_0^1(\Omega))$ and $\partial_t\bfq\in\rL^2(0,T_c;\rH^{-1}(\Omega))$.
The Lions--Magenes lemma then yields $\overline\bfd\in\rC([0,T_c];\rL^2(\Omega))$, so $\overline\bfd(t)$ converges strongly in $\rL^2(\Omega)$ as $t\to T_c$ and, being bounded in $\rH^1(\Omega)$, also weakly in $\rH^1(\Omega)$.
This proves \eqref{eq:constrained collision director trace}. 
The limit $\bfd_c$ takes values in $\bS_+^2$ by the strong convergence, and $\bfd_c(\bfx)=\bfd_0(\bfO_c^\top(\bfx-\bfh_c))$ for a.e.\ $\bfx\in\cS_c$, completing the proof of~(i).

For~(ii), consider $\bfzeta\in\cD_{\Omega_c}$. 
By the convergence $\cS(t)\to\cS_c$, the support of $\widetilde\bfzeta$ has positive distance from $\cS(t)$ for all $t$ close to $T_c$. 
Multiplied by a cut-off function in time, $\widetilde\bfzeta$ is therefore an admissible test field in the momentum equation. 
Thanks to \eqref{eq:energy-main}, the convective term and the Ericksen stress are bounded in $\rL^\infty_t\rL^1_{\bfx}$, and the viscous term is bounded in $\rL^2_{t,\bfx}$, so the distributional derivative of~$t\mapsto(\bfu(t),\widetilde{\bfzeta})_{\rL^2(\Omega)}$ belongs to $\rL^1$ near $T_c$.
Thus, the limit as $t\to T_c$ exists, and we get the bound
\[
    \left|\lim_{t\to T_c}(\bfu(t),\widetilde{\bfzeta})_{\rL^2(\Omega)}\right| \le \sup_{t<T_c}\|\bfu(t)\|_{\rL^2(\Omega)}\,\|\bfzeta\|_{\rL^2(\Omega_c)}.
\]
Since $\cD_{\Omega_c}$ is dense in $\cH_{\Omega_c}$, the limit defines a bounded linear functional on $\cH_{\Omega_c}$, and the Riesz representation theorem gives a unique $\bfu_c^+\in\cH_{\Omega_c}$ and thus proves~(ii).

With regard to~(iii), let $\bfu(t_j)\rightharpoonup\bfu_c^-$ be as in the statement. Then $\bfu_c^-\in\cH(\Omega)$ by weak closedness. On $\cS(t_j)$, we write
\[
 \bfu(t_j,\bfx)=\bfell_j+\bfomega_j\times(\bfx-\bfh(t_j)).
\]
By \eqref{eq:energy-main} and the proof of \autoref{lem:rigid velocity energy control}, the sequence $(\bfell_j,\bfomega_j)$ is bounded, and after extraction it converges. 
Since for sufficiently large $j$, every compact subset of $\cS_c$ lies in $\cS(t_j)$, testing the weak convergence against fields with compact support in the set $\cS_c$ reveals that $\bfu_c^-\in\cH_{\cS_c}(\Omega)$ in the sense of \eqref{eq:rigid kinetic space general}. 
Moreover, (ii) yields $(\bfu_c^+,\bfzeta)_{\rL^2(\Omega_c)}=(\bfu_c^-|_{\Omega_c},\bfzeta)_{\rL^2(\Omega_c)}$ for every~$\bfzeta\in\cD_{\Omega_c}$.
By density of $\cD_{\Omega_c}$ in $\cH_{\Omega_c}$ and $\bfu_c^+\in\cH_{\Omega_c}$, this shows $\bfu_c^+=P_{\cH_{\Omega_c}}(\bfu_c^-|_{\Omega_c})$.

Finally, choose times $t_j\to T_c$, in the set with full measure on which the energy inequality holds and~$\bfu(t_j)$ is rigid on $\cS(t_j)$, such that $E(t_j)\to\operatorname*{ess\,liminf}_{t\to T_c}E(t)$. 
After extraction, we infer that $\bfu(t_j)\rightharpoonup\bfu_c^-$ in $\rL^2(\Omega)$.
In view of (iii), we have the orthogonal decomposition
\[
    \|\bfu_c^-\|_{\rL^2(\Omega)}^2=\|\bfu_c^+\|_{\rL^2(\Omega_c)}^2+\|\bfu_c^-\|_{\rL^2(\cS_c)}^2+\|\bfu_c^-|_{\Omega_c}-\bfu_c^+\|_{\rL^2(\Omega_c)}^2.
\]
Thus, the weak lower semicontinuity of the norm yields the estimate
\[
    \|\bfu_c^+\|_{\rL^2(\Omega_c)}^2+\|\bfu_c^-\|_{\rL^2(\cS_c)}^2+\|\bfu_c^-|_{\Omega_c}-\bfu_c^+\|_{\rL^2(\Omega_c)}^2\le \liminf_{j\to\infty}\|\bfu(t_j)\|_{\rL^2(\Omega)}^2.
\]
Moreover, the norm of the gradient of the part of $\overline\bfd(t)$ defined on the solid equals $\|\nabla\bfd_0\|_{\rL^2(\cS_0)}$ for every $t<T_c$, and the same holds for $\bfd_c$ on $\cS_c$. 
The weak convergence from \eqref{eq:constrained collision director trace} and the weak lower semicontinuity of the $\rH^1(\Omega)$-norm then give
\[
    \|\nabla\bfd_c\|_{\rL^2(\Omega_c)}^2 \le \liminf_{j\to\infty}\|\nabla\bfd(t_j)\|_{\rL^2(\cF(t_j))}^2.
\]
Combining the last two estimates, we find \eqref{eq:constrained collision energy drop strengthened}. 
Dropping the two nonnegative impact terms yields~\eqref{eq:constrained collision energy drop}, finishing the proof.
\end{proof}

\begin{proof}[Completion of the proof of \autoref{exist:main}]
From \autoref{ssec:maximal constrained}, we obtain the solution to the moving-body problem on $[0,T_{\rm col})$, proving \autoref{exist:main}\,(i).
If $T_{\rm col}=\infty$, nothing remains to be shown, so assume~$T_{\rm col}<\infty$. 
Parts (i) and (ii) of \autoref{lem:constrained collision traces} provide the rigid body position at contact $\cS_c$, the contact value~$\bfd_c\in\rH^1(\Omega;\bS_+^2)$ of the material extension, and the post-impact velocity $\bfu_c^+\in\cH_{\Omega_c}$. 
This proves \autoref{exist:main}\,(ii).
With regard to \autoref{exist:main}\,(iii), we apply \autoref{prop:postcontact constrained arbitrary domain} on $\Omega_c$ with initial data $(\bfu_c^+,\bfd_c|_{\Omega_c})$, keep the body at $\cS_c$ for~$t\ge T_{\rm col}$, and extend the post-impact velocity by zero to $\cS_c$.
The energy inequality \eqref{eq:global constrained energy main} follows from the sharpened inequality in \autoref{prop:impact law}, proved below, since~$\mathcal J_{\rm imp}\ge0$.
\end{proof}

\begin{proof}[Proof of \autoref{prop:impact law}]
\autoref{lem:constrained collision traces}\,(iii) gives \eqref{eq:impact projection main} for every $\bfu_c^-$, and the orthogonal decomposition
\[
    \|\bfu_c^-\|_{\rL^2(\Omega)}^2=\|\bfu_c^+\|_{\rL^2(\Omega_c)}^2+\|\bfu_c^-\|_{\rL^2(\cS_c)}^2+\|\bfu_c^-|_{\Omega_c}-\bfu_c^+\|_{\rL^2(\Omega_c)}^2
\]
corresponds precisely to \eqref{eq:impact energy decomposition main}. 
For the energy inequalities, let $\mathcal D$ be the dissipation from \autoref{exist:main}\,(iii), and 
choose~${t_j\to T_{\rm col}}$ and $\bfu_c^-$ as in \autoref{lem:constrained collision traces}\,(iv), so that $E(t_j)\to\operatorname*{ess\,liminf}_{t\to T_{\rm col}}E(t)$ and $\bfu(t_j)\rightharpoonup\bfu_c^-$. By the projection formula \eqref{eq:impact projection main}, the two impact terms in \eqref{eq:constrained collision energy drop strengthened} add up to $\mathcal J_{\rm imp}(\cS_c,\bfu_c^-)$, and therefore
\[
 \mathcal E_{\Omega_c}(\bfu_c^+,\bfd_c)+\mathcal J_{\rm imp}(\cS_c,\bfu_c^-)\le\lim_{j\to\infty}E(t_j)=\operatorname*{ess\,liminf}_{t\to T_{\rm col}}E(t),
\]
which is \eqref{eq:impact energy drop main}.

It remains to prove the sharpened global energy inequality. Before contact, the energy inequality~\eqref{eq:energy-main} at the times $t_j$ reads $E(t_j)+\mathcal D(0,t_j)\le E_0$. As $j\to\infty$, we have $E(t_j)\to\operatorname*{ess\,liminf}_{t\to T_{\rm col}}E(t)$ and, by monotone convergence, $\mathcal D(0,t_j)\to\mathcal D(0,T_{\rm col})$. Together with \eqref{eq:impact energy drop main}, this gives
\[
 \mathcal E_{\Omega_c}(\bfu_c^+,\bfd_c)+\mathcal J_{\rm imp}(\cS_c,\bfu_c^-)+\mathcal D(0,T_{\rm col})\le E_0.
\]
After contact, the solution satisfies \eqref{eq:postcontact constrained energy} on $\Omega_c$ with initial energy $\mathcal E_{\Omega_c}(\bfu_c^+,\bfd_c)$, that is,
\[
 \mathcal E_{\Omega_c}(\bfu(t),\bfd(t))+\mathcal D(0,t)-\mathcal D(0,T_{\rm col})\le\mathcal E_{\Omega_c}(\bfu_c^+,\bfd_c)
 \qquad\text{for a.e.\ }t>T_{\rm col},
\]
since $\mathcal D(0,t)-\mathcal D(0,T_{\rm col})$ is the dissipation after contact. Adding the two inequalities yields \eqref{eq:global constrained energy main} with $\mathcal J_{\rm imp}(\cS_c,\bfu_c^-)$ added to its left-hand side.
\end{proof}

\medskip
{\bf Acknowledgements.}
F.B.~acknowledges the support by the German National Academy of Sciences Leopoldina through the Leopoldina Fellowship Program with grant number~LPDS 2024-07. V.N. acknowledges partial support from an AMS--Simons travel grant.

\appendix
\section{Proofs of auxiliary results}\label{sec:proofs aux lemmas}

\subsection{Approximation of the director fields}

The lemma below asserts that an $\rH^1$-director field with values in $\bS_+^2$ can be approximated by smooth director fields with values in $\bS_+^2$.

\begin{lem}\label{lem:hemisphere preparation}
Let $\Omega\subset\R^3$ be a bounded $\rC^2$-domain, and let $\bfd\in \rH^1(\Omega;\bS_+^2)$. 
Then there are maps~$\bfd_j\in \rC^\infty(\overline\Omega;\bS_+^2)\subset \rH^2(\Omega;\bS_+^2)$ such that $\bfd_j\rightarrow\bfd$ strongly in $\rH^1(\Omega;\R^3)$.

The assertion remains valid under fixed rotations of the target hemisphere.
\end{lem}

\begin{proof}
It is enough to consider $\bS_+^2$
To this end, we write $\bfd=(\bfp,q)$ with $\bfp=(d_1,d_2)$ and $q=d_3\ge0$.  
Then $|\bfp|^2+q^2=1$ and, a.e., we get $\bfp\cdot\nabla\bfp+q\nabla q=0$.
For $0<\kappa<1$, we set $\bfp^\kappa=(1-\kappa)\bfp$, $q^\kappa=\sqrt{1-|\bfp^\kappa|^2}$, and $\bfd^\kappa=(\bfp^\kappa,q^\kappa)$, so we obtain $q^\kappa\ge\sqrt{1-(1-\kappa)^2}>0$ as well as
\[
 \nabla q^\kappa
 =\frac{(1-\kappa)^2q}{\sqrt{q^2+(2\kappa-\kappa^2)|\bfp|^2}}\nabla q.
\]
Note that the coefficient of $\nabla q$ is bounded by one and converges to one wherever $q>0$. 
On the other hand, on the set where the nonnegative Sobolev function $q$ vanishes, one has $\nabla q=0$ a.e.  
From dominated convergence, it then follows that $\bfd^\kappa\to\bfd$ strongly in $\rH^1(\Omega)$ as $\kappa\to0$.

Now we fix $\kappa$, extend $\bfp^\kappa$ to an $\rH^1(\R^3;\R^2)$-map, and compose the extension with the metric projection onto the closed disk of radius $1-\kappa/2$.  
Note that this does not change the map on $\Omega$.  
By mollification, we get smooth maps $\bfp_{\kappa,m}$ taking values in that same convex disk and converging to~$\bfp^\kappa$ in $\rH^1(\Omega)$. 
From the smoothness of the map $\bfz\longmapsto\bigl(\bfz,\sqrt{1-|\bfz|^2}\bigr)$ with bounded derivative on $\{|\bfz|\le1-\kappa/2\}$, we infer the strong convergence of $\bfd_{\kappa,m}$ to $\bfd^\kappa$ in $\rH^1(\Omega)$.  
Choosing a suitable diagonal sequence with $\kappa\to0$, $m\to\infty$ proves the claim.
\end{proof}

\subsection{Test-space approximation and material gluing}\label{ssec:test transport auxiliary}

We first prove density of $\cT_{\bfu}^{\rm col}$ as introduced in \eqref{eq:collar momentum tests} in $cT_{\bfu}$

\begin{lem}\label{lem:collar density}
Assume that the rigid body $\cS(t)=\bfh(t)+\bfO(t)\cS_0$ maintains a positive distance from $\partial\Omega$, where $\bfh$ and~$\bfO$ belong to $\rW^{1,\infty}(0,T)$ and $\partial\cS_0$ is of class $\rC^2$. 
Then for every $\bfzeta\in\cT_{\bfu}$ there are $r_k\to0$ and $\bfzeta_k\in\cT_{\bfu}^{\rm col}$, with~$\bD(\bfzeta_k)=0$ on $\cS(t)_{r_k}$, such that
\begin{equation}\label{eq:collar density convergence}
    \|\bfzeta_k-\bfzeta\|_{\rL^\infty(0,T;\rW^{1,\infty}(\Omega))}+\|\partial_t(\bfzeta_k-\bfzeta)\|_{\rL^1(0,T;\rL^2(\Omega))}\rightarrow0.
\end{equation}
\end{lem}

\begin{proof}[Proof of \autoref{lem:collar density}]
We start by fixing $\bfzeta\in\cT_{\bfu}$.  
Since $\cS(t)$ is connected and
$\bD(\bfzeta(t))=0$ there, it follows that there exists a unique affine rigid field $\bfR(t,\bfx)=\bfa(t)+\bfB(t)\bfx$, with $\bfB(t)^\top=-\bfB(t)$, and which agrees with $\bfzeta(t)$ on $\cS(t)$.  
The coefficients $\bfa(t)$ and $\bfB(t)$ can be written by the mean and angular-momentum formulas on $\cS(t)$, showing that they belong to $\rW^{1,\infty}(0,T)$. 
We now set~$\bfq=\bfzeta-\bfR$.  
Then $\bfq=0$ on the open set $\cS(t)$.  
By the smoothness of $\bfzeta$ in space, every first spatial derivative of $\bfq$ has zero trace on $\partial\cS(t)$.  
For $d_t(\bfx):=\dist(\bfx,\cS(t))$ and~$d_t(\bfx)<r_0$, it then follows that
\begin{equation}\label{eq:flat rigid test}
    |\bfq(t,\bfx)|\le C d_t(\bfx)^2, \qquad |\nabla\bfq(t,\bfx)|\le C d_t(\bfx).
\end{equation}
Differentiating the identity $\bfq(t,\bfh(t)+\bfO(t)\bfy)=0$ for $\bfy\in\partial\cS_0$ at a.e.\ time and using $\nabla\bfq=0$ on the interface, for a.e.\ $t$, we find
\begin{equation}\label{eq:flat rigid test time}
    |\partial_t\bfq(t,\bfx)|\le C d_t(\bfx)
\end{equation}
in the same tubular neighborhood. 
Note that the rigidity of the motion and its Lipschitz regularity yield that the constants are uniform in time.

Next, we choose a smooth cutoff $\vartheta\in \rC^\infty(\R)$ with $\vartheta=0$ on $(-\infty,1]$ and $\vartheta=1$ on $[2,\infty)$, and we set $\vartheta_r(t,\bfx)=\vartheta(d_t(\bfx)/r)$.  
The field $\bfzeta_r^0=\bfR+\vartheta_r\bfq$ equals $\bfR$ on $\cS(t)_r$ and equals $\bfzeta$ outside $\cS(t)_{2r}$, but $\mdiv\bfzeta_r^0=f_r:=\nabla\vartheta_r\cdot\bfq$. 
The defect is supported in the tubular annulus $A_r(t)=\cS(t)_{2r}\setminus\overline{\cS(t)_r}$.
From~\eqref{eq:flat rigid test}, we deduce that $\|f_r\|_{\rL^\infty((0,T)\times\Omega)}\le Cr$ and $\int_{A_r(t)}f_r(t) \srd \bfx=0$, where we note that the second identity follows from the divergence theorem.
The cutoff vanishes on the inner parallel boundary, while the flux of $\bfq$ through the outer parallel boundary equals its flux through~$\partial\cS(t)$ and is therefore zero because $\mdiv\bfq=0$ and $\bfq=0$ on $\partial\cS(t)$.

Fix $r_0>0$ smaller than the uniform tubular radius and the distance from the moving body to~$\partial\Omega$. 
Instead of solving the divergence equation on the shrinking annulus, solve it on the collar with fixed-width collar $A_r^*(t):=\cS(t)_{r_0}\setminus\overline{\cS(t)_r}$.
Then the defect $f_r$ is supported in $A_r(t)\subset A_r^*(t)$, and its mean over $A_r^*(t)$ is zero by the divergence theorem.
The estimates \eqref{eq:flat rigid test} imply, for every fixed $p\in(3,\infty)$, the bounds $\|f_r(t)\|_{\rL^\infty(A_r^*(t))}\le Cr$ and $\|f_r(t)\|_{\rW^{1,p}(A_r^*(t))}\le C r^{1/p}$, since first derivatives of $f_r=\nabla\vartheta_r\cdot\bfq$ are $O(1)$ and are supported in a layer of volume $O(r)$.

We observe that the collars $A_r^*(t)$, $0<r<r_0/4$, form a uniformly $\rC^2$-family. 
After undoing the rigid motion, their inner and outer boundaries are parallel surfaces of $\partial\cS_0$ at offsets in $[0,r_0]$, with $r_0$ fixed.  
The Bogovski\u{\i} construction on a finite $\rC^2$-atlas then gives operators $\mathcal B_r(t)$ with constants independent of $r$ and $t$ such that $\mdiv\mathcal B_r(t)f_r=f_r$, $\mathcal B_r(t)f_r\in \rW^{2,p}(A_r^*(t);\R^3)\cap \rW^{1,p}_0(A_r^*(t);\R^3)$, and
\begin{equation}\label{eq:bog collar estimates}
    \|\mathcal B_r f_r\|_{\rW^{2,p}(A_r^*(t))}\le C\|f_r\|_{\rW^{1,p}(A_r^*(t))}\le C r^{1/p},
\end{equation}
where we refer to \cite{Bog:79,Gal:94} for the Bogovski\u{\i} operator and its higher-order estimates.  
Due to $p>3$, Sobolev embedding yields $\|\mathcal B_r f_r\|_{\rW^{1,\infty}(A_r^*(t))}
 \le C r^{1/p}\rightarrow0$.
We then extend the correction by zero outside $A_r^*(t)$. 
As its trace vanishes on the inner boundary, this extension does not alter the rigid field on $\cS(t)_r$.

The time derivative is estimated separately. 
A differentiation of $f_r$ produces only $O(1)$ terms on the~$O(r)$ transition layer, by \eqref{eq:flat rigid test} and \eqref{eq:flat rigid test time}. 
Moreover, we have $\nabla f_r=O(1)$ on the same layer, so for a.e.\ $t$, we get $\|\partial_t f_r(t)\|_{\rL^2(\Omega)}+\|\nabla f_r(t)\|_{\rL^2(\Omega)}\le Cr^{1/2}$.

For a uniform control of the time derivative, we write the Bogovski\u{\i} correction on a fixed reference collar. 
To this end, let $\bfZ_t(\bfy):=\bfh(t)+\bfO(t)\bfy$, $A_r^{*,0}:=\cS_{0,r_0}\setminus\overline{\cS_{0,r}}$, and invoke the Bogovski\u{\i} operator~$\mathcal B_r^0$ on $A_r^{*,0}$ with the uniform bounds above.
We set $\widetilde f_r(t,\bfy):=f_r(t,\bfZ_t(\bfy))$ as well as~$\bfw_r(t,\bfZ_t(\bfy)):=\bfO(t)\,\mathcal B_r^0\widetilde f_r(t,\bfy)$.
Since the divergence is preserved under rotation, we get~$\mdiv_{\bfx}\bfw_r=f_r$. 
With $\mathbf V_R(t,\bfx):=\bfh'(t)+\bfO'(t)\bfO(t)^\top(\bfx-\bfh(t))$, it then follows that
\[
  \partial_t\widetilde f_r
  =\bigl(\partial_t f_r+\mathbf V_R\cdot\nabla f_r\bigr)\circ \bfZ_t.
\]
The rigid velocity $\mathbf V_R$ is uniformly bounded on the fixed collar, and $\bfZ_t$ maps $A_r^{*,0}$ isometrically onto~$A_r^*(t)$. Therefore
\[
 \|\partial_t\widetilde f_r(t)\|_{\rL^2(A_r^{*,0})}
 =\bigl\|\bigl(\partial_t f_r+\mathbf V_R\cdot\nabla f_r\bigr)(t)\bigr\|_{\rL^2(A_r^*(t))}\le Cr^{1/2}.
\]
Differentiating the formula for $\bfw_r$ at fixed $\bfx$ gives a rotation term, a term containing $\mathcal B_r^0\partial_t\widetilde f_r$, and the Eulerian transport term containing $\nabla\mathcal B_r^0\widetilde f_r$. 
The uniform $\rL^2\to\rH_0^1$-estimate for the Bogovski\u{\i} operators $\mathcal B_r^0$ on the reference collars, see \cite{Bog:79} and \cite[Chap.~III.3]{Gal:94}, thus gives $\|\partial_t\bfw_r(t)\|_{\rL^2(\Omega)}\le Cr^{1/2}$, and we may take $\mathcal B_r(t)f_r=\bfw_r$. In particular, the time estimate does not rely on differentiating an abstract operator on the moving collar.

Define $\widehat\bfzeta_r:=\bfzeta_r^0-\mathcal B_r(t)f_r$.
Then $\mdiv\widehat\bfzeta_r=0$ and it is rigid on $\cS(t)_r$. The Bogovski\u{\i} correction is supported in the fixed collar~$\cS(t)_{r_0}$ and tends to zero there in $\rW^{1,\infty}$. 
It need not vanish on all of~$\cS(t)_{r_0}\setminus\cS(t)_{2r}$.  
From the estimates~\eqref{eq:flat rigid test}--\eqref{eq:bog collar estimates}, it follows that
\begin{equation}\label{eq:pre collar convergence}
    \|\widehat\bfzeta_r-\bfzeta\|_{\rL^\infty(0,T;\rW^{1,\infty}(\Omega))}\to0,\qquad \|\partial_t(\widehat\bfzeta_r-\bfzeta)\|_{\rL^1(0,T;\rL^2(\Omega))}\to0,
\end{equation}
where all fields are defined on the whole container $\Omega$.
The field $\widehat\bfzeta_r$ is not necessarily smooth yet because the tubular coordinates inherit only the $\rC^2$-regularity of $\partial\cS_0$. 
We note that $\widehat\bfzeta_r$ coincides with the affine rigid field $\bfR(t)$ on $\cS(t)_r$ and has compact spatial support that has a fixed distance from $\partial\Omega$.
We then extend it by zero to~$\R^3$ and convolve in space with a symmetric mollifier of radius $\rho_r=o(r)$, where we invoke that convolution preserves the divergence-free condition.  
Moreover, on~$\cS(t)_{r/2}$ the convolution samples only the affine field~$\bfR(t)$. Symmetry of the mollifier preserves affine functions, so the spatially mollified field is still rigid there.  The estimates in \eqref{eq:pre collar convergence} are preserved. Indeed, writing $\eta_{\rho_r}$ for the spatial mollifier,
\[
 \|\eta_{\rho_r}*\widehat\bfzeta_r-\bfzeta\|_{\rW^{1,\infty}(\Omega)}
 \le
 \|\widehat\bfzeta_r-\bfzeta\|_{\rW^{1,\infty}(\Omega)}
 +\|\eta_{\rho_r}*\bfzeta-\bfzeta\|_{\rW^{1,\infty}(\Omega)},
\]
and the second term tends to zero because $\bfzeta$ is spatially smooth. The same contraction argument in~$\rL^1_t\rL^2_{\bfx}$, applied after differentiating in time, gives the second convergence in \eqref{eq:pre collar convergence}.

At this stage the time regularity is only $\rW^{1,\infty}$ because the rigid trajectory has that regularity. 
We then extend the fields slightly beyond the temporal endpoints and use mollification in time on a second scale~${\theta_r=o(r)}$.  
From the Lipschitz continuity of $\cS(t)$, we deduce that $\cS(t)_{r/4}\subset\cS(s)_{r/2}$ whenever $|t-s|\le c r$.  
Note that the time-mollified field remains rigid on~$\cS(t)_{r/4}$. 
A final time cutoff, chosen in a region where the original test already vanishes, produces~$\bfzeta_r\in\cT_{\bfu}^{\rm col}$, and \eqref{eq:pre collar convergence} is preserved.  
A choice of $r=r_k\to0$ proves \eqref{eq:collar density convergence} and thus completes the proof.
\end{proof}

The density result has a formulation independent of the liquid-crystal equations.

\begin{prop}\label{prop:rigid test equivalence}
Let $\cS(t)$ satisfy the hypotheses of \autoref{lem:collar density}.  
Suppose $\bfu$ is solenoidal and satisfies $\bfu\in \rL^\infty(0,T;\rL^2(\Omega;\R^3))\cap \rL^2(0,T;\rH_0^1(\Omega;\R^3))$ as well as $\mathbb A\in \rL^1((0,T)\times\Omega;\R^{3\times3})$.
Assume that the momentum identity
\[
 \int_0^T\int_\Omega
 \bigl[-\bfu\cdot\partial_t\bfzeta-(\bfu\otimes\bfu):\nabla\bfzeta
 +2\bD(\bfu):\bD(\bfzeta)\bigr] \srd \bfx \srd t
 =\int_\Omega\bfu_0\cdot\bfzeta(0) \srd \bfx
 +\int_0^T\int_\Omega\mathbb A:\nabla\bfzeta \srd \bfx \srd t
\]
holds for every $\bfzeta\in\cT_{\bfu}^{\rm col}$.  Then it holds for every
$\bfzeta\in\cT_{\bfu}$.
\end{prop}

\begin{proof}[Proof of \autoref{prop:rigid test equivalence}]
Let $\bfzeta_k$ be the sequence from \autoref{lem:collar density}.  The convergence \eqref{eq:collar density convergence} gives
\[
 \left|\int_0^T\int_\Omega \bfu\cdot\partial_t(\bfzeta_k-\bfzeta) \srd \bfx \srd t\right|
 \le \|\bfu\|_{\rL^\infty_t\rL^2_{\bfx}}
      \|\partial_t(\bfzeta_k-\bfzeta)\|_{\rL^1_t\rL^2_{\bfx}}\to0.
\]
Similarly,
\[
 \left|\int_0^T\int_\Omega (\bfu\otimes\bfu):\nabla(\bfzeta_k-\bfzeta) \srd \bfx \srd t\right|
 \le T\|\bfu\|_{\rL^\infty_t\rL^2_{\bfx}}^2
      \|\nabla(\bfzeta_k-\bfzeta)\|_{\rL^\infty_{t,\bfx}}\to0,
\]
and the viscous term tends to zero by Cauchy--Schwarz and the uniform
$\rL^\infty_t\rW^{1,\infty}_{\bfx}$-convergence.  The stress term satisfies
\[
 \left|\int_0^T\int_\Omega
 \mathbb A:\nabla(\bfzeta_k-\bfzeta) \srd \bfx \srd t\right|
 \le \|\mathbb A\|_{\rL^1((0,T)\times\Omega)}
 \|\nabla(\bfzeta_k-\bfzeta)\|_{\rL^\infty((0,T)\times\Omega)}\to0.
\]
The initial term converges uniformly in space.  Passing $k\to\infty$ proves the claim.
\end{proof}

The next lemma gives the transport identity for the material extension and shows that no distribution supported on the moving interface is produced.

\begin{lem}\label{lem:material gluing}
Let $\cS(t)=\bfh(t)+\bfO(t)\cS_0\Subset\Omega$ be a Lipschitz rigid trajectory and $\cF(t)=\Omega\setminus\overline{\cS(t)}$.  Let~$\bfd_0\in\rH^1(\Omega;\R^3)$, and let
$\bfu\in \rL^2(0,T;\rH^1_0(\Omega;\R^3))$ be solenoidal and equal in $\cS(t)$ to the rigid velocity generating the trajectory.  Suppose
$\bfd_F\in \rL^\infty(0,T;\rH^1(\cF(t);\R^3))$ has trace
\[
 \operatorname{Tr}_{\partial\cS(t)}\bfd_F(t,\bfx)
 =\bfd_0(\bfO(t)^\top(\bfx-\bfh(t))),
\]
and satisfies
$\partial_t\bfd_F+(\bfu\cdot\nabla)\bfd_F=\bfM$ in the moving fluid region, with $\bfM\in \rL^2(\cQ_{\cF})$, in the energy weak formulation.  Define
\[
 \bfd_S(t,\bfx):=\bfd_0(\bfO(t)^\top(\bfx-\bfh(t)))
 \quad(\bfx\in\cS(t)),
 \qquad
 \overline\bfd:=\bfd_F\dsone_{\cF(t)}+\bfd_S\dsone_{\cS(t)},
\]
and let $\widetilde\bfM$ be the zero extension of $\bfM$ through the solid.  
Then $\overline\bfd(t)\in \rH^1(\Omega;\R^3)$ for a.e.\ $t$ and
\begin{equation}\label{eq:material gluing identity}
 \partial_t\overline\bfd+(\bfu\cdot\nabla)\overline\bfd
 =\widetilde\bfM
 \quad\text{in }\cD'((0,T)\times\Omega).
\end{equation}
In particular, if $\overline\bfd\in \rL^\infty(0,T;\rH^1(\Omega))$ and $\widetilde\bfM\in\rL^2((0,T)\times\Omega)$, then $\partial_t\overline\bfd\in \rL^2(0,T;\rH^{-1}(\Omega))$.
\end{lem}

\begin{proof}[Proof of \autoref{lem:material gluing}]
For a.e.\ time, integration by parts in the two phases shows that the distributional spatial gradient of the glued field is its piecewise $\rL^2$-gradient, because the possible interface term is the jump of its trace multiplied by the normal, and this jump is zero.
This shows that~$\overline\bfd(t)\in\rH^1(\Omega;\R^3)$.

We justify the material equation on a compact time interval on which the gap is positive. 
To this end, let~$\bfU$ be the affine rigid velocity generating $\cS(t)$, and let $\rho(t,\bfx)=\dist(\bfx,\cS(t))$ in a fixed fluid-side tubular neighborhood.
From rigidity, it follows that $(\partial_t+\bfU\cdot\nabla)\rho=0$.
We then choose $\vartheta\in\rC^\infty(\R)$, equal to zero on $(-\infty,1]$ and to one on $[2,\infty)$, and set $\eta_r=\vartheta(\rho/r)$ in the collar, extended by one in the remaining fluid.
This shows $(\partial_t+\bfu\cdot\nabla)\eta_r=(\bfu-\bfU)\cdot\nabla\eta_r$ and $|\nabla\eta_r|\le C/r$., and the derivative is supported in $A_r(t)=\{\bfx\in\cF(t):r<\rho(t,\bfx)<2r\}$.
The trace of~$\bfu-\bfU$ vanishes at the inner interface, so that the one-dimensional Hardy estimate in tubular coordinates gives
\[
    \left\|\frac{\bfu-\bfU}{\rho}\right\|_{\rL^2(\{(t,\bfx):0<\rho(t,\bfx)<r_0\})}<\infty.
\]
Let $\bfxi\in\rC_c^\infty((0,T)\times\Omega;\R^3)$.
We then use $\eta_r\bfxi$ in the fluid molecular-field equation and observe that it is an admissible test function.
Integration by parts in the convective term, we get
\[
    \int_0^T\int_{\cF(t)}\eta_r\bigl[-\bfd_F\cdot(\partial_t\bfxi+\bfu\cdot\nabla\bfxi)-\bfM\cdot\bfxi\bigr]\srd\bfx\srd t=
 \int_0^T\int_{A_r(t)}
 \bfd_F\cdot\bfxi\,(\partial_t+\bfu\cdot\nabla)\eta_r
 \srd\bfx\srd t.
\]
In this identity, the absolute value of the right-hand side is bounded by
\[
    C\|\bfxi\|_{\rL^\infty((0,T)\times\Omega)}\|\bfd_F\|_{\rL^2(\{(t,\bfx):\bfx\in A_r(t)\})}\left\|\frac{\bfu-\bfU}{\rho}\right\|_{\rL^2(\{(t,\bfx):0<\rho<r_0\})},
\]
which tends to zero by absolute continuity of the first $\rL^2$ integral. 
For the other terms, we may use dominated convergence.
This proves the fluid material identity for arbitrary $\bfxi$ without an interface term. 
The solid branch satisfies the analogous identity with zero source by the rigid change of variables. 
Adding
them proves \eqref{eq:material gluing identity}.
The argument is local in time. 
The resulting distributional identity holds on the whole open pre-contact interval.

Finally, we use $\bfu\in \rL^2_t\rL^6_{\bfx}$ and $\nabla\overline\bfd\in \rL^\infty_t\rL^2_{\bfx}$ to get $\bfu\cdot\nabla\overline\bfd\in \rL^2(0,T;\rL^{3/2}(\Omega))\hookrightarrow \rL^2(0,T;\rH^{-1}(\Omega))$.
The last assertion of this lemma then follows from \eqref{eq:material gluing identity}.
\end{proof}

\subsection{Auxiliary estimates for the Galerkin approximation}\label{ssec:aux ests Galerkin}

The Galerkin construction uses the following maximum principle for the relaxed director equation.

\begin{lem}\label{lem:relaxed dir cond}
Let $T>0$, let $\bfu\in\rL^2(0,T;\cV(\Omega))$, let $a\in\rL^\infty((0,T)\times\Omega)$ satisfy $a\ge0$, and let~$\bb\in\rL^\infty((0,T)\times\Omega;\R^3)$ satisfy $|\bb|\le1$ and $(\bb)_3\ge0$ a.e.
Suppose that $\bfd$ is a sufficiently regular weak solution of
\begin{equation}\label{eq:penalized maximum principle equation}
 \del_t\bfd+\bfu\cdot\nabla\bfd
 =\Delta\bfd+\frac1{\eps^2}(1-|\bfd|^2)\bfd-a(\bfd-\bb)
\end{equation}
with initial and Dirichlet boundary data satisfying $|\bfd|\le1$ and $d_3\ge0$. 
Then for a.e.\ $(t,\bfx)\in(0,T)\times\Omega$, it holds that $|\bfd(t,\bfx)|\le1$ and $d_3(t,\bfx)\ge0$.
In particular, the conclusion applies to the Galerkin director equation \eqref{eq:k-director}, with $a=n\varphi^k$ and~$\bb=\bb^k$.
\end{lem}

\begin{proof}
We set $w=|\bfd|^2$. 
From \eqref{eq:penalized maximum principle equation}, it follows that
\[
 (\del_t+\bfu\cdot\nabla-\Delta)w=-2|\nabla\bfd|^2+\frac{2}{\eps^2}(1-w)w-2a(w-\bfd\cdot\bb).
\]
On the set $\{w>1\}$, the bound $|\bb|\le1$ gives
$w-\bfd\cdot\bb\ge w-\sqrt w\ge0$, while $(1-w)w\le0$. 
Testing the scalar inequality by $(w-1)_+$, using $\mdiv\bfu=0$ and the fact that $(w-1)_+$ has zero trace and zero initial value, we conclude that $(w-1)_+=0$.

For $q=d_3$, after the first step, we have $\del_t+\bfu\cdot\nabla-\Delta)q=\frac1{\eps^2}(1-|\bfd|^2)q-aq+a(\bb)_3$.
Testing this identity by $-q_-$ and using $(\bb)_3\ge0$, $a\ge0$, and $0\le1-|\bfd|^2\le1$, we find that
\[
    \frac12\frac{\rd}{\rd t}\|q_-\|_{\rL^2(\Omega)}^2+\|\nabla q_-\|_{\rL^2(\Omega)}^2 \le \frac1{\eps^2}\|q_-\|_{\rL^2(\Omega)}^2.
\]
Since $q_-(0)=0$ and the negative part has zero boundary trace, an application of Gr\"onwall's lemma yields $q_-=0$. 
The argument is first justified at the Galerkin level and then passes to the weak limits by the strong $\rL^2$-convergence available in the construction. 
We also refer here to \cite[Lemmas~2.1--2.2]{LW:16}.
\end{proof}

\begin{proof}[Proof of \autoref{prop:k-weak sols}]
We construct solutions to the semi-Galerkin problem introduced in \autoref{ssec:semi-Galerkin approx}, verify the energy identity~\eqref{eq:k-energy-identity}, and prove the separation estimate~\eqref{eq:lipschitz flow separation}.
Note that the construction uses the semi-Galerkin scheme of \cite[Sec.~2]{LL:95}, coupled to the transported solid and lifting introduced in \autoref{ssec:semi-Galerkin approx}.

We start with a given velocity field $\overline{\bfu}^k\in \rC([0,T];\cV_k)$.
Next, we construct the flow map $\overline{\bfX}^k$ as a solution of \eqref{eq:k-flow-map}, the inverse flow map $\overline{\bfY}^k$ by \eqref{eq:k-inverse flow-map}, the transported solid indicator $\overline{\varphi}^k$ by~\eqref{eq:k-phi}, and the lifting vector field $\overline{\bb}^k$ by \eqref{eq:k-lift}.
With these quantities fixed, the semilinear parabolic problem~\eqref{eq:k-director} has a unique energy solution on a short time interval. We denote it by $\overline{\bfd}^k$ and define $\overline{\bfM}^{n,k}$ by \eqref{eq:k-M-def}.

The velocity field $\bfu^k$ is constructed by the following fixed-point argument.
We use the $\rL^2$-orthonormal basis $\{\bfv_j\}_{j\in\N}\subset X_u=\cV(\Omega)\cap\rH^3(\Omega;\R^3)$ of $\cH(\Omega)$ made precise in \autoref{ssec:semi-Galerkin approx}. 
Its finite spans are dense in $X_u$ in the $\rH^3$-norm, and every fixed basis vector belongs to $\rW^{1,\infty}(\Omega;\R^3)$.
We recall that $\cV_k:=\Span\{\bfv_1,\ldots,\bfv_k\}$ and consider the associated orthogonal projection $\cP_k \colon \cH(\Omega) \to \cV_k$.
We then look for a solution of the form as in \eqref{eq:approx vel}, i.e., $\bfu^k(t,\bfx)=\sum_{j=1}^k a_j^k(t)\bfv_j(\bfx)$, satisfying for all $l=1,\dots,k$, the identity
\begin{equation}\label{eq:k-momentum fixed point}
    \begin{aligned}
        &(\del_t\bfu^k,\bfv_l)_{\rL^2(\Omega)}
        +\int_{\Omega} ((\bfu^k\cdot\nabla)\bfu^k)\cdot\bfv_l\srd\bfx
        +\int_{\Omega} (2+n\overline{\varphi}^k)\bD(\bfu^k):\bD(\bfv_l)\srd\bfx +\delta\,a_u(\bfu^k,\bfv_l) \\
        &\quad=
        -\int_{\Omega}
        (\bfv_l\cdot\nabla)\overline{\bfd}^k\cdot
        \bigl(\overline{\bfM}^{n,k}+n\overline{\varphi}^k(\overline{\bfd}^k-\overline{\bb}^k)\bigr)\srd\bfx.
    \end{aligned}
\end{equation}
We formulate the fixed-point step with the transported characteristic function and the director penalization included. 
Fix $R>2\|\cP_k\bfu_0\|_{\rL^2(\Omega)}$ and let
\[
   B_R:=\{\overline{\bfu}^k\in\rC([0,T];\cV_k):
        \overline{\bfu}^k(0)=\cP_k\bfu_0,
        \ \|\overline{\bfu}^k\|_{\rC_t\rL^2_{\bfx}}\le R\}.
\]
Note that all norms on the finite-dimensional space $\cV_k$ are equivalent.
In particular, for $\overline{\bfu}^k\in B_R$, the associated flow and inverse flow exist on a time interval depending only on $k$ and $R$, and are uniformly bounded in $\rC^1$ there. 
For the transported indicator $\varphi^k$, it follows that $0\le\overline\varphi^k\le1$, while $\overline\bb^k=\bfd_0\circ\overline\bfY^k$ is bounded in all Sobolev norms needed in the parabolic problem, amd the constants depend on $k$, $R$, $T$ as well as the data but not on the particular element of $B_R$.

For these prescribed coefficients, the director equation gives rise to a semilinear parabolic Dirichlet problem with bounded zeroth-order coefficient $n\overline\varphi^k$. 
Employing standard Galerkin or monotonicity theory, we find a unique energy solution on the chosen interval. 
Moreover, \autoref{lem:relaxed dir cond} yields $|\overline{\bfd}^k|\le1$ and $(\overline{\bfd}^k)_3\ge0$. 
Testing the director equation by its molecular field $\overline\bfM^{n,k}$ and using the bounds for the prescribed drift and lift yields, on sufficiently short intervals, the estimate
\[
    \|\overline\bfd^k\|_{\rL^\infty_t\rH^1_{\bfx}}+\|\overline\bfd^k\|_{\rL^2_t\rH^2_{\bfx}}+\|\overline\bfM^{n,k}\|_{\rL^2_{t,{\bfx}}}\le C_{k,n,\delta,\eps,R,\bfu_0,\bfd_0}.
\]
We observe that the right-hand side of \eqref{eq:k-momentum fixed point} belongs to
$\rL^2(0,T)$, while the remaining coefficients are continuous or bounded in time. 
Hence \eqref{eq:k-momentum fixed point} is a Carath\'eodory system of finite-dimensional ODEs for the coefficient vector $(a_1^k,\ldots,a_k^k)$, and it has a unique absolutely continuous local solution. 
The estimates below give $\bfu^k\in \rW^{1,2}(0,T;\cV_k)\subset \rC([0,T];\cV_k)$. 
Norm equivalence and the preceding bounds give
\[
    \sup_{0\le t\le T}\|\bfu^k(t)-\cP_k\bfu_0\|_{\rL^2(\Omega)} \le C_{k,n,\delta,\eps,R,\bfu_0,\bfd_0}\,T^{\theta}
\]
for some $\theta>0$. Thus, after reducing $T$ if necessary, the map $\cT:B_R\to B_R$, $\cT(\overline\bfu^k):=\bfu^k$, is well-defined.

We verify continuity of the fixed-point map without assuming additional
smoothness of the spectral basis. Fix $k,n,\delta,\eps$ and suppose
$\overline{\bfu}_j^k\to\overline{\bfu}^k$ in
$\rC([0,T];\cV_k)$. Since $\cV_k\subset\rH^3(\Omega;\R^3)$ is
finite dimensional, this convergence also holds in
$\rC([0,T];\rH^3(\Omega;\R^3))$ and hence in
$\rC([0,T];\rW^{1,\infty}(\Omega;\R^3))$.
Stability of the flow and its inverse gives $\overline\varphi_j^k\to\overline\varphi^k$ in $\rC([0,T];\rL^p(\Omega))$ whenever $p<\infty$ and $\overline\bb_j^k\to\overline\bb^k$ in~$\rC([0,T];\rL^\infty(\Omega;\R^3))$.
Only these convergences of the transported coefficients are needed in
the physical director equation.

For brevity in this paragraph, denote the prescribed velocities by $\bfv_j,\bfv$, their transported coefficients by $\phi_j,\phi$ and $\bb_j,\bb$, and the corresponding director solutions by $\bfd_j,\bfd$.
Set $\bfr_j=\bfd_j-\bfd$.
Then $\bfr_j$ has zero initial value and
zero outer trace, and satisfies
\[
    \begin{aligned}
        \partial_t\bfr_j+(\bfv_j\cdot\nabla)\bfr_j-\Delta\bfr_j+n\phi_j\bfr_j
         &=-((\bfv_j-\bfv)\cdot\nabla)\bfd-\bigl(\bff_\eps(\bfd_j)-\bff_\eps(\bfd)\bigr)\\
        &-n(\phi_j-\phi)(\bfd-\bb)+n\phi_j(\bb_j-\bb).
    \end{aligned}
\]
Testing by $\bfr_j$, using $\mdiv\bfv_j=0$, the maximum principle,
and Lipschitz continuity of $\bff_\eps$ on the closed unit ball,
gives
\[
\begin{aligned}
 \frac12\frac{\rd}{\rd t}\|\bfr_j\|_{\rL^2(\Omega)}^2
 +\|\nabla\bfr_j\|_{\rL^2(\Omega)}^2
 \le{}& C_{n,\eps}\|\bfr_j\|_{\rL^2(\Omega)}^2
 +C\|\bfv_j-\bfv\|_{\rL^\infty(\Omega)}^2
          \|\nabla\bfd\|_{\rL^2(\Omega)}^2\\
 &+C n^2\bigl(\|\phi_j-\phi\|_{\rL^2(\Omega)}^2
                    +\|\bb_j-\bb\|_{\rL^2(\Omega)}^2\bigr).
\end{aligned}
\]
Gr\"onwall's inequality then yields $\bfd_j\to\bfd$ in $\rC([0,T];\rL^2(\Omega;\R^3)) \cap\rL^2(0,T;\rH^1(\Omega;\R^3))$.
For each fixed velocity basis vector $\bfv_l$, write the director force in stress form. 
With $F_{j,l}(t):=\int_\Omega(\nabla\bfd_j)^\top\nabla\bfd_j:\nabla\bfv_l\srd\bfx$,  we have
\[
 \|F_{j,l}-F_l\|_{\rL^1(0,T)}
 \le{}\|\nabla\bfv_l\|_{\rL^\infty(\Omega)}
  \bigl(\|\nabla\bfd_j\|_{\rL^2((0,T)\times\Omega)}
       +\|\nabla\bfd\|_{\rL^2((0,T)\times\Omega)}\bigr)\cdot
   \|\nabla\bfr_j\|_{\rL^2((0,T)\times\Omega)}\rightarrow0.
\]
Thanks to the convergence of $\phi_j$ and the boundedness of the gradients of the fixed basis functions, the viscosity matrices also converge.
Continuous dependence in the integral form of the finite-dimensional momentum ODE now gives $\cT(\overline\bfu_j^k)\to\cT(\overline\bfu^k)$ in $\rC([0,T];\cV_k)$.
The coefficient derivatives are uniformly bounded in $\rL^2(0,T)$ on the fixed-point ball, which implies equicontinuity. 
We may thus apply the Arzel\`a--Ascoli theorem, which yields that $\cT(\overline\bfu^k)$ is compact, so Schauder's theorem implies in turn that there is a fixed point.

At the fixed point, the flow $\bfX^k$, indicator $\varphi^k$, lift $\bb^k$, director $\bfd^k$, and velocity $\bu^k$ satisfy the equations \eqref{eq:k-phi}--\eqref{eq:k-momentum}. 
The exact cancellation shown in \eqref{eq:k-energy-proof-fluid}--\eqref{eq:k-penalty-energy-derivative} in the proof of \autoref{prop:k-energy-identity} gives \eqref{eq:k-energy-identity}.
We observe that these bounds are independent of the local fixed-point radius. 
Since the director problem is posed on the fixed container and the velocity takes values in the finite-dimensional space $\cV_k$, the construction can be restarted, and the Galerkin solution exists on every bounded time interval.

Finally, we prove the separation estimate \eqref{eq:lipschitz flow separation} for the solids $\cS^k(t)=\bfX^k_{0,t}(\cS_0)$ resulting from the Galerkin approximation.
Let~${\bfx\in\Omega}$, and let $\bfz\in\partial\Omega$ be a nearest boundary point. The half-open segment from $\bfx$ to $\bfz$ lies in $\Omega$, and $\bfu^k(t)$ vanishes on $\partial\Omega$, so $|\bfu^k(t,\bfx)|\le\|\nabla\bfu^k(t)\|_{\rL^\infty(\Omega)}\dist(\bfx,\partial\Omega)$.
For $\bfy\in\cS_0$, the function $g(t):=\dist(\bfX^k_{0,t}(\bfy),\partial\Omega)$ is Lipschitz and satisfies, for a.e.\ $t$,the estimate $|g'(t)|\le|\bfu^k(t,\bfX^k_{0,t}(\bfy))|\le\|\nabla\bfu^k(t)\|_{\rL^\infty(\Omega)}\,g(t)$.
Since $g(0)\ge\delta_0$ and $\delta\|\bfu^k\|_{\rL^2(0,t;X_u)}^2\le E_k^{n,\delta,\eps}(0)\le E_{\rm pen,0}$ by \eqref{eq:k-energy-identity}, by Gr\"onwall's lemma and the Cauchy--Schwarz inequality we have
\[
    g(t)\ge\delta_0\exp\Bigl(-\int_0^t\|\nabla\bfu^k(s)\|_{\rL^\infty(\Omega)}\srd s\Bigr)\ge\delta_0\exp\Bigl(-C_X\sqrt t\,\|\bfu^k\|_{\rL^2(0,t;X_u)}\Bigr)\ge\rho_\delta(t).
\]
Taking the infimum over $\bfy\in\cS_0$ proves the separation estimate for every $k$ and $n$. In particular, the approximate solid never reaches the container wall.
\end{proof}

The time-derivative estimates in \autoref{claim:t-regs u and d} are proved next.

\begin{proof}[Proof of \autoref{claim:t-regs u and d}]
From~\eqref{eq:k-momentum fixed point}, for every $\bfv_l\in\cV_k$ we have
\begin{equation}\label{eq:recalled uk}
    \begin{aligned}
        \biggl|\int_{\Omega}\del_t\bfu^k\cdot\bfv_l\srd\bfx\biggr|
        &\le \int_{\Omega}\Bigl(|(\bfu^k\otimes\bfu^k):\nabla\bfv_l| + |(2+n\varphi^k)\bigl(\bD(\bfu^k):\bD(\bfv_l)\bigr)|\\
        &\qquad + |(\bfv_l\cdot\nabla)\bfd^k\cdot(\bfM^{n,k} + n\varphi^k(\bfd^k - \bb^k))|\Bigr)\srd\bfx
        +\delta\,|a_u(\bfu^k,\bfv_l)|.
    \end{aligned}
\end{equation}
Applying H\"older's inequality together with the Sobolev embedding $\rH^2(\Omega)\hookrightarrow\rL^\infty(\Omega)$, we estimate the first three terms on the right-hand side of \eqref{eq:recalled uk} by
\begin{align}
    \int_{\Omega}|(\bfu^k\otimes\bfu^k)(t):\nabla\bfv_l|\srd\bfx 
    &\le \|\bfu^k(t)\|_{\rL^2(\Omega)}^2\,\|\bfv_l\|_{\rH^3(\Omega)}, \label{eq:ut1}\\
    \int_{\Omega}|(2+n\varphi^k)\bigl(\bD(\bfu^k(t)):\bD(\bfv_l)\bigr)|\srd\bfx 
    &\le C\,\|\nabla\bfu^k(t)\|_{\rL^2(\Omega)}\,\|\bfv_l\|_{\rH^1(\Omega)}, \tand \notag\\
    \delta\,|a_u(\bfu^k,\bfv_l)| 
    &\le \delta\,\|\bfu^k\|_{X_u}\,\|\bfv_l\|_{X_u}. \notag
\end{align}
For the last term, we recall from \eqref{eq:k-M-def} that $\bfM^{n,k}+n\varphi^k(\bfd^k-\bb^k)
    =\Delta\bfd^k+\frac{1}{\eps^2}(1-|\bfd^k|^2)\bfd^k$.
Thus, additionally using the embedding $\rH^1(\Omega) \hookrightarrow \rL^6(\Omega)$, we find that
\begin{equation}\label{eq:ut3}
    \begin{aligned}
         \int_{\Omega}|(\bfv_l\cdot\nabla)\bfd^k\cdot(\bfM^{n,k}+n\varphi^k(\bfd^k - \bb^k))|\srd\bfx
        &= \int_{\Omega}\left|(\bfv_l\cdot\nabla)\bfd^k\cdot\left(\Delta\bfd^k+\frac{1}{\eps^2}(1-|\bfd^k|^2)\bfd^k\right)\right|\srd\bfx\\
        &\le C \| \bfd^k \|_{\rH^1(\Omega)}
        \bigl(\|\bfd^k\|_{\rH^2(\Omega)} + \| \bfd^k \|_{\rL^6(\Omega)}^3\bigr) \,\|\bfv_l\|_{\rH^2(\Omega)}\\
        &\le C \| \bfd^k \|_{\rH^1(\Omega)}
        \bigl(\|\bfd^k\|_{\rH^2(\Omega)} + \| \bfd^k \|_{\rH^1(\Omega)}^3\bigr) \,\|\bfv_l\|_{\rH^2(\Omega)}.
    \end{aligned}
\end{equation}

The estimates \eqref{eq:recalled uk}--\eqref{eq:ut3} are linear in the test
field and therefore give, for every $\bfv\in\cV_k$,
\[
   |(\del_t\bfu^k(t),\bfv)_{\rL^2(\Omega)}|
   \le C\,a_k(t)\,\|\bfv\|_{X_u}
   \qquad\text{for a.e.\ }t\in(0,T),
\]
where $C$ is independent of $k$ (with $n,\delta,\eps$ fixed) and
\begin{equation*}
    \begin{aligned}
        a_k(t)
        &\coloneqq\|\bfu^{k}(t)\|_{\rL^2(\Omega)}^2 + \|\nabla \bfu^{k}(t)\|_{\rL^2(\Omega)} + \delta\|\bfu^{k}(t)\|_{X_u} + \| \bfd^k \|_{\rH^1(\Omega)}
        \bigl(\|\bfd^k\|_{\rH^2(\Omega)} + \| \bfd^k \|_{\rH^1(\Omega)}^3\bigr).
    \end{aligned}
\end{equation*}
For arbitrary $\bfzeta\in X_u$, use the spectral projection $\cP_k\bfzeta\in\cV_k$. 
Since $\del_t\bfu^k(t)\in\cV_k$ and $\cP_k$ is the~$\rL^2$-orthogonal projection, we find that $\langle\del_t\bfu^k,\bfzeta\rangle=(\del_t\bfu^k,\cP_k\bfzeta)_{\rL^2(\Omega)}$, and \eqref{eq:velocity spectral projection} yields the estimate
\[
 |\langle\del_t\bfu^k(t),\bfzeta\rangle|
 \le C a_k(t)\|\cP_k\bfzeta\|_{X_u}
 \le C a_k(t)\|\bfzeta\|_{X_u}.
\]
The energy bound in~\eqref{recall:energy} yields $a_k\in\rL^2(0,T)$, and
therefore \eqref{eq:t-regs u & d} follows with a constant independent of $k$.

Writing the director equation \eqref{eq:k-director}$_1$ in terms of $\bfM^{n,k}$ as in \eqref{eq:k-material-d-H}, for $\bfe \in \rL^2(\Omega;\R^3)$, we get
\begin{equation*}
    -\int_{\Omega}\del_t\bfd^k\cdot\bfe\srd\bfx = \int_{\Omega}\Bigl[(\bfu^k\cdot\nabla)\bfd^k\cdot\bfe - \bfM^{n,k} \cdot \bfe\Bigr] \srd \bfx.
\end{equation*}
Estimating each term via H\"older's inequality and Sobolev embeddings (cf.~\eqref{eq:ut1}--\eqref{eq:ut3} for the analogous velocity estimates), we obtain
\begin{equation*}
    \biggl|\int_{\Omega}\del_t\bfd^k\cdot\bfe\srd\bfx\biggr| \le C\,b_k(t)\,\|\bfe\|_{\rL^2(\Omega)} \quad\text{for all } \bfe\in\rL^2(\Omega;\R^3),
\end{equation*}
where $b_k(t) \coloneqq \| \bfu^k \|_{\rH^3(\Omega)} \| \bfd^k \|_{\rH^1(\Omega)} + \| \bfM^{n,k} \|_{\rL^2(\Omega)}$.
By the bounds for $\bfd^k$ and $\bfu^k$ in \eqref{recall:energy}, we conclude that $b_k\in\rL^2(0,T)$ and thus \eqref{eq:t-regs u & d}, finishing the proof.
\end{proof}

\subsection{Auxiliary results for the rigidity and artificial-viscosity limits}\label{ssec:aux ests for rigidity limit}

\begin{proof}[Proof of \autoref{lem:strong conv d}]
It is sufficient to consider $\sigma>0$ smaller than a fixed tubular radius of the rigid body and than one quarter of its distance from $\partial\Omega$.
The assertion for larger $\sigma$ then follows by restriction. 
By the uniform Hausdorff convergence \eqref{eq:conv of Hausdorff dist n} and the uniform continuity of the limiting rigid trajectory, we can cover $[0,T]$ by finitely many intervals $I_j$ and choose smooth fixed sets $A_j\Subset B_j\Subset\Omega$ such that, with $G_j:=\Omega\setminus\overline{B_j}$, for all sufficiently large $n$ and all $t\in I_j$,
\begin{equation}\label{eq:n nested localization sets}
 \cS(\varphi^n(t))\Subset A_j\subset\cS(\varphi(t))_{\sigma/2},
 \qquad
 \widetilde{\cF}_\sigma(t)\subset G_j,
 \qquad
 \dist(A_j,\Omega\setminus B_j)\ge c\sigma.
\end{equation}
For example, after choosing a center time $t_j$ so that the oscillation of $\cS(t)$ on $I_j$ is less than $\sigma/16$, one may take smooth tubular neighborhoods of $\cS(t_j)$ of radii $\sigma/4$ and $3\sigma/4$, with an arbitrarily small adjustment of the radii if needed.

Now we choose $\chi_j\in \rC^\infty(\overline\Omega)$ with $0\le\chi_j\le1$, $\chi_j=0$ on a neighborhood of $A_j$, and $\chi_j=1$ on a neighborhood of $\overline{G_j}$, and with derivatives bounded by constants depending only on $\sigma$ and  the fixed geometry. 
We then set $\tbfd^n=\bfd^n-\bfd_0$. 
Due to $\operatorname{Tr}_{\partial\Omega}\tbfd^n=0$, $\chi_j\tbfd^n$ belongs to $\rH^1_0(\Omega)$.
Moreover, on~$\operatorname{supp}\chi_j$ the penalization vanishes for all large $n$, and the molecular-field identity \eqref{eq:n-M-def} gives
\[
    \chi_j\Delta\bfd^n=\chi_j\bigl(\bfM^n+\bff_\eps(\bfd^n)\bigr),
\]
where we recall $\bff_\eps(\bfd)=-\eps^{-2}(1-|\bfd|^2)\bfd$.
The energy inequality and $\rH^1(\Omega)\hookrightarrow \rL^6(\Omega)$ imply, for fixed $\eps$, the bound $\sup_n\|\bff_\eps(\bfd^n)\|_{\rL^2(0,T;\rL^2(\Omega))}\le C_\eps$.
This shows the uniform boundedness of 
\[
    \Delta(\chi_j\tbfd^n)=\chi_j\bigl(\bfM^n+\bff_\eps(\bfd^n)-\Delta\bfd_0\bigr)+2\nabla\chi_j\cdot\nabla \tbfd^n+(\Delta\chi_j)\tbfd^n
\]
in $\rL^2(I_j;\rL^2(\Omega))$. 
An $\rH^2$-estimate on $\Omega$ and restricting to $G_j$ then yield $\sup_n\|\bfd^n\|_{\rL^2(I_j;\rH^2(G_j))} \le C_{\sigma,\eps}$.
By the choice of $B_j$ as smooth, $G_j$ is a fixed $\rC^2$-domain, so the compact embedding~$\rH^2(G_j)\Subset \rH^1(G_j)$ used below is standard.

For $\partial_t\bfd^n$, we invoke the distributional form of the director equation $\partial_t\bfd^n+(\bfu^n\cdot\nabla)\bfd^n=\bfM^n$.
Combining the boundedness of $\bfM^n$ in $\rL^2_{t,\bfx}$, $\bfu^n$ in $\rL^2_t\rH^3_{\bfx}\hookrightarrow \rL^2_t\rL^\infty_{\bfx}$, and $\nabla\bfd^n$ in $\rL^\infty_t\rL^2_{\bfx}$, we get $\sup_n\|\partial_t\bfd^n\|_{\rL^2(0,T;\rL^2(\Omega))}\le C_{\delta,\eps}$.
Applying the Aubin--Lions lemma on each fixed domain $G_j$, where we exploit $\rH^2(G_j)\Subset \rH^1(G_j)\hookrightarrow \rL^2(G_j)$, gives $\bfd^n\to\bfd$ strongly in $\rL^2(I_j;\rH^1(G_j))$.
By \eqref{eq:n nested localization sets}, the spacetime $(0,T)\times\widetilde\cF_\sigma(t)$ is contained in the union of these fixed cylinders. 
Summing over the finite covering proves $\bfd^n\to\bfd$ strongly in $\rL^2(0,T;\rH^1(\widetilde\cF_\sigma(t)))$, which is the assertion of \autoref{lem:strong conv d}.
\end{proof}

We next discuss the strong convergence of velocity fields $(\bfv_j)$.

\begin{lem}\label{lem:abstract velocity compactness}
Let $(\bfv_j)$ be bounded in $\rL^\infty(0,T;\rL^2(\Omega;\R^3))$ with $\bfv_j\rightharpoonup\bfv$ weakly in $\rL^2((0,T)\times\Omega;\R^3)$. 
Suppose that for every $\eta>0$ and $\gamma\in(0,T/2)$ there are a partition $\gamma=s_0<s_1<\dots<s_N=T-\gamma$ with intervals $I_\ell:=(s_{\ell-1},s_\ell)$, connected open sets $A_\ell\Subset\Omega$ with $\rC^2$-boundary, an index $j_\eta$, a constant~$C_\eta$, and fields $\bfz_{j,\ell}\in\rL^2(I_\ell;\cK(A_\ell))$, $1\le\ell\le N$, such that, for all $j\ge j_\eta$,
\begin{equation}\label{eq:abstract restriction}
    \sum_{\ell=1}^N\int_{I_\ell}\|\bfv_j-\bfz_{j,\ell}\|_{\rL^2(\Omega)}^2\srd t\le\eta, \qquad \sum_{\ell=1}^N\int_{I_\ell}\|\bfz_{j,\ell}\|_{\rH^1(\Omega)}^2\srd t\le C_\eta,
\end{equation}
and such that, for every $\ell$ and every $\bfe\in\cK(A_\ell)\cap\rC_c^\infty(\Omega;\R^3)$, there is a constant $C_{\eta,\ell,\bfe}$ with
\begin{equation}\label{eq:abstract temporal}
 \left|\int_{I_\ell}(\bfv_j(t),\bfe)_{\rL^2(\Omega)}\vartheta'(t)\srd t\right|
 \le C_{\eta,\ell,\bfe}\|\vartheta\|_{\rL^\infty(I_\ell)}
 \tfor \vartheta\in\rC_c^\infty(I_\ell) \tand j\ge j_\eta.
\end{equation}
Then $\bfv_j\to\bfv$ strongly in $\rL^2((0,T)\times\Omega;\R^3)$.
\end{lem}

\begin{proof}
We define $\overline V:=\sup_j \|\bfv_j\|_{\rL^\infty(0,T;\rL^2(\Omega))}<\infty$.
Next, we fix $\gamma\in(0,T/2)$ and $\eta>0$, and choose the partition, sets, and approximating fields from the hypotheses.

For each $\ell$, the space $\cK(A_\ell)\cap\rC_c^\infty(\Omega;\R^3)$ is dense in $\cK(A_\ell)$ in the $\rH^1$-norm.
Indeed, choose a function $\chi_\ell$ in $\rC_c^\infty(\Omega)$ that equals one near
$\overline{A_\ell}$. 
For an affine rigid field $\bfR(\bfx)=\bfa+\bfomega\times\bfx$, set
\[
    E_\ell\bfR:=\nabla\times\left[\chi_\ell\left(\frac12\bfa\times\bfx-\frac12|\bfx|^2\bfomega\right)\right].
\]
This field is smooth, solenoidal, compactly supported in $\Omega$, and equal to $\bfR$ near $\overline{A_\ell}$. For $\bfz\in\cK(A_\ell)$, we invoke that $A_\ell$ is connected to identify the  restriction there with a rigid field $\bfR_{\bfz}$. 
As a consequence, $\bfz-E_\ell\bfR_{\bfz}$ vanishes in $A_\ell$ and belongs to the solenoidal subspace of $\rH_0^1(\Omega\setminus\overline{A_\ell};\R^3)$.
It readily follows from a density argument on this $\rC^2$ open set that there are approximations in~$\rC_c^\infty(\Omega\setminus\overline{A_\ell};\R^3)$.
Extending them by zero and adding $E_\ell\bfR_{\bfz}$ proves the claimed density.

We then orthonormalize a countable subset that is dense in $\rH^1$ to find  an $\rL^2$-orthonormal basis $(\bfe_{\ell,m})_{m\ge1}$ of the $\rL^2$-closure of $\cK(A_\ell)$, with $\bfe_{\ell,m}\in\cK(A_\ell)\cap\rC_c^\infty(\Omega;\R^3)$.
Let $P_{\ell,M}$ be the orthogonal projection onto~$\Span\{\bfe_{\ell,1},\dots,\bfe_{\ell,M}\}$. 
Applying the Rellich--Kondrachov compactness theorem and by strong convergence of $P_{\ell,M}$, we have
\[   \alpha_{\ell,M}:=\sup_{\substack{\bfz\in\cK(A_\ell)\\\|\bfz\|_{\rH^1(\Omega)}\le1}}\|(I-P_{\ell,M})\bfz\|_{\rL^2(\Omega)}\longrightarrow0 \tfor M\to\infty.
\]

For every $\ell$ and $m$, the scalar functions
$t\mapsto(\bfv_j(t),\bfe_{\ell,m})_{\rL^2(\Omega)}$ are uniformly
bounded in $\rL^\infty(I_\ell)$ and, by
\eqref{eq:abstract temporal}, in $\mathrm{BV}(I_\ell)$ for $j\ge j_\eta$.
Compactness in $\mathrm{BV}(I_\ell)$, together with the uniform
$\rL^\infty$-bound, yields relative compactness in
$\rL^2(I_\ell)$. Consequently, for fixed $M$, the fields
\[
 \boldsymbol V_{j,M}(t):=
 \begin{cases}
  P_{\ell,M}\bfv_j(t),&t\in I_\ell,\quad 1\le\ell\le N,\\
  0,&t\in(0,\gamma)\cup(T-\gamma,T),
 \end{cases}
\]
form a relatively compact family in
$\rL^2((0,T)\times\Omega;\R^3)$.

Since $I-P_{\ell,M}$ is a contraction and
$\bfz_{j,\ell}(t)\in\cK(A_\ell)$, for almost every $t\in I_\ell$,
\[
 \|\bfv_j-P_{\ell,M}\bfv_j\|_{\rL^2(\Omega)}
 \le
 \|\bfv_j-\bfz_{j,\ell}\|_{\rL^2(\Omega)}
 +\alpha_{\ell,M}\|\bfz_{j,\ell}\|_{\rH^1(\Omega)}.
\]
Squaring, integrating, and using \eqref{eq:abstract restriction}
and the definition of $\overline V$, we obtain
\[
 \|\bfv_j-\boldsymbol V_{j,M}\|_
      {\rL^2((0,T)\times\Omega)}^2
 \le
 2\eta+2C_\eta\max_{1\le\ell\le N}\alpha_{\ell,M}^2
 +2\gamma\overline V^2 \tfor j\ge j_\eta.
\]
Thus, by first choosing $\gamma$ and $\eta$ sufficiently small
and then $M$ sufficiently large, the tail of $(\bfv_j)$ can be
approximated arbitrarily closely by a relatively compact family.
Adding the finitely many omitted fields shows that $(\bfv_j)$
is totally bounded, hence relatively compact, in
$\rL^2((0,T)\times\Omega;\R^3)$. Every strong cluster point equals
$\bfv$ by the assumed weak convergence. Therefore the whole
sequence converges strongly to $\bfv$.
\end{proof}

\begin{proof}[Proof of \autoref{lem:penalization velocity compactness}]
We combine approximate rigidity with temporal compactness in fixed spaces
of solenoidal fields. 
The argument follows the restriction-operator strategy of \cite[Theorem~3.2]{MST:02} and \cite{F3}. 
Throughout the proof, $C$ is a constant that is independent of $n$, the comparison-set index $j$, and sufficiently small $\sigma$, whereas constants denoted by $C_\sigma$ may depend on $\sigma$.
All pointwise-in-time estimates hold for almost every $t$.

We first show that the sets $S_n(t)$ have a uniformly controlled geometry and that $\bfv_n$ is close to a rigid velocity field on them.
By \eqref{eq:uniform bilipschitz penal compactness} and Morrey's embedding, the sets $S_n(t)$ enjoy a uniform $\rC^{1,\alpha}$-H\"older regularity for every $0<\alpha<1/2$. 
Let $\bfR_n(t,\bfx)=\bfa_n(t)+\bfomega_n(t)\times\bfx$ be the $\rL^2(S_n(t))$-orthogonal projection of $\bfv_n(t)$ onto the space of affine rigid velocities, and set~${\bfq_n:=\bfv_n-\bfR_n}$.
Since the Gram matrices associated with the rigid motion are uniformly coercive, we get the estimate $|\bfa_n(t)|+|\bfomega_n(t)|\le C\|\bfv_n(t)\|_{\rL^2(S_n(t))}$.
By Korn's inequality, with the rigid motion chosen by this projection, we obtain
\begin{equation}\label{eq:korn rigid projection penal compactness}
 \|\bfq_n(t)\|_{\rH^1(S_n(t))}\le C\|\bD(\bfv_n(t))\|_{\rL^2(S_n(t))}.
\end{equation}
To justify the choice of projection, one may first apply Korn's
inequality with any minimizing rigid motion and then use the
$\rL^2$-minimality of $\bfR_n$ and the uniform equivalence of norms
on the six-dimensional rigid-motion space.

On this uniformly Lipschitz family, the Sobolev extension operators
are bounded simultaneously at the $\rH^1$- and $\rH^3$-levels, with
uniform constants. 
By interpolation of the extensions, the embedding~$\rH^2(\R^3)\hookrightarrow\rL^\infty(\R^3)$, and \eqref{eq:korn rigid projection penal compactness}, we have
\[
    \|\bfq_n(t)\|_{\rL^\infty(S_n(t))} \le C\|\bfq_n(t)\|_{\rH^1(S_n(t))}^{1/2} \|\bfq_n(t)\|_{\rH^3(S_n(t))}^{1/2} \le C\|\bD(\bfv_n(t))\|_{\rL^2(S_n(t))}^{1/2}\|\bfv_n(t)\|_{\rH^3(\Omega)}^{1/2}.
\]
Here the rigid-coefficient bound controls
$\|\bfq_n(t)\|_{\rH^3(\Omega)}$ by $C\|\bfv_n(t)\|_{\rH^3(\Omega)}$. By the same bound,
together with Sobolev embedding on the fixed container, we have
\[
 \|\bfq_n(t)\|_{\rL^\infty(\Omega)}
 +[\bfq_n(t)]_{\rC^{0,1/2}(\overline\Omega)}\le C\|\bfv_n(t)\|_{\rH^3(\Omega)}.
\]
Applying the penalty estimate \eqref{eq:penalty est} and H\"older's inequality in time, together with the uniform~$\rL^2_t\rH^3_{\bfx}$-bound recorded in \eqref{eq:L2H3-bound}, we obtain
\begin{equation}\label{eq:integrated penalty penal compactness}
    \int_0^T\|\bD(\bfv_n(t))\|_{\rL^2(S_n(t))}^2\srd t\le Cn^{-1}, \enspace \int_0^T\|\bD(\bfv_n(t))\|_{\rL^2(S_n(t))}\|\bfv_n(t)\|_{\rH^3(\Omega)}\srd t \le Cn^{-1/2}.
\end{equation}
These estimates will control the integrated restriction error.

We now fix a sufficiently small $\sigma>0$ and take the intervals $I_j$ and comparison sets $S_j$ supplied by the hypotheses. 
Upon refining the covering, we may assume that the $I_j$ are consecutive intervals with disjoint interiors. Increasing $n_\sigma$ if necessary, we also find that $d_{\rm H}(\overline{S_n(t)},\overline{S(t)})\le\sigma$ uniformly in~$t$. 
Since the volume of the fluid is preserved, we have $|S_n(t)|=|S(t)|=|S_0|$.
On the other hand, we invoke that $S_n(t)\subset S_j\subset S(t)_\sigma$ and $(S_j)_\sigma\subset S(t)_{2\sigma}$ for $t\in I_j$, 
\[
 S_n(t)\subset S_j\subset S(t)_\sigma,
 \qquad (S_j)_\sigma\subset S(t)_{2\sigma}
 \qquad(t\in I_j),
\]
as well as the associated volume estimates
\begin{equation}\label{eq:thin shell penal compactness}
    |S_j\setminus S_n(t)|\le C\sigma,\qquad |(S_j)_\sigma\setminus S_n(t)|\le C\sigma.
\end{equation}
Moreover, every point of $\partial S_j$ is within distance $2\sigma$ of $\overline{S_n(t)}$. 
The preceding $\rL^\infty$ and H\"older estimates, and the uniform surface-area bound for $\partial S_j$,
thus yield
\begin{equation}\label{eq:boundary mismatch penal compactness}
 \begin{aligned}
 \|\bfq_n(t)\|_{\rL^2(\partial S_j)}
 &\le C\|\bD(\bfv_n(t))\|_{\rL^2(S_n(t))}^{1/2}\|\bfv_n(t)\|_{\rH^3(\Omega)}^{1/2}+C\sigma^{1/2}\|\bfv_n(t)\|_{\rH^3(\Omega)}.
 \end{aligned}
\end{equation}

Next, on each interval $I_j$, we replace $\bfv_n$ by a field that is exactly rigid on the fixed set $S_j$ and remains close to $\bfv_n$.
Choose a smooth cutoff $\eta_{j,\sigma}$ that is zero near $\overline{S_j}$, equals one outside $(S_j)_\sigma$, and satisfies $0\le\eta_{j,\sigma}\le1$ and $|\nabla\eta_{j,\sigma}|\le C\sigma^{-1}$. 
Set $\bfz^0_{n,j}:=\bfR_n+\eta_{j,\sigma}\bfq_n$ and $\mdiv\bfz^0_{n,j}=\nabla\eta_{j,\sigma}\cdot\bfq_n$.
Then $\bfz^0_{n,j}$ equals $\bfR_n$ on $S_j$ and $\bfv_n$
outside $(S_j)_\sigma$. Its divergence is compactly supported in
the fixed-width collar
$D_j=(S_j)_{r_*}\setminus\overline{S_j}$, provided $\sigma<r_*$.
The integral of $\mdiv\bfz^0_{n,j}$ vanishes on every connected component of
$D_j$. Indeed, the field equals $\bfR_n$ on the inner boundary and
$\bfv_n$ on the outer boundary, and both have zero flux through
each closed boundary surface. For $\bfv_n$, this follows by
extending the~$\rH_0^1(\Omega)$-field by zero to $\R^3$,
which preserves its divergence-free property.

On $D_j$, choose a Bogovski\u{\i} right inverse $\mathcal B_j$.
Besides the usual $\rL^2_0\to\rH_0^1$-bound, we use its
negative-order estimate, see \cite{GHH:06}, with $s=-1$ and $p=2$, namely
\begin{equation}\label{eq:neg order Bogovskii est}
 \|\mathcal B_j f\|_{\rL^2(D_j)}\le C|f|_{*,D_j},\qquad
 |f|_{*,D_j}:=
 \sup_{\substack{\phi\in\rH^1(D_j),\ \int_{D_j}\phi\srd\bfx=0\\
                  \|\nabla\phi\|_{\rL^2(D_j)}\le1}}
 \Big|\int_{D_j}f\phi\srd\bfx\Big|.
\end{equation}
This is the case $s=-1$, $p=2$ of \cite[Theorem~2.5]{GHH:06}. There, the source space is the dual of the full space $\rH^1(D_j)$, restricted to functionals that annihilate constants (and not $\rH^{-1}(D_j)=(\rH_0^1(D_j))'$), and Poincar\'e's inequality identifies its norm with $|\cdot|_{*,D_j}$.
If $D_j$ is disconnected, the operator and the zero-mean
normalization are understood componentwise. The constants in
both Bogovski\u{\i} estimates, and in the Poincar\'e and trace
inequalities below, are uniform, since the fixed-width collars considered here admit a uniformly controlled finite atlas (see also \cite{Bog:79} and \cite[Chap.~III.3]{Gal:94}).

Let $\bfn_j$ be the outward normal to $S_j$. Since
$\mdiv\bfq_n=0$ and $\eta_{j,\sigma}=1$ near the outer boundary
of~$D_j$, integration by parts gives
\[
 \int_{D_j}\mdiv\bfz^0_{n,j}\phi\srd\bfx
 =-\int_{D_j}(\eta_{j,\sigma}-1)\bfq_n\cdot\nabla\phi\srd\bfx
   +\int_{\partial S_j}(\bfq_n\cdot\bfn_j)\phi\srd\sigma.
\]
Consequently, retaining the inner-boundary term and applying the
trace and Poincar\'e inequalities, we obtain
\[
 |\mdiv\bfz^0_{n,j}|_{*,D_j}
 \le C\left(
 \|\bfq_n\|_{\rL^2((S_j)_\sigma\setminus\overline{S_j})}
 +\|\bfq_n\|_{\rL^2(\partial S_j)}\right).
\]
The first term is at most $C\sigma^{1/2}\|\bfv_n(t)\|_{\rH^3(\Omega)}$ by
\eqref{eq:thin shell penal compactness} and the global
$\rL^\infty$-bound for $\bfq_n$. Combining this with
\eqref{eq:boundary mismatch penal compactness} and
\eqref{eq:neg order Bogovskii est} yields
\begin{equation}\label{eq:bog correction penal compactness}
 \begin{aligned}
 \|\mathcal B_j\mdiv\bfz^0_{n,j}(t)\|_{\rL^2(D_j)}
 &\le C\|\bD(\bfv_n(t))\|_{\rL^2(S_n(t))}^{1/2}\|\bfv_n(t)\|_{\rH^3(\Omega)}^{1/2}+C\sigma^{1/2}\|\bfv_n(t)\|_{\rH^3(\Omega)}.
 \end{aligned}
\end{equation}
Extend $\mathcal B_j \mdiv\bfz^0_{n,j}$ by zero to $\Omega$ and define $\bfz_{n,j}:=\bfz^0_{n,j}-\mathcal B_j \mdiv\bfz^0_{n,j}$.
The zero trace of the correction ensures $\bfz_{n,j}(t)\in\cK(S_j)$.
The ordinary Bogovski\u{\i}
estimate joint with the estimate $\|\mdiv\bfz^0_{n,j}(t)\|_{\rL^2(D_j)}\le C\sigma^{-1}\|\bfv_n(t)\|_{\rH^3(\Omega)}$ gives
\begin{equation}\label{eq:restriction H1 penal compactness}
 \|\bfz_{n,j}(t)\|_{\rH^1(\Omega)}\le C_\sigma \|\bfv_n(t)\|_{\rH^3(\Omega)}.
\end{equation}
Furthermore, $\bfv_n-\bfz^0_{n,j}=(1-\eta_{j,\sigma})\bfq_n$ is supported in $(S_j)_\sigma$. 
On $S_n(t)$, its $\rL^2$-norm is controlled by~\eqref{eq:korn rigid projection penal compactness}. 
On the remaining set, \eqref{eq:thin shell penal compactness} and the global $\rL^\infty$-bound give $C\sigma^{1/2}\|\bfv_n(t)\|_{\rH^3(\Omega)}$. 
Together with the above bound \eqref{eq:bog correction penal compactness}, for $t\in I_j$ and $n\ge n_\sigma$, this proves
\begin{equation}\label{eq:pointwise restriction penal compactness}
    \begin{aligned}
        \|\bfv_n(t)-\bfz_{n,j}(t)\|_{\rL^2(\Omega)}
        &\le C\|\bD(\bfv_n(t))\|_{\rL^2(S_n(t))}+C\|\bD(\bfv_n(t))\|_{\rL^2(S_n(t))}^{1/2}\|\bfv_n(t)\|_{\rH^3(\Omega)}^{1/2}\\
        &\quad+C\sigma^{1/2}\|\bfv_n(t)\|_{\rH^3(\Omega)}.
    \end{aligned}
\end{equation}

It remains to verify the hypotheses of \autoref{lem:abstract velocity compactness} with $A_j:=S_j$ and the fields $\bfz_{n,j}$ constructed above.
Since the hypotheses of the present lemma hold on every subinterval, they may be restricted to $(\gamma,T-\gamma)$ for a given $\gamma\in(0,T/2)$. 
Squaring \eqref{eq:pointwise restriction penal compactness}, integrating over the intervals $I_j$, which have disjoint interiors, and using the integrated estimates \eqref{eq:integrated penalty penal compactness} together with \eqref{eq:L2H3-bound} and \eqref{eq:restriction H1 penal compactness}, we obtain, for $n\ge n_\sigma$, the bounds
\[
    \sum_j\int_{I_j}\|\bfv_n-\bfz_{n,j}\|_{\rL^2(\Omega)}^2\srd t \le C\bigl(n^{-1}+n^{-1/2}+\sigma\bigr) \tand \sum_j\int_{I_j}\|\bfz_{n,j}\|_{\rH^1(\Omega)}^2\srd t\le C_\sigma.
\]
Given $\eta>0$, we first choose $\sigma$ with $C\sigma\le\eta/2$ and then enlarge $n_\sigma$ such that $C(n^{-1}+n^{-1/2})\le\eta/2$ for $n\ge n_\sigma$. This gives \eqref{eq:abstract restriction}. Finally, for $\bfe\in\cK(S_j)\cap\rC_c^\infty(\Omega;\R^3)$ and $\vartheta\in\rC_c^\infty(I_j)$, the assumed temporal estimate applied to $\bfzeta(t,\bfx)=\vartheta(t)\bfe(\bfx)$ gives
\[
 \left|\int_{I_j}(\bfv_n(t),\bfe)_{\rL^2(\Omega)}\vartheta'(t)\srd t\right|
 \le C_{\sigma,j}\bigl(1+T^{1/2}\bigr)\|\bfe\|_{\rH^3(\Omega)}\|\vartheta\|_{\rL^\infty(I_j)},
\]
which is \eqref{eq:abstract temporal}. Hence \autoref{lem:abstract velocity compactness} yields $\bfv_n\to\bfv$ strongly in $\rL^2(0,T;\rL^2(\Omega))$.
\end{proof}

\begin{proof}[Proof of \autoref{lem:strong conv u}]
The argument follows the compactness strategy of~\cite[Theorem~3.2]{MST:02} and~\cite{F3}.
To estimate the temporal functional on $I_j$, we test \eqref{eq:penalized fluid} with functions compactly supported in time in $I_j$, which gives
\begin{align}
        \int_{I_j}\int_{\Omega} \bfu^n\cdot \del_t\bfzeta \srd \bfx \srd t
        &= \int_{I_j}\int_{\Omega} \Bigl[- (\bfu^n\otimes \bfu^n): \nabla\bfzeta + (2+n \varphi^n)\bD(\bfu^n):\bD(\bfzeta)\Bigr] \srd \bfx \srd t \notag\\
        &\quad + \delta\int_{I_j} a_u(\bfu^n,\bfzeta)\srd t
        +\int_{I_j}\int_{\Omega} (\bfzeta \cdot \nabla) \bfd^n \cdot \Delta \bfd^n \srd \bfx \srd t.\label{eq:pen fluid with Cc test fct}
\end{align}
The director-dependent term will be treated in stress form, so that only the energy-level $\rH^1$-bound is used.
We estimate the terms on the right-hand side of \eqref{eq:pen fluid with Cc test fct} for
\[
    \bfzeta\in \rC_c^1(I_j\times\Omega;\R^3)\cap \rL^2(I_j;\cK(A_j)\cap \rH^3(\Omega;\R^3)).
\]
Then $\varphi^n\bD(\bfzeta)=0$ for all sufficiently large $n$, because $\cS(\varphi^n(t))\Subset A_j$, so the penalized viscosity term drops out. A global $\rH^2$-bound for $\bfd^n$ is not available uniformly in $n$. The elastic term is therefore rewritten in stress form.
From $\mdiv\bfzeta=0$ and $\bfzeta=0$ on $\partial\Omega$, it follows that
\[
    -\int_\Omega(\bfzeta\cdot\nabla)\bfd^n\cdot\left(\Delta\bfd^n+\frac1{\eps^2}(1-|\bfd^n|^2)\bfd^n\right)\srd\bfx=\int_\Omega (\nabla\bfd^n)^\top\nabla\bfd^n:\nabla\bfzeta\srd\bfx.
\]
The contribution from the Ginzburg--Landau regularization is a spatial gradient and therefore vanishes. 
The energy bound allows us to infer that
\[
    \left|\int_{I_j}\int_\Omega(\nabla\bfd^n)^\top\nabla\bfd^n:\nabla\bfzeta\srd\bfx\srd t\right|\le C\|\bfzeta\|_{\rL^2(I_j;\rW^{1,\infty}(\Omega))}\le C\|\bfzeta\|_{\rL^2(I_j;\rH^3(\Omega))}.
\]
The convective term is bounded in the same norm by $\|\bfu^n\|_{\rL^\infty_t\rL^2_{\bfx}}^2$, the baseline viscous term by~$\|\nabla\bfu^n\|_{\rL^2(I_j\times\Omega)}$, and the penalization term vanishes because $\varphi^n\bD(\bfzeta)=0$. Finally,
\[
 \biggl|\delta\int_{I_j}a_u(\bfu^n,\bfzeta)\srd t\biggr|
 \le \sqrt\delta\,\|\sqrt\delta\,\bfu^n\|_{\rL^2(I_j;X_u)}\,\|\bfzeta\|_{\rL^2(I_j;X_u)}.
\]
Thus, for $0<\delta\le1$,
\begin{equation*}
    \biggl|\int_{I_j}\int_{\Omega}\bfu^n\,\del_t\bfzeta\srd\bfx\srd t\biggr|
    \le C(\bfu_0,\bfd_0,\eps,T)\,\|\bfzeta\|_{\rL^2(I_j;\rH^3(\Omega))},
\end{equation*}
with a constant independent of $n$ (and bounded for $0<\delta\le1$). The nested localization construction~\eqref{eq:n nested localization sets}, with the fixed comparison set $S_j:=A_j$ on each $I_j$, gives the fixed comparison sets required in \autoref{lem:penalization velocity compactness}.  Moreover, the artificial-viscosity estimate gives
$\bfu^n$ uniformly bounded in~$\rL^2(0,T;X_u)$ at fixed $\delta$.  Since
$X_u\hookrightarrow\rW^{1,\infty}(\Omega)$,
\[
 \sup_n\int_0^T\|\nabla\bfu^n(t)\|_{\rL^\infty(\Omega)}\srd t<\infty.
\]
Applying Gr\"onwall's inequality for the forward and inverse flow equations we get  uniform bi-Lipschitz bounds for
$\Phi_n(t):=\bfX^n_{0,t}$ and $\Phi_n(t)^{-1}=\bfY^n(t)$. The maps are volume preserving because $\mdiv\bfu^n=0$.  In addition,
$X_u\hookrightarrow\rH^3(\Omega)\hookrightarrow\rW^{2,6}(\Omega)$, so the differentiated-flow estimate established at the beginning of \autoref{ssec:limit n to infinity} gives
\[
 \sup_{n,t}\bigl(
 \|\Phi_n(t)\|_{\rW^{2,6}(\Omega)}+
 \|\Phi_n(t)^{-1}\|_{\rW^{2,6}(\Omega)}
 \bigr)<\infty.
\]
Thus, the transported sets $\cS(\varphi^n(t))=\Phi_n(t,\cS_0)$ satisfy the full geometric hypothesis \eqref{eq:uniform bilipschitz penal compactness}.  
The sets $A_j$ can be chosen as small tubular neighborhoods of a rigid copy of $\cS_0$. 
From the fixed positive distance of the limiting body from the wall, it follows that the comparison-set collar hypothesis of \autoref{lem:penalization velocity compactness} is also satisfied.
The energy inequality also yields the uniform bound
\[
n\int_0^T\int_\Omega\varphi^n|\bD(\bfu^n)|^2\srd\bfx\srd t\le C,
\]
so all hypotheses of \autoref{lem:penalization velocity compactness} are satisfied.
An application of \autoref{lem:penalization velocity compactness} then implies that $\bfu^n\to\bfu$ strongly in~$\rL^2(0,T;\rL^2(\Omega))$, as claimed.
\end{proof}

\begin{proof}[Proof of \autoref{lem:density director tests}]
Because the limiting solid stays a positive distance from $\partial\Omega$ on the interval under consideration and moves by a rigid motion, there exists a family of bi-Lipschitz diffeomorphisms $\Psi(t)\colon\Omega\to\Omega$, with $\Psi(t)(\cS_0)=\cS(t)$. This map equals the rigid map $\bfh(t)+\bfO(t)\,\cdot$ in a fixed neighborhood of $\cS_0$, equals the identity in a fixed neighborhood of $\partial\Omega$, and is uniformly bounded in $\rW^{1,\infty}(0,T;\rW^{1,\infty}(\Omega))$ together with its inverse. 
For the construction of such a map, we extend the instantaneous rigid velocity smoothly from a fixed collar of the moving solid to an Eulerian vector field that vanishes in a fixed collar of $\partial\Omega$, with a uniform spatial Lipschitz bound. 
As the separation is positive, this extension can be chosen uniformly in time.
The flow of this ALE field equals the prescribed rigid motion near the solid and the identity near the container wall. 
Since the rigid translational and angular velocities are essentially bounded in time, the flow and its inverse have the stated $\rW^{1,\infty}_t\rW^{1,\infty}_{\bfx}$-bounds.

For $\bfpsi\in\cX_{\bfd}(\varphi)$, we set $\widehat{\bfpsi}(t,\bfy)=\bfpsi(t,\Psi(t,\bfy))$ on $\cF_0$. 
By the Sobolev chain rule, the following identities hold a.e. $\nabla_{\bfy}\widehat{\bfpsi}=(\nabla_{\bfx}\bfpsi)\circ\Psi\,\nabla_{\bfy}\Psi$ and $\partial_t\widehat{\bfpsi}=(\partial_t\bfpsi)\circ\Psi+(\nabla_{\bfx}\bfpsi)\circ\Psi\,\partial_t\Psi$.
The corresponding inverse formulas are obtained by composing with $\Psi^{-1}$. 
Uniform bounds for $\Psi$, $\Psi^{-1}$ and their first derivatives, together with uniform upper and lower Jacobian bounds, show that pullback is an isomorphism, with equivalent norms, between $\cX_{\bfd}(\varphi)$ and $\rH^1(0,T;\rL^2(\cF_0))\cap\rL^2(0,T;\rH_0^1(\cF_0))$ with $\widehat{\bfpsi}(T)=0$.
By construction, $\Psi(t,\cdot)$ maps the boundary components of $\cF_0$ onto those of $\cF(t)$, this composition also preserves the variational zero-trace condition defining $\rH_0^1$.
Next, we observe that $\rC_c^\infty([0,T)\times\cF_0)$ is dense in this fixed cylinder space. 
Choose $\widehat{\bfpsi}_m$ converging to $\widehat{\bfpsi}$ in the indicated norm. 
Pushing forward gives functions supported a positive spatial distance from the moving boundary and vanishing on a terminal time interval. 
Note that these pushforwards only have Sobolev regularity in time because  the map $\Psi$ is only Lipschitz in time.  
To smooth them without losing the value at $t=0$, we start by using the standard bounded extension operator for~$\rH^1(0,T;\rL^2(\Omega))\cap\rL^2(0,T;\rH^1(\Omega))$ across $t=0$ and extend by zero across the terminal interval.  
Because the support at $t=0$ is compactly contained in $\cF_0$ and the rigid trajectory is continuous, the extension can be chosen on a short negative-time interval with the same positive spatial clearance.  
A spacetime mollification with radius smaller than that clearance then produces $\bfpsi_m\in\cD_{\cF}$ without changing the limit. 
The continuity of the trace map $\rH^1(0,T;\rL^2(\Omega))\hookrightarrow\rC([0,T];\rL^2(\Omega))$ implies that $\bfpsi_m(0)\to\bfpsi(0)$ in $\rL^2(\Omega)$.
\end{proof}

\begin{proof}[Proof of \autoref{lem:FSI compactness criterion}]
The task is essentially to verify the hypotheses of \autoref{lem:abstract velocity compactness}. 
For this purpose, we fix $\gamma\in(0,T/2)$ and $\sigma>0$ smaller than one quarter of the uniform distance to $\partial\Omega$. 
The limiting rigid trajectory is uniformly continuous, and the rigid motions converge uniformly, so there exist a partition $\gamma=s_0<s_1<\dots<s_N=T-\gamma$ and fixed smooth connected sets $A_\ell\Subset\Omega$ such that, with $I_\ell:=(s_{\ell-1},s_\ell)$, for all large $j$ and all $t\in I_\ell$, we have $S_j(t)\Subset A_\ell$, $S(t)_{2\sigma}\subset A_\ell\subset S(t)_{3\sigma}$ as well as $A_\ell\subset (S_j(t))_{C\sigma}$.
The sets can be chosen as tubular enlargements of the limiting solid at a time in~$I_\ell$. 
After refining the partition, $\partial A_\ell$ and $\partial S_j(t)$ are uniformly normal graphs over one another, with graph height between $c\sigma$ and $C\sigma$. 
All constants below are independent of $j$, $t$, and small $\sigma$.

For a.e.\ $t\in I_\ell$, let $\bfR_j(t)$ be the unique affine rigid velocity agreeing with $\bfv_j(t)$ on $S_j(t)$ and put~$\bfq_j:=\bfv_j-\bfR_j$.  Thus $\bfq_j=0$ in $S_j(t)$.  
On the other hand, the rigid coefficients of $\bfR_j$ are bounded by $\|\bfv_j\|_{\rL^2(S_j(t))}$.
This implies that
\begin{equation}\label{eq:est of q_j by v_j}
    \| \bfq_j \|_{\rH^1(\Omega)} \le C \| \bfv_j \|_{\rH^1(\Omega)}.
\end{equation}
Uniform tubular coordinates for the rigid copies of $\cS_0$, followed by the one-dimensional fundamental theorem of calculus along the normal fibers, give the thin-layer Poincar\'e and trace estimates
\begin{equation}\label{eq:Poincare & trace q_j}
 \|\bfq_j\|_{\rL^2((S_j(t))_{C\sigma}\setminus S_j(t))}
 \le C\sigma\|\nabla\bfq_j\|_{\rL^2((S_j(t))_{C\sigma})},
 \qquad
 \|\bfq_j\|_{\rL^2(\partial A_\ell)}
 \le C\sigma^{1/2}\|\nabla\bfq_j\|_{\rL^2((S_j(t))_{C\sigma})}.
\end{equation}
These estimates extend from smooth fields to $\rH^1$ by density.

Choose a cutoff $\eta_{\ell,\sigma}$ equal to zero on $A_\ell$, equal to one outside $(A_\ell)_\sigma$, and satisfying $|\nabla\eta_{\ell,\sigma}|\le C\sigma^{-1}$.  Set $\bfz^0_{j,\ell}:=\bfR_j+\eta_{\ell,\sigma}\bfq_j$.
Take a fixed-width smooth collar $D_\ell$ of $A_\ell$, independent of the transition width, and let $\mathcal B_\ell$ be its Bogovski\u{\i} right inverse.  The flux of $\mdiv\bfz^0_{j,\ell}$ is zero.  As in \eqref{eq:neg order Bogovskii est} in the proof of \autoref{lem:penalization velocity compactness}, the negative-order estimate
\begin{equation}\label{eq:Bogovskii est}
    \|\mathcal B_\ell \mdiv\bfz^0_{j,\ell}\|_{\rL^2(D_\ell)}
 \le C |\mdiv\bfz^0_{j,\ell}|_{*,D_\ell}
\end{equation}
and integration by parts give
\begin{equation}\label{eq:est of div z_j,l}
 |\mdiv\bfz^0_{j,\ell}|_{*,D_\ell}
 \le C\left(
 \|\bfq_j\|_{\rL^2((A_\ell)_\sigma\setminus A_\ell)}
 +\|\bfq_j\|_{\rL^2(\partial A_\ell)}
 \right)
 \le C\sigma^{1/2}\|\bfq_j\|_{\rH^1(\Omega)}.
\end{equation}
The usual $\rL^2_0\to\rH^1_0$ Bogovski\u{\i} estimate gives an $\rH^1$-bound for the correction for each fixed $\sigma$.  Extending the correction by zero and setting $\bfz_{j,\ell}:=\bfz^0_{j,\ell}-\mathcal B_\ell \mdiv\bfz^0_{j,\ell}$, we obtain
\[
 \bfz_{j,\ell}(t)\in \cK(A_\ell)
 :=\{\bfz\in\cV(\Omega):\bD(\bfz)=0\text{ on }A_\ell\}.
\]
Thus, invoking that $\bfv_j-\bfz_{j,\ell} = \bfq_j - \eta_{\ell,\sigma}\bfq_j + \mathcal B_\ell \mdiv\bfz^0_{j,\ell}$, and using the estimates \eqref{eq:est of q_j by v_j}--\eqref{eq:est of div z_j,l}, for a.e.\ $t\in I_\ell$, we find that
\[
    \|\bfv_j(t)-\bfz_{j,\ell}(t)\|_{\rL^2(\Omega)} \le C\sigma^{1/2}\|\bfv_j(t)\|_{\rH^1(\Omega)}.
\]
Squaring, integrating over $I_\ell$, and summing over $\ell$, we obtain the estimate
\[
    \sum_\ell\int_{I_\ell}\|\bfv_j-\bfz_{j,\ell}\|_{\rL^2(\Omega)}^2\srd t \le C\sigma\sup_i\|\bfv_i\|_{\rL^2(0,T;\rH^1(\Omega))}^2 .
\]
For fixed $\sigma$, the $\rL^2_0\to\rH^1_0$ Bogovski\u{\i} estimate also gives $\sum_\ell\int_{I_\ell}\|\bfz_{j,\ell}\|_{\rH^1(\Omega)}^2\srd t\le C_\sigma$, uniformly in~$j$. 
Thus, \eqref{eq:abstract restriction} holds with $\eta=C\sigma$. 
We also argue that every $\bfe\in\cK(A_\ell)\cap\rC_c^\infty(\Omega;\R^3)$ is rigid on $A_\ell\supset S(t)_{2\sigma}$ for $t\in I_\ell$, so that $\bfzeta(t,\bfx)=\vartheta(t)\bfe(\bfx)$ with $\vartheta\in\rC_c^\infty(I_\ell)$ is admissible in \eqref{eq:FSI temporal functional bound} with $2\sigma$ in place of $\sigma$, showing \eqref{eq:abstract temporal}. 
As $\sigma$ can be chosen arbitrarily small, \autoref{lem:abstract velocity compactness} yields the claim.
\end{proof}

\bibliography{lcfsi_main_final}
\bibliographystyle{plain}

\end{document}